\documentclass[11pt,preprint]{elsarticle}
\usepackage{geometry}
\usepackage{algorithm}
\usepackage{algpseudocode}
\usepackage{stmaryrd}
\usepackage[colorlinks]{hyperref}
\usepackage{graphicx}
\usepackage{epstopdf}
\usepackage{subfigure}
\usepackage{color}
\usepackage{amsmath}
\usepackage{amssymb,bm}
\usepackage{amsfonts}
\usepackage{amstext}
\usepackage{amsxtra}
\usepackage{mathrsfs}
\usepackage{amscd}
\usepackage{amsgen}
\usepackage{amsopn}
\usepackage{cases}
\usepackage{ntheorem}
\usepackage[utf8]{inputenc}
\usepackage{booktabs}
\usepackage{listings}
\usepackage[normalem]{ulem}
\numberwithin{equation}{section}
\allowdisplaybreaks[4]

\newcommand{\normdg}[1]{{\left\vert\kern-0.25ex\left\vert\kern-0.25ex\left\vert #1\right\vert\kern-0.25ex\right\vert\kern-0.25ex\right\vert}}

\newcommand{\topcaption}{%
\setlength{\abovecaptionskip}{0pt}%
\setlength{\belowcaptionskip}{10pt}%
\caption}

\newtheorem{theorem}{Theorem}[section]

\newtheorem{lemma}[theorem]{Lemma}

\newtheorem{remark}{Remark}[section]
\theorembodyfont{\normalfont}
\newproof{pf}{Proof}
\newproof{pot}{Proof of Theorem \ref{thm2}}

\journal{Elsevier}

\begin{document}

\begin{frontmatter}
\title{A discontinuous Galerkin method by patch reconstruction for convection-diffusion-reaction problems over polytopic meshes\tnoteref{mytitlenote}}
\tnotetext[mytitlenote]{This research was supported in part by National Science Foundation Grant DMS-11771348, Natural Science Foundation of Shaan Xi Province in 2019 (No.2019JQ--755) and Natural Science Foundation of Shaanxi Provincial Department of Education in 2019 (19JK0462).}
\author[add1]{Di Yang}\ead{yangdi0226@163.com}
\author[add1]{Yinnian He\corref{correspondingauthor}}\ead{heyn@mail.xjtu.edu.cn}
\cortext[correspondingauthor]{Corresponding author.}
\address[add1]{School of Mathematics and Statistics, Xi'an Jiaotong University, Xi'an, Shaanxi 710049, P.R. China }

\begin{abstract}
  In this article, using the weighted discrete least-squares, we propose a patch reconstruction finite element space with only one degree of freedom per element. As the approximation space, it is applied to the discontinuous Galerkin methods with the upwind scheme for the steady-state convection-diffusion-reaction problems over polytopic meshes. The optimal error estimates are provided in both diffusion-dominated and convection-dominated regimes. Furthermore, several numerical experiments are presented to verify the theoretical error estimates, and to well approximate boundary layers and/or internal layers.
\end{abstract}

\begin{keyword}
discontinuous Galerkin methods\sep patch reconstruction\sep polytopic meshes\sep
convection-dominated regime\sep optimal error estimates\sep boundary layers
\end{keyword}
\end{frontmatter}


\section{Introduction}
\label{sec:Intro}

In this work, we consider the following convection-diffusion-reaction problem:
\begin{align}
\label{p1}
-\nu\Delta u+\nabla\cdot(\bm{b}u)+c u&=f\qquad\text{in}~\Omega\subset\mathbb{R}^d,\\
\label{p2}
u&=g\qquad\text{on}~\Gamma=\partial\Omega,
\end{align}
where $d=2$ or $3$, $\Omega$ is a bounded open polygonal (when $d=2$) or polyhedral (when $d=3$) domain with boundary $\Gamma$, $\nu$ is a positive constant diffusivity coefficient, $\bm{b}$ is the velocity vector field defined on $\overline{\Omega}=\Omega\cup\Gamma$ with $\bm{b}\in[W^{1,\infty}(\Omega)]^d$, $c\in L^\infty(\Omega)$ is the reaction coefficient, $f$ is a source term, and $g$ is a function describing the
essential boundary conditions.

Though \eqref{p1}--\eqref{p2} looks simple, it is well known that the equation \eqref{p1} shows the hyperbolic characteristics when its P\'{e}clet number $Pe=|\bm{b}|/\nu\gg 1$. It results to the famous convection-dominated regime. If we let $|\bm{b}|\sim O(1)$, then the condition $Pe\gg1 $ is equivalent to $0<\nu\ll 1$. The above case originates from amounts of physical applications, such as chemical engineering, fluid mechanics, aerodynamics, etc., so it has been attracted many engineers and scientists from different fields over nearly half a century. Mathematically, in a convection-dominated regime, the solution and/or its derivatives change rapidly near the boundary lines, which is called boundary layer phenomenons. In addition, the thickness of boundary layers becomes thinner as the value of $\nu$ decreases. Because we require a numerical solution to fit the exact solution in all cases, the dominated convection needs more numerical innovations added to the standard finite element methods. If only the standard finite element methods are applied, terrible spurious oscillations will appear near boundary layers unless the computational mesh is the magnitude of boundary layers \cite{Lin20181482}.

In the last two decades, the discontinuous Galerkin (DG) finite element methods \cite{Arnold1982742, Arnold20021749, Babuska1973221, Babuska1973863, Bassi1997267, Baumann1999311, Cockburn19982440, Reed1973, Wheeler1978152}, firstly introduced in the early 1970's, have prompted vigorous developments. DG methods have prominent advantages of high-order accuracy, high parallelizability, local mass conservation, allowed hanging nodes and easy handling of complicated geometries due to discontinuous functions employed. However, one of the shortcomings is well-known more degrees of freedom than the standard finite element methods. Effectively reducing degrees of freedom, several variants of DG methods such as HDG \cite{Cockburn20093827, Nguyen20093232, Shin2015292} and weak Galerkin methods \cite{Chen2017107, Lin20181482} have been developed to solve \eqref{p1}--\eqref{p2}. We note that most numerical experiments in the above works are implemented over traditional triangulations or several constructed meshes.

Recently, there are quite a few works concerning the DG method on polytopic (polygonal or polyhedral) meshes \cite{Antonietti20131417, Cangiani2016699, Cangiani20142009, DiPietro2014461, Dong20201309, Larsson201265, Li2019268}, which provide more flexibility in implementation, in particular for domain with microstructures or problems with certain physical constraints. Among these beautiful works, the patch reconstruction finite element space adopting to general polytopic meshes, which is proposed by Li et al. \cite{Li2019268}, attracts our attentions because of its arbitrary-order accuracy and only one degree of freedom per element. They will help reducing the scales of the linear equation system. However, following the rules introduced in \cite{Li2019268} of constructing patch reconstruction finite element space, we shall use the weighted discrete least-squares instead of the conventional discrete least-squares.

The rest of this article is organized as follows. In Section \ref{sec:Const}, we focus on details of constructing a patch reconstruction finite element space. In Section \ref{sec:DG:appro}, using the approximation space proposed in Section \ref{sec:Const}, we present a symmetric interior penalty DG method with the upwind scheme to solve the convection-diffusion-reaction problem \eqref{p1}--\eqref{p2}, and provide the discrete stability. In Section \ref{prior:error}, we prove the theoretical results of error estimates in different norm regimes. In Section \ref{numer:ex}, four numerical experiments are presented to verify the results in Section \ref{prior:error} and to examine performance of approximating several popular and interesting phenomenons including boundary layers and internal layers. Finally, in Section \ref{conclusion}, we summarize our work, draw conclusions and point out some plans to improve.

\section{Reconstruction operator}
\label{sec:Const}

In this section, we shall propose a patch reconstruction finite element space with only one degree of freedom per element. It can adapt to general polytopic meshes.

\subsection{Preliminaries}

Throughout this paper we use the following standard function spaces. For a Lipschitz domain $D\subset\mathbb{R}^d$, $d\geqslant 1$, we denote by $W^{s,p}(D)$ the Sobolev space with indexes $s\geqslant 0$, $1\leqslant p\leqslant\infty$ of real-valued functions defined on $D$, endowed with the seminorm $|\cdot|_{W^{s,p}(D)}$ denoted by $|\cdot|_{s,p,D}$ and norm $\|\cdot\|_{W^{s,p}(D)}$ denoted by $\|\cdot\|_{s,p,D}$. When $p=2$, $H^s(D)$ is denoted as $W^{s,2}(D)$ and the corresponding seminorm and norm are written as $|\cdot|_{s,D}$ and $\|\cdot\|_{s,D}$, respectively. Furthermore, with $|D|$ we denote the $d$-dimensional Hausdorff measure of $D$.

As mentioned in Section \ref{sec:Intro}, $\Omega\subset\mathbb{R}^d~(d=2,3)$ is a bounded open polygonal (when $d=2$) or polyhedral (when $d=3$) domain.
The mesh $\mathcal T_h$ is a subdivision of $\Omega$ with disjoint open polygonal(polyhedral) elements $K$, which may not be convex or even star-shaped.
In the absence of hanging nodes(edges), we denote by $\mathcal E_h$ the set of $(d-1)$--dimensional edges(interfaces) of all element in $\mathcal T_h$, $\mathcal E_h^o$ the set of interior edges(faces), and $\mathcal E_h^b$ the set of boundary edges(faces). It is obvious that $\mathcal E_h=\mathcal E_h^o\cup\mathcal E_h^b$. We set
\[
h=\max_{K\in\mathcal T_h}h_K,
\quad h_K=\text{diam}(K),\quad\forall K\in\mathcal T_h.
\]
Let $P_k(D)$ be the space of all polynomials with total degree not greater than $k$ defined on domain $D$. Then we define two piecewise function spaces with respect to $\mathcal T_h$ as:
\begin{align*}
\mathcal V&=\big\{
v\in L^2(\Omega):v|_K\in H^s(K),~\forall K\in\mathcal T_h,~s\geqslant 3/2+\epsilon,\;\epsilon>0\big\},\\
V_h^k&=\big\{v\in L^2(\Omega):v|_K\in P_k(K),~\forall K\in\mathcal T_h\big\}.
\end{align*}

We assume that $\mathcal T_h$ satisfies the following shape regularity conditions, also shown in \cite{BdaVeiga2014, Brezzi2009277, Cangiani20142009, Li2019268}, in the sense of that: there exist
\begin{itemize}
  \item three positive numbers $\mathcal N_0$, $N$, $\gamma$ independent of the mesh size $h$;
  \item a compatible sub-decomposition $\widetilde{\mathcal T}_h$ consisting of shape-regular $d$-simplexes,
\end{itemize}
such that
\begin{enumerate}[({A}1)]
 \item any element $K\in\mathcal T_h$ has a finite number of edges(faces). Mathematically, there is a positive integer $\mathcal N_0$ such that
     \begin{equation*}
      \max_{K\in\mathcal T_h}\#\big\{e\in\mathcal E_h: e\subset\partial K\big\}\leqslant \mathcal N_0,
      \end{equation*}
      where the cardinality of $U$ is denoted by $\# U$ for any finite-number set $U$;
  \item any element $K\in\mathcal T_h$ admits a decomposition $\widetilde{\mathcal T}_h|_K$ that consists of at most $N$ shape-regular $d$-simplexes;
  \item for any $d$-simplex $\tau\in\widetilde{\mathcal T}_h$, the ratio $h_\tau/\rho_\tau$ is bounded by $\gamma$, where $h_\tau$ is the diameter of $\tau$ and $\rho_\tau$ is the radius of the largest ball inscribed in $\tau$.
\end{enumerate}

The above assumptions lead to the following useful properties, which are extensively used in the later analysis.
\begin{enumerate}[({M}1)]
  \item There exists a constant $\kappa_0\geqslant 1$ that depends only on $N$ and $\gamma$, such that $h_K\leqslant \kappa_0 h_\tau$ for any $K\in\mathcal T_h$ and any $\tau\in\widetilde{\mathcal T}_h|_K$.
  \item {[\emph{Trace inequality}]} $\forall v\in H^s(K)$, $s\geqslant 1$,  there exists a positive constant $\kappa_1$ that depends only on $N$ and $\gamma$, such that
      \begin{equation}\label{polygon:tr}
          \|v\|_{0,e}\leqslant\kappa_1 |e|^{\frac{1}{2}}|K|^{-\frac{1}{2}}
          (\|v\|_{0,K}+h_K|v|_{1,K})
          \quad\forall e\subset\partial K.
      \end{equation}
  \item {[\emph{Approximation property}]} $\forall v\in H^{s}(K)$, $s\geqslant1$, there exists a positive constant $\kappa_2$ that depends only on $N$, $k$, $\gamma$, and an approximation polynomial $\widetilde{v}\in P_k(K)$, such that
      \begin{equation}\label{polygon:appro}
      \|v-\widetilde{v}\|_{t,K}
      \leqslant \kappa_2 h_K^{\min\{k+1,s\}-t}|v|_{s,K}\quad
      \forall~0\leqslant t\leqslant s.
      \end{equation}
  \item {[\emph{Local inverse inequality}]} $\forall v\in P_k(K)$, there exists a positive constant $\kappa_3$ that depends only on $N$, $\gamma$ and $k$, such that
      \begin{equation}\label{polygon:inverse}
      |v|_{1,K}\leqslant\kappa_3 h_K^{-1}\|v\|_{0,K}.
      \end{equation}
\end{enumerate}

The properties (M1) and (M4) have been proved in \cite{Li2019268}, and (M2)--(M3) can be found in \cite{BdaVeiga2014, Riviere2008}. Moreover, based on (M2) and (M4), the following local inverse trace inequality can be easily derived.
\begin{enumerate}[(M5)]
    \item {[\emph{Local inverse trace inequality}]} $\forall v\in P_k(K)$, there exists a positive constant $\kappa_4$ that depends only on $N$, $\gamma$ and $k$, such that
        \begin{equation}\label{inverse:trace}
        \|v\|_{0,e}\leqslant\kappa_4 |e|^{\frac{1}{2}}|K|^{-\frac{1}{2}}
        \|v\|_{0,K}\quad\forall e\subset\partial K.
        \end{equation}
\end{enumerate}
Moreover, for any $K\in\mathcal T_h$ and any $\tau\in \widetilde{\mathcal T}_h|_K$, it follows from (A3) and (M1) that
\begin{equation*}
|K|\geqslant
\begin{cases}
\pi\rho_\tau^2\geqslant
\frac{\pi}{\gamma^2\kappa_0^2}h_K^2&d=2,\\
\frac{4}{3}\pi\rho_\tau^3\geqslant
\frac{4\pi}{3\gamma^3\kappa_0^3}h_K^3,&d=3,
\end{cases}
\end{equation*}
which implies that there exists a constant $C_{g}$ that depends only on $d$, $N$, and $\gamma$, such that
\begin{equation}\label{polygon:area}
|K|^{-1}\leqslant C_{g} h_K^{-d},\quad
\forall K\in\mathcal T_h.
\end{equation}
Collecting \eqref{inverse:trace}, \eqref{polygon:area}, and the fact
$|e|\leqslant h_K^{d-1}~(\forall e\subset\partial K)$, we obtain another version of (M5) as follows:
\begin{enumerate}[(M6)]
    \item {[\emph{Another version of local inverse trace inequality}]} $\forall v\in P_k(K)$, there exists a positive constant $\kappa_5$ that depends only on $d$, $N$, $\gamma$ and $k$, such that
        \begin{equation}\label{inverse:trace2}
        \|v\|_{0,e}\leqslant\kappa_5 h_K^{-\frac{1}{2}}\|v\|_{0,K}\quad\forall e\subset\partial K.
        \end{equation}
\end{enumerate}

\subsection{Element patch}
\label{subsec:ep}

If given a polytopic mesh $\mathcal T_h$, for each element $K\in\mathcal T_h$, using the Von Neumann neighbor (adjacent face-neighboring elements) technique \cite{Sullican200101}, we shall recursively construct an element patch $S(K)$ which is an agglomeration of elements that contain $K$ itself and some elements around $K$.

To fix ideas, we denote by $\bm{a}_K$ the barycenter coordinates of element $K\in\mathcal T_h$, and we define
\begin{equation}\label{re:patch}
\begin{split}
S_0(K)&:=\{K\},\\
S_j(K)&:=\big\{K'\in\mathcal T_h~|~
K'\notin S_{j-1}(K),~\text{and}\\
&\qquad\exists E\in S_{j-1}(K)~\text{such~that~}K'\cap E\in\mathcal E_h\big\},\quad j=1,2,\cdots.
\end{split}
\end{equation}
Then, we set a uniform threshold $M$ for the cardinality of $S(K)$ so that there always exists some integer $N_0\geqslant 1$ such that
\begin{equation}\label{form:patch}
\bigcup_{j=0}^{N_0-1}S_j(K)\subseteq S(K)\subsetneqq\bigcup_{j=0}^{N_0}S_j(K).
\end{equation}
Generally speaking, we enlarge $S(K)$ element by element, and stop the recursive procedure \eqref{re:patch} if the unique integer $N_0$ satisfies \eqref{form:patch}.

In order to well implement, following the similar technique shown in \cite{Li2020112902}, in Algorithm \ref{alg:element:patch} we provide the details of constructing $S(K)$ for each element $K\in\mathcal T_h$. Furthermore, one can refer to \cite{Li2019524, Li2019268, Li201901} and the references therein for more discussion.

\begin{algorithm}
\caption{Constructing Element Patch}
\label{alg:element:patch}
\begin{algorithmic}[1]
\Require
a polytopic mesh $\mathcal T_h$ and a uniform threshold $M$;
\Ensure
the element patch $S(K)$ for all $K\in\mathcal T_h$;
\For{each element $K\in\mathcal T_h$}
\State set $j=0$, $S_j(K)=\{K\}$, $\mathcal I(K)=\{\bm{a}_K\}$;
\While{$\# S_j(K)<M$}
\State initialize sets $S_{j+1}(K)=S_j(K)$;
\For{each element $K\in S_j(K)$}
\State add all adjacent face-neighbouring elements of $K$ to $S_{j+1}(K)$;
\State add the barycenter coordinates of all elements in $S_{j+1}(K)\setminus S_j(K)$ to $\mathcal I(K)$;
\EndFor
\State let $j=j+1$ and delete $S_j(K)$;
\EndWhile
\State sort the Euclidean distances between points in $\mathcal I(K)$ and $\bm{a}_K$;
\State select the $M$ smallest values and collect the corresponding elements to form $S(K)$;
\EndFor
\end{algorithmic}
\end{algorithm}



\subsection{Constructing reconstruction operator with one unknown per element}
\label{operator:construct}

To ease the presentation we shall restrict ourselves
to the two-dimensional case, although the results here presented also hold in three
dimensions.

For any polygonal element $K\in\mathcal T_h$, we denote by $E_0$, $E_1$, $\cdots$, $E_{M-1}$, total $M$ elements involved in the element patch $S(K)$, and denote by $\bm{a}_i=(x_i, y_i)$ the barycenter coordinates of $E_i$, $i=0,\cdots, M-1$. We assume that $E_0=K$, and $\bm{a}_0=\bm{a}_K$ for convenience. Then, we introduce the following scaled monomial defined on $S(K)$:
\begin{equation}\label{eq:scaled:monomial}
m_{\bm{\alpha}}(\bm{x})=\frac{(x-x_0)^{\alpha_0}(y-y_0)^{\alpha_1}}
{h_K^{|\bm{\alpha}|}}\quad\forall\,\bm{x}=(x,y)\in S(K),
\end{equation}
where the double-index, $\bm{\alpha}$, is a two-tuple of non-negative integers, $\alpha_0$, $\alpha_1$, and the length of $\bm{\alpha}$ is given by $|\bm{\alpha}|:=\alpha_0+\alpha_1$. We set $n_t=(t+1)(t+2)/2$ for any non-negative integer $t$. It is easily observed that the following set of $n_k$ monomials:
\[
\big\{m_{(0,0)}(\bm{x}),~m_{(1,0)}(\bm{x}),~m_{(0,1)}(\bm{x}),~\cdots,~m_{(k,0)}(\bm{x}),
~\cdots,~m_{(k-j,j)}(\bm{x}),~\cdots,~m_{(0,k)}(\bm{x})\big\},\quad\forall\,\bm{x}\in S(K),
\]
is a basis of $P_k(S(K))$. The above $n_k$ monomials can also form a vector-valued function $\bm{m}(\bm{x})$ as
\begin{equation}\label{basis:pk}
\bm{m}(\bm{x})=
\big(1,~m_{(1,0)}(\bm{x}),~ m_{(0,1)}(\bm{x}),~\cdots,~m_{(k,0)}(\bm{x}),~
\cdots,~m_{(0,k)}(\bm{x})\big).
\end{equation}
Hence, we can define a matrix $\bm{X}_k$ composed of the barycenter point set $\{\bm{a}_i\}_{i=1}^{M-1}$ as
\[
\bm{X}_k=\left(\bm{m}(\bm{a}_1)^\top,\bm{m}(\bm{a}_2)^\top,\cdots,\bm{m}(\bm{a}_{M-1})^\top\right)^\top.
\]
In addition, for any $v\in V_h^0$, we set $v_i=v(\bm{a}_i)$, $i=0,1,\cdots,M-1$, and define a vector $\bm{v}$ as $\bm{v}=(v_1,v_2,\cdots,v_{M-1})^\top$.

Now for any given $v\in V_h^0$, we consider the following weighted discrete least-square problem with constraints:
\begin{equation}\label{discrete:wls}
\underset{\bm{\beta}\in\mathbb{R}^{n_k}}{\arg\min}~
\|\bm{W}(\bm{v}-\bm{X}_k\bm{\beta})\|^2,\quad
\text{s.t.}\quad \beta_0=v_0,
\end{equation}
where $\bm{\beta}=(\beta_0,\beta_1,\cdots,\beta_{n_k-1})^\top$ is the unknown coefficient vector, and the weighted matrix $\bm{W}$ is given by
\[
\bm{W}=\text{diag}\big\{\sqrt{w_1},\cdots,\sqrt{w_{M-1}}\big\}.
\]
Based on Algorithm \ref{alg:element:patch} and the weight function expressions from \cite{Liu2005, Zhang2017185}, we define the weight $w_j$ as:
\begin{equation*}
w_j=\frac{1}
{\displaystyle|\bm{a}_0-\bm{a}_j|^2
\cdot\sum_{j=1}^{M-1}\frac{1}{|\bm{a}_0-\bm{a}_j|^2}},
\quad j=1,2,\cdots,M-1,
\end{equation*}
where $|\bm{a}_0-\bm{a}_j|$ is the Euclidean distance between $\bm{a}_0$ and $\bm{a}_j$, $j=1,2,\cdots,M-1$. Then by substituting the constraint $\beta_0=v_0$ into the objective function, we can rewrite the primal problem \eqref{discrete:wls} as the following unconstrained optimization problem:
\begin{equation}\label{discrete:wls:2}
\underset{\underline{\bm{\beta}}\in\mathbb{R}^{n_k-1}}{\arg\min}~
\|\bm{W}(\underline{\bm{v}}-\underline{\bm{X}}_k\underline{\bm{\beta}})\|^2,
\end{equation}
where $\underline{\bm{v}}=(v_1-v_0,v_2-v_0,\cdots,v_{M-1}-v_0)^\top$, $\underline{\bm{\beta}}=(\beta_1,\beta_2,\cdots,\beta_{n_k-1})^\top$, and
\[
\underline{\bm{X}}_k=
\left[
  \begin{array}{cccccc}
    m_{(1,0)}(\bm{a}_1) & m_{(0,1)}(\bm{a}_1) & \cdots & m_{(k,0)}(\bm{a}_1) & \cdots & m_{(0,k)}(\bm{a}_1) \\
    m_{(1,0)}(\bm{a}_2) & m_{(0,1)}(\bm{a}_2) & \cdots & m_{(k,0)}(\bm{a}_2) & \cdots & m_{(0,k)}(\bm{a}_2) \\
    \vdots & \vdots & \ddots & \vdots & \ddots & \vdots \\
    m_{(1,0)}(\bm{a}_{M-1}) & m_{(0,1)}(\bm{a}_{M-1}) & \cdots & m_{(k,0)}(\bm{a}_{M-1}) & \cdots & m_{(0,k)}(\bm{a}_{M-1}) \\
  \end{array}
\right].
\]
It is well known that \eqref{discrete:wls:2} has a unique solution $\widehat{\underline{\bm{\beta}}}$ with the following explicit expression:
\begin{equation*}
\widehat{\underline{\bm{\beta}}}
=[(\bm{W}\underline{\bm{X}}_k)^\top
(\bm{W}\underline{\bm{X}}_k)]^{-1}
(\bm{W}\underline{\bm{X}}_k)^\top
\bm{W}\underline{\bm{v}}.
\end{equation*}
In addition, the matrix $[(\bm{W}\underline{\bm{X}}_k)^\top
(\bm{W}\underline{\bm{X}}_k)]^{-1}
(\bm{W}\underline{\bm{X}}_k)^\top
\bm{W}$ is denoted as $\bm{G}$ for convenience. Finally, by adding $\beta_0$ to $\widehat{\underline{\bm{\beta}}}$, we obtain the unique solution $\widehat{\bm{\beta}}$ to \eqref{discrete:wls} with
\begin{equation}\label{solution:wls}
\widehat{\bm{\beta}}
=\left[
   \begin{array}{cc}
     1 & \bm{0} \\
     -\bm{G}\bm{I} & \bm{G} \\
   \end{array}
 \right]\left[
          \begin{array}{c}
            v_0 \\
            \bm{v} \\
          \end{array}
        \right],
\end{equation}
where $\bm{0}=(0,0,\cdots,0)$ with $\dim(\bm{0})=M-1$, and
$\bm{I}=(1,1,\cdots,1)^\top$ with $\dim(\bm{I})=M-1$. As a result, for any $K\in\mathcal T_h$, we can define a local $k$-th reconstruction operator $\mathcal P_K^k:V_h^0|_{S(K)}\rightarrow P_k(S(K))$ as
\begin{equation}\label{re:poly}
\mathcal P_K^kv(\bm{x})=\bm{m}(\bm{x})\widehat{\bm{\beta}},
\quad\forall v\in V_h^0|_{S(K)},~\forall\bm{x}\in S(K).
\end{equation}

In order to know how the operator $\mathcal P_K^k$ works, we denote the vector-valued function $\bm{m}(\bm{x})$ given in \eqref{basis:pk} by another form $(1,m_1(\bm{x}),m_2(\bm{x}),\cdots,m_{n_k-1}(\bm{x}))$, and denote each matrix element of $\bm{G}$ by $g_{i,j}$, $i=1,\cdots,n_k-1$, $j=1,\cdots,M-1$. Then, by collecting the above notations, \eqref{solution:wls}, and \eqref{re:poly}, we rewrite $\mathcal P_K^k$ as:
\begin{align}
&\mathcal P_K^kv(\bm{x})\notag\\
&=\left[
    \begin{array}{cccc}
      1, & m_1(\bm{x}), & \cdots,& m_{n_k-1}(\bm{x}) \\
    \end{array}
  \right]
\left[
\begin{array}{cccc}
1       &    0    & \cdots &    0    \\
\displaystyle-\sum_{j=1}^{M-1}g_{1,j} & g_{1,1} & \cdots & g_{1,M-1}\\
\vdots  & \vdots  & \ddots &  \vdots  \\
\displaystyle-\sum_{j=1}^{M-1}g_{n_k-1,j} & g_{n_k-1,1} & \cdots & g_{n_k-1,M-1}
\end{array}\right]
\left[
\begin{array}{c}
v_0\\
v_1\\
\vdots\\
v_{M-1}
\end{array}\right]\notag\\
&=\left[
    \begin{array}{cccc}
      \displaystyle 1-\sum_{i=1}^{n_k-1}\sum_{j=1}^{M-1}
       m_i(\bm{x})g_{i,j},& \displaystyle \sum_{i=1}^{n_k-1}m_i(\bm{x})g_{i,1}, & \cdots, &\displaystyle
       \sum_{i=1}^{n_k-1}m_i(\bm{x})g_{i,M-1} \\
    \end{array}
  \right]
\left[
\begin{array}{c}
v_0\\
v_1\\
\vdots\\
v_{M-1}
\end{array}\right]\notag\\
&=\sum_{j=0}^{M-1}v(\bm{a}_j)\lambda_{j}(\bm{x}),
\quad\forall\bm{x}\in S(K),\label{re:poly:2}
\end{align}
where
\[
\lambda_0(\bm{x})=1-\sum_{i=1}^{n_k-1}\sum_{j=1}^{M-1}
       m_i(\bm{x})g_{i,j},\qquad
\lambda_j(\bm{x})=\sum_{i=1}^{n_k-1}m_i(\bm{x})g_{i,j},\quad
j=1,\cdots,M-1.
\]
Let $\mathcal H=\mathcal V\cup C^0(\Omega)$, and therefore, we can define a global $k$-th reconstruction operator $\mathcal P^k:\mathcal H\rightarrow V_h^k$ as
\begin{equation}\label{global:op}
(\mathcal P^kv)|_K=(\mathcal P_K^k\widetilde{v})|_K
\quad\text{with}\quad\widetilde{v}\in V_h^0,\quad\widetilde{v}(\bm{a}_K)=v(\bm{a}_K),\quad
\forall v\in \mathcal H.
\end{equation}
Actually, it follows from \eqref{global:op} that the image space $\mathcal P^k \mathcal H$ is equivalent to $\mathcal P^k V_h^0$. which is denoted as $S_h^k$  hereafter.

Furthermore, we shall outline a basis of $S_h^k$ to present more details. For any $K\in\mathcal T_h$, we introduce the characteristic function $\chi_K\in C^0(\Omega)$ satisfying
\begin{equation*}
\chi_K(\bm{x})=
\begin{cases}
1,&\bm{x}=\bm{a}_K,\\
0,&\bm{x}\in K'\in\mathcal T_h,~K'\neq K.
\end{cases}
\end{equation*}
We stress that in $K$ it is not necessary to know the value of $\chi_K(\bm{x})$ at $\bm{x}\neq\bm{a}_K$.
Then we define $\phi_K=\mathcal P^k\chi_K$, and it follows from \eqref{re:poly:2} and Lemma 3 in \cite{Li2019arxiv} that the set $\{\phi_K\}_{K\in\mathcal T_h}$ is a basis of $S_h^k$. Note that the support of $\phi_K$ is related with all the element patches which includes $K$:
\[
\text{supp}(\phi_K)=\bigcup_{E\in\mathcal T_h,K\in S(E)}E,\quad \forall\,K\in\mathcal T_h.
\]
Hence, for any $v\in \mathcal H$, we obtain the corresponding $k$-th polynomial $\mathcal P^kv$ described as
\begin{equation}\label{global:exp}
\mathcal P^k v(\bm{x})=\sum_{K\in\mathcal T_h} v(\bm{a}_K)\phi_K(\bm{x}),\quad \forall\bm{x}\in\Omega.
\end{equation}
As \eqref{global:exp} shows, $\mathcal P^k$ is also called a reconstruction operator with one unknown per element.

\begin{remark}\label{why:dg}
We note that $\mathcal P^k v(\bm{x})|_K$ is a piecewise $k$-th polynomial for any $K\in\mathcal T_h$. More specifically, $\phi_K$ is discontinuous across the interior edges/faces $e\in\mathcal E_h^o$, which whets our appetites of employing DG methods for solving the primal problem \eqref{p1}--\eqref{p2}.
\end{remark}

\subsection{Some properties and estimates of the reconstruction operator}
\label{useful:pk}

Recalling the optimization problem \eqref{discrete:wls:2}, we can find that \eqref{discrete:wls:2} is a matrix version of the following weighted discrete least-square problem with constraints: to find a $k$-th polynomial $\mathcal P_K^k v$ defined on $S(K)$ for any $v\in V_h^0$, such that
\begin{equation}\label{discrete:ls}
\mathcal P_K^k v=\underset{p\in P_k(S(K))}{\arg\min}
\sum_{j=1}^{M-1} w_j|v(\bm{a}_j)-p(\bm{a}_j)|^2,\quad\text{s.t.}\quad
p(\bm{a}_0)=v(\bm{a}_0).
\end{equation}
To ensure the unisolvence of \eqref{discrete:ls}, we need the following assumption:
\begin{enumerate}[(A4)]
  \item For any $K\in\mathcal T_h$ and any $p\in P_k(S(K))$, one has that
  \[
  p|_{\mathcal I_K}=0\quad\text{implies}\quad p|_{S(K)}\equiv 0,
  \]
  where $\mathcal I_K$ is the barycenter set $\{\bm{a}_i\}_{i=0}^{M-1}$.
\end{enumerate}
This assumption can be also found in \cite{Li2019268, Li2012259, Li2020112902}, where without proof the authors demanded that
the number $M$ should be greater than $n_k=\dim(P_k)$. In \ref{app:proof:A4}, we explicate the demand in the view of matrix analysis.

Further it stems from \eqref{discrete:ls} that
\begin{equation}\label{poly:projection}
(\mathcal P^k v)|_K=v|_K,\quad\forall v\in P_k(S(K)),\ \forall\,K\in\mathcal T_h,
\end{equation}
which implies that $\mathcal P^k$ can also be regarded as a $L^2$-projection operator on $V_h^k$. Hence, based on the above projection \eqref{poly:projection}, some related estimates of $\mathcal P^k$ from \cite{Li2019524, Li2019268, Li201901} can be employed here, which are presented as follows. In \ref{append:pf:24} we provide proofs of the follwing lemmas in this part.

\begin{lemma}\label{estimates:op}
Let $v\in H^{k+1}(\Omega)$ and $K\in\mathcal T_h$. If (A4) holds, then there exist two constants $\kappa_6$, $\kappa_7$ that both depend only on $N$, $\gamma$, and $k$, such that
\begin{align}
\|v-\mathcal P^kv\|_{0,K}
&\leqslant \kappa_6 h^{k+1}|v|_{k+1,S(K)}.
\label{estimate:L2op}\\
|v-\mathcal P^kv|_{1,K}
&\leqslant \kappa_7 h^k|v|_{k+1,S(K)}.
\label{estimate:H1op}
\end{align}
\end{lemma}

In order to get an appropriate trace estimate, we should first rewrite the trace inequality \eqref{polygon:tr}. Based on \eqref{polygon:tr}, \eqref{polygon:area}, and the fact $|e|\leqslant h_K^{d-1}$, we have
\begin{equation}\label{re:polygontr}
\|v\|_{0,e}\leqslant \widetilde{\kappa}_1|e|^{\frac{1}{2(d-1)}}\big(h_K^{-1}\|v\|_{0,K}+|v|_{1,K}
\big),
\quad\forall e\subset\partial K,
\end{equation}
where $\widetilde{\kappa}_1$ is a constant that depends only on $N$, $\gamma$, and $d$. Hence, the following lemma shows the appropriate trace estimate we need.

\begin{lemma}\label{estimates:op2}
Let $v\in H^{k+1}(\Omega)$ and $K\in\mathcal T_h$. If (A4) holds, then there exists a constant $\kappa_8$ that depends only on $N$, $\gamma$, $k$, and $d$, such that
\begin{equation}\label{estimate:trace}
\|v-\mathcal P^kv\|_{0,e}
\leqslant \kappa_8 |e|^{\frac{1}{2(d-1)}}h^k|v|_{k+1,S(K)},
\quad\forall e\subset\partial K.
\end{equation}
\end{lemma}


\section{DG with the reconstruction operator for the model problem}
\label{sec:DG:appro}

\subsection{Setting of the problem}
We first define the inflow and outflow parts of $\Gamma$ as follows:
\begin{align}
\Gamma^{-}&:=\big\{\bm{x}\in\Gamma:
~\bm{b}(\bm{x})\cdot\bm{n}(\bm{x})<0\big\}
=\text{inflow},\label{inflow:boundary}\\
\Gamma^{+}&:=\big\{\bm{x}\in\Gamma:
~\bm{b}(\bm{x})\cdot\bm{n}(\bm{x})\geqslant0\big\}
=\text{outflow},\label{outflow:boundary}
\end{align}
where $\bm{n}(\bm{x})$ is the unit outward normal vector to $\Gamma$ at $\bm{x}\in\Gamma$. We also define an ``effective" reaction function $r(\bm{x})$, and assume there exists a constant $r_0$ such that
\begin{equation}\label{asm:r}
r(\bm{x})=c(\bm{x})+\frac{1}{2}\nabla\cdot
\bm{b}(\bm{x})\geqslant
r_0\geqslant 0,\quad\forall\bm{x}\in\Omega.
\end{equation}

The well-posedness of the model problem \eqref{p1}--\eqref{p2} has been proved in \cite{Houston200199} under homogeneous boundary conditions. Further if the domain $\Omega$ is convex, $f\in H^s(\Omega),~s\geqslant 0$, $g=0$, and \eqref{asm:r} also holds, then \eqref{p1}--\eqref{p2} has the unique solution $u$ satisfies $u\in H^2(\Omega)\cap
H^{2+s}(\Omega)$, and the following priori estimates \cite{Roos2008}:
\begin{align*}
&\nu^{3/2}\|u\|_{2,\Omega}+\nu^{1/2}\|u\|_{1,\Omega}
+\|u\|_{0,\Omega}\leqslant C(\bm{b},r,d,\Omega)\|f\|_{0,\Omega},\\
&\nu^2\|u\|_{2+s,\Omega}
\leqslant C(d,\Omega)
(\nu+\|\bm{b}\|_{1,\infty,\Omega}+
\|r\|_{0,\infty,\Omega})\|f\|_{s,\Omega},\quad
s>0.
\end{align*}

Let $L$ be the diameter of $\Omega$. In order to keep homogeneity of dimensions, we define
\begin{equation}\label{Linfty:norm}
\|v\|_{k,\infty,\Omega}:=
\sum_{s=0}^k L^s|v|_{s,\infty,\Omega}
\quad \forall v\in W^{k,\infty}(\Omega).
\end{equation}
For the convection velocity $\bm{b}$, we shall provide the following assumptions \cite{Ayuso20091391}:
\begin{enumerate}[(H1)]
  \item there exists $\xi\in W^{1,\infty}(\Omega)$ such that
      \[
      \bm{b}\cdot\nabla\xi\geqslant
      2b_0:=2\frac{
      \|\bm{b}\|_{0,\infty,\Omega}}{L};
      \]
  \item there exists $C_b>0$ such that
      \[
      |\bm{b}(\bm{x})|\geqslant
      C_b
      \|\bm{b}\|_{1,\infty,\Omega},
      \quad\forall \bm{x}\in\Omega;
      \]
  \item there exists $C_r>0$ such that
      \[
      \|r\|_{0,\infty,K}\leqslant
      C_r\big(\min_{\bm{x}\in K}r(\bm{x})+b_0\big),
      \quad
      \forall K\in\mathcal T_h.
      \]
\end{enumerate}
In addition, it follows from \eqref{Linfty:norm} and (H1)--(H2) that
\begin{equation}\label{velocity:assumption}
C_b\frac{\|\bm{b}\|_{1,\infty,\Omega}}{L}
\leqslant b_0\leqslant\frac{\|\bm{b}\|_{1,\infty,\Omega}}{L},
\quad
|\bm{b}|_{1,\infty,\Omega}\leqslant
\frac{\|\bm{b}\|_{1,\infty,\Omega}}{L}
\leqslant\frac{b_0}{C_b}.
\end{equation}

Let $K^1$ and $K^2$ be two adjacent elements of the partition $\mathcal T_h$. We write $\bm{n}^1$ and $\bm{n}^2$ to denote the outward unit normal vectors on $e=\partial K^1\cap\partial K^2$, relative to $\partial K^1$ and $\partial K^2$, respectively. Let $v$ and $\bm{q}$ be a scalar-valued function and a vector-valued function, respectively. Both of them are smooth inside $K^1$ and $K^2$. Then, we write $v^i$ and $\bm{q}^i$ to denote the restrictions of $v$ and $\bm{q}$ to $K^i$, $i=1,2$, respectively. The averages of $v$ and $\bm{q}$ on $e$ are given by
\begin{equation}\label{average:dg}
\{\!\{v\}\!\}=\frac{1}{2}(v^1+v^2),\quad
\{\!\{\bm{q}\}\!\}=\frac{1}{2}(\bm{q}^1+\bm{q}^2),
\quad \text{on}~e\in\mathcal E_h^o.
\end{equation}
The jumps of $v$ and $\bm{q}$ across $e$ are given by
\begin{equation}\label{jump:dg}
\llbracket v \rrbracket =v^1\bm{n}^1+v^2\bm{n}^2,\quad
\llbracket \bm{q} \rrbracket =\bm{q}^1\cdot\bm{n}^1+\bm{q}^2\cdot\bm{n}^2,
\quad \text{on}~e\in\mathcal E_h^o.
\end{equation}
On a boundary edge/face $e\in\mathcal E_h^b\subset\partial\Omega$, we set
\begin{equation}\label{boundeage:dg}
\{\!\{v\}\!\}=v,\quad
\{\!\{\bm{q}\}\!\}=\bm{q},\quad
\llbracket v \rrbracket=v\bm{n},\quad
\llbracket \bm{q} \rrbracket=\bm{q}\cdot\bm{n} \quad
\text{on}~e\in\mathcal E_h^b.
\end{equation}

For future purposes we shall replace the average $\{\!\{\bm{b}v\}\!\}$ on each interior edge/face $e\in\mathcal E_h^o$ by the upwind value of $\bm{b}v$, which is defined as follows \cite{Brezzi20063293}:
\begin{equation}\label{upwinding}
\{\!\{\bm{b}v\}\!\}_{up}=
\begin{cases}
\bm{b}v^i,&\text{if}~\bm{b}\cdot\bm{n}^i>0,\\
\bm{b}v^j,&\text{if}~\bm{b}\cdot\bm{n}^i<0,\\
\bm{b}\{\!\{v\}\!\},&\text{if}~\bm{b}\cdot\bm{n}^i=0,\\
\end{cases}
\end{equation}
where $i,j\in\{1,2\}$ and $i\neq j$. However, on each boundary edge/face $e\in\mathcal E_h^b$, we already have $\{\!\{\bm{b}v\}\!\}=\bm{b}v$, and we leave it unchanged. Further one can find that the definition \eqref{upwinding} misses conciseness and convenience for our analysis and implementations. To solve it, we have the following formula via some simple derivations.
\begin{equation}\label{upwind:value}
\{\!\{\bm{b}v\}\!\}_{up}\cdot\bm{n}^i
=\Big(\{\!\{\bm{b}v\}\!\}+
\frac{|\bm{b}\cdot\bm{n}^i|}{2}
\llbracket v \rrbracket\Big)\cdot\bm{n}^i,\quad i=1,2.
\end{equation}

The following crucial formula \cite{Arnold20021749} will be extensively employed:
\begin{equation}\label{magic:formula}
\sum_{K\in\mathcal T_h}\int_{\partial K}
v\bm{q}\cdot\bm{n}ds
=\sum_{e\in\mathcal E_h}\int_e
\{\!\{\bm{q}\}\!\}\cdot\llbracket v \rrbracket ds+\sum_{e\in\mathcal E_h^o}\int_e\llbracket \bm{q} \rrbracket
\{\!\{v\}\!\} ds,
\end{equation}
where both $\bm{q}$ and $v$ are piecewise function over $\mathcal T_h$. We shall also  extensively use the following Poincar\'{e}--Friedrichs inequality for piecewise $H^1$ functions \cite{Brenner2003306, Riviere2008}:
\begin{equation}\label{Poincare:dg}
\|v\|_{0,\Omega}\leqslant
LC_P\bigg(|v|_{1,h}^2+\sum_{e\in\mathcal E_h}\frac{1}{|e|^{1/(d-1)}}
\|\llbracket v \rrbracket\|_{0,e}^2\bigg)^{\frac{1}{2}},
\end{equation}
where $C_P$ is a positive constant depending on the minimum angle of sub-decomposition $\widetilde{\mathcal T}_h$, and $|\cdot|_{1,h}$ is the broken $H^1$-seminorm, which is defined as:
\[
|v|_{1,h}=\bigg(\sum_{K\in\mathcal T_h}
\int_K|\nabla v|^2d\bm{x}\bigg)^{\frac{1}{2}}.
\]

\subsection{DG variational problem with the reconstructed operator}

Motivated by the statements of Remark \ref{why:dg}, we define a piecewise bilinear form $\mathcal A_h:\mathcal V\times \mathcal V\rightarrow\mathbb{R}$ as:
\begin{equation}\label{bfc}
\begin{split}
\mathcal A_h(v,w)=&\int_\Omega(c vw+\nu\nabla_hv\cdot\nabla_hw-
v\bm{b}\cdot\nabla_hw)\\
&+\sum_{e\in\mathcal E_h\backslash\Gamma^-}
\int_e\{\!\{\bm{b}v\}\!\}_{up}\cdot
\llbracket w \rrbracket
+\sum_{e\in\mathcal E_h}
\frac{\sigma_e\nu}{|e|^{1/(d-1)}}\int_e\llbracket v \rrbracket\cdot\llbracket w \rrbracket\\
&-\sum_{e\in\mathcal E_h}\int_e\big(\{\!\{\nu
\nabla_h v\}\!\}\cdot\llbracket w \rrbracket+\{\!\{\nu
\nabla_h w\}\!\}\cdot\llbracket v \rrbracket\big),
\end{split}
\end{equation}
where $\nabla_h$ denotes the gradient element by element, and $\sigma_e$ is the penalty parameter satisfying $0<\eta_0\leqslant \sigma_e\leqslant c_0$, $\forall e\in\mathcal E_h$. In addition, both $c_0$ and $\eta_0$ are global constants independent of any edge/face $e\in\mathcal E_h$ (see \cite{Arnold20021749}). We also define a piecewise linear form $\mathcal L:\mathcal V\rightarrow\mathbb{R}$ as:
\begin{equation}\label{lfc}
\mathcal L(v)=\int_\Omega fv+
\sum_{e\in\Gamma}
\int_e\bigg(\frac{c_e\nu}{|e|^{1/(d-1)}} v-\nu\nabla v\cdot\bm{n}\bigg)g
-\sum_{e\in\Gamma^-}\int_e\bm{b}\cdot\bm{n}gv.
\end{equation}
Then, the general DG variational problem with respect to the model problem \eqref{p1}--\eqref{p2} is: to find $u\in\mathcal V$ such that
\begin{equation}\label{vfc}
\mathcal A_h(u,v)=\mathcal L(v),\quad
\forall v\in\mathcal V.
\end{equation}

Now, using the reconstructed operator $\mathcal P^k$ defined in \eqref{global:op}, the discrete DG variational problem is: to find $u_h\in V_h^0$ such that
\begin{equation}\label{vfd}
\mathcal A_h(\mathcal P^ku_h, v)=\mathcal L(v),\quad\forall v\in V_h^k.
\end{equation}
Let $u$ be the solution of \eqref{vfc}, and it is a simple matter to check the consistency in the sense that
\begin{equation}\label{con:dg}
\mathcal A_h(u-\mathcal P^k u_h, v)=0,\quad\forall v\in V_h^k.
\end{equation}

\begin{remark}
We note that $\mathcal P^k V_h^0$ is equivalent to $V_h^k$ because $(\mathcal P^k v)|_K$ is a $k$-th polynomial for any $v\in V_h^0$.
\end{remark}

Next, we define the DG-energy norm as:
\begin{equation}\label{dg:norm}
\normdg{v}^2=\normdg{v}_{D}^2+\normdg{v}_{RC}^2,
\end{equation}
with
\begin{align*}
\normdg{v}_{D}^2&=\nu|v|_{1,h}^2+\sum_{e\in
\mathcal E_h}\frac{\nu}{|e|^{1/(d-1)}}
\|\llbracket v \rrbracket\|_{0,e}^2,\\
\normdg{v}_{RC}^2&=
\|(\overline{r}+b_0)^{1/2}v\|_{0,\Omega}^2
+\sum_{e\in\mathcal E_h}\||\bm{b}\cdot\bm{n}|^{1/2}
\llbracket v \rrbracket\|_{0,e}^2,
\end{align*}
where $\overline{r}$ is the piecewise constant function defined as:
\begin{equation}\label{reaction:innorm}
\overline{r}(\bm{x})|_K=\min_{\bm{x}\in K}r(\bm{x}),\quad\forall K\in\mathcal T_h.
\end{equation}
Inspired by the definition of the DG-energy norm \eqref{dg:norm}, we can also separate the bilinear form $\mathcal A_h$ into two parts. To fix ideas,
\[
\mathcal A_h(v,w)=a_h^{D}(v,w)+a_h^{RC}(v,w),
\]
with the diffusion part:
\begin{equation}\label{d:part}
\begin{split}
a_h^{D}(v,w)=&\int_\Omega\nu\nabla_hv\cdot\nabla_hw
+\sum_{e\in\mathcal E_h}
\frac{c_e\nu}{|e|^{1/(d-1)}}\int_e\llbracket v \rrbracket\cdot\llbracket w \rrbracket\\
&-\sum_{e\in\mathcal E_h}\int_e\big(\{\!\{\nu
\nabla_h v\}\!\}\cdot\llbracket w \rrbracket+\{\!\{\nu
\nabla_h w\}\!\}\cdot\llbracket v \rrbracket\big),
\end{split}
\end{equation}
and the reaction-convection part:
\begin{equation}\label{rc:part}
a_h^{RC}(v,w)=\int_\Omega(\gamma vw-
v\bm{b}\cdot\nabla_hw)
+\sum_{e\in\mathcal E_h\backslash\Gamma^-}
\int_e\{\!\{\bm{b}v\}\!\}_{up}\cdot
\llbracket w \rrbracket.
\end{equation}

\subsection{Stability}
\label{sec:dis:stable}

Now, based on the ingenious techniques in \cite{Ayuso20091391}, we shall verify the stability of \eqref{vfd}. In \ref{append:pf:stable} we provide proofs of all the following lemmas and theorems.

First, the following inequality will be extensively employed.

\begin{lemma}\label{useful:inequality}
For any $v_h\in V_h^k$ and any $w\in \mathcal V$, there exists a positive constant $\kappa_9$ that depends only on $\mathcal N_0$, $N$, $\gamma$, $k$, and $d$, such that
\begin{equation}\label{dg:inequality}
\sum_{e\in\mathcal E_h}\int_e\Big|
\{\!\{\nu\nabla_h v_h\}\!\}
\cdot\llbracket w \rrbracket\Big|
\leqslant \kappa_9\nu
|v_h|_{1,h}
\bigg(\sum_{e\in
\mathcal E_h}\frac{1}{|e|^{1/(d-1)}}
\|\llbracket w \rrbracket\|_{0,e}^2\bigg)^{\frac{1}{2}}.
\end{equation}
\end{lemma}

Then, we introduce a weighted function $\mu=\exp(-\xi)+\delta$ where $\xi$ is defined in (H1) and $\delta$ is a positive constant. In addition, from (H1) and \eqref{velocity:assumption} we propose the demands for $\mu$ that
\begin{align}
\label{wc:property}
&\mu_1\leqslant \exp(-\xi)\leqslant\mu_2,\qquad
|\nabla\mu|\leqslant\mu_3,\\
\label{wc:property2}
&\mu_1+\delta>\eta_1LC_P\kappa_{10}\mu_3,
\quad
\mu_1+\delta>\eta_2(\mu_2+\delta),
\end{align}
where $\mu_1$, $\mu_2$, $\mu_3$ are three positive constants, $\kappa_{10}=\max\{\kappa_9,1\}$, $C_P$ is the Poincar\'{e} constant in \eqref{Poincare:dg}, and $\eta_1$, $\eta_2$ are two positive constants with $\eta_2<1$. The following lemma can be regarded as a ``weak" version of coercivity because the test function is chosen as $\mu v_h$ instead of $v_h$.

\begin{lemma}\label{lemma:coercivity}
Suppose there exists a constant $\theta_0>2/\eta_2$ such that
\begin{equation}\label{assum:eta0}
\eta_0\geqslant \max\bigg\{\frac{\theta_0^2}{4}\kappa_9^2,1\bigg\}
\quad\text{and}\quad
\frac{\sqrt{2}}{\eta_1}+\frac{2}{\eta_2\theta_0}<1,
\end{equation}
then for any $v_h\in V_h^k$, the following inequalities hold:
\begin{align}
&a_h^{D}(v_h,\mu v_h)\geqslant \kappa_{11}(\mu_1+\delta)\normdg{v_h}_{df}^2,
\label{coer:diff}\\
&a_h^{RC}(v_h,\mu v_h)\geqslant \kappa_{12}\normdg{v_h}_{rc}^2,
\label{coer:rc}\\
&\normdg{\mu v_h}\leqslant \kappa_{13}(\mu_1+\delta)\normdg{v_h}.\label{bound:norm}
\end{align}
where $\kappa_{11}=1-\frac{\sqrt{2}}{\eta_1}-\frac{2}{\eta_2\theta_0}$, $\kappa_{12}=\min\{\mu_1,\frac{\mu_1+\delta}{2}\}$ and $\kappa_{13}=\frac{\sqrt{2(\eta_1^2\kappa_{10}^2+\eta_2^2)}}{\eta_1\eta_2\kappa_{10}}$.
\end{lemma}

In \cite{Ayuso20091391} the mesh is triangulation, however, we consider general polytopic meshes throughout our work. Hence, it is necessary to recount the following polynomial interpolation results.

\begin{lemma}\label{super:appro}
Suppose $\mu\in W^{k+1,\infty}(\Omega)$. For any $v_h\in V_h^k$, let $\widetilde{\mu v_h}$ be the $L^2$-projection of $\mu v_h$ in $V_h^k$, then we have
\begin{align}
&\|\mu v_h-\widetilde{\mu v_h}\|_{0,\Omega}\leqslant\kappa_{14}\frac{\|\overline{\xi}\|_{k+1,\infty,\Omega}}{L}h
\|v_h\|_{0,\Omega},\label{project:L2}\\
&|\mu v_h-\widetilde{\mu v_h}|_{1,h}\leqslant\kappa_{15}\frac{\|\overline{\xi}\|_{k+1,\infty,\Omega}}{L}
\|v_h\|_{0,\Omega},\label{project:H1}\\
&\bigg(\sum_{e\in \mathcal E_h}\|\mu v_h-\widetilde{\mu v_h}\|_{0,e}^2\bigg)^{\frac{1}{2}}
\leqslant\kappa_{16}\frac{\|\overline{\xi}\|_{k+1,\infty,\Omega}}{L}h^{\frac{1}{2}}
\|v_h\|_{0,\Omega},\label{project:trace}
\end{align}
where $\overline{\xi}=\exp(-\xi)$, and $\kappa_{14}$, $\kappa_{15}$, $\kappa_{16}$ are three positive constants that depend only on $\mathcal N_0$, $N$, $\gamma$, and $k$.
\end{lemma}

In addition, we also need the following lemma to ensure the boundedness of $\mathcal A_h$ with $\mu v_h-\widetilde{\mu v_h}$ as the test function.

\begin{lemma}\label{lemma:contin}
In the hypotheses of Lemma \ref{lemma:coercivity}, for any $v_h\in V_h^k$, there exist two positive constants $\kappa_{17}$ and $\kappa_{18}$, such that
\begin{align}
& \big|a_h^{D}(v_h,\mu v_h-\widetilde{\mu v_h})\big|\leqslant \kappa_{17}\normdg{v_h}_{D}^2,
\label{quasi:contin:df}\\
& \big|a_h^{RC}(v_h,\mu v_h-\widetilde{\mu v_h})\big|\leqslant \kappa_{18}\bigg(\frac{h}{L}
\bigg)^{\frac{1}{2}}\normdg{v_h}_{RC}^2.
\label{quasi:contin:rc}
\end{align}
\end{lemma}

Finally, based on Lemma \ref{lemma:coercivity} -- Lemma \ref{lemma:contin}, the next theorem provides the stability result for the discrete DG variational problem \eqref{vfd} in the norm $\normdg{\cdot}$.

\begin{theorem}\label{first:stability}
In the case of $\kappa_{11}(\mu_1+\delta)\geqslant 2\kappa_{17}$, there exist two positive constants $\alpha_1$ and $h_0$ (independent of $h$), such that
for $h<h_0$,
\begin{equation}\label{first:infsup}
\sup_{w_h\in V_h^k}\frac{\mathcal A_h(v_h, w_h)}{\normdg{w_h}}\geqslant
\alpha_1\normdg{v_h},\quad\forall v_h\in V_h^k.
\end{equation}
\end{theorem}

\subsection{Stability in a convection-dominated regime}
\label{sec:cd}

In a convection-dominated regime, we
desire to have a control also on the streamline derivative. It means that the stability is needed in another norm including a term of SUPG type. Hence, based on the norm $\normdg{\cdot}$, we define a new DG-energy norm of SUPG type as:
\begin{equation}\label{supg:norm}
\normdg{v}_{S}=\big(\normdg{v}^2+\|v\|_{b}^2\big)^{\frac{1}{2}},
\quad
\|v\|_{b}=\bigg(\sum_{K\in\mathcal T_h}\frac{h_K}{\|\bm{b}\|_{0,\infty,K}}
\|\mathcal P^k(\bm{b}\cdot\nabla v)\|_{0,K}^2\bigg)^{\frac{1}{2}}.
\end{equation}
In addition, when the convection dominates, i.e., $0<\nu\ll 1$ described in Section \ref{sec:Intro}, we assume that
\begin{equation}\label{convection:dominate:assum}
\nu<\frac{h_K\|\bm{b}\|_{0,\infty,K}}{2},\quad
\|c+\nabla\cdot\bm{b}\|_{0,\infty,K}\leqslant\frac{\|\bm{b}\|_{0,\infty,K}}{h_K},
\quad\forall K\in\mathcal T_h.
\end{equation}

\begin{remark}\label{remark:mesh:convection}
We must stress that the assumption \eqref{convection:dominate:assum} is valid only in a convection-dominated regime, while it does not hold when $\nu$ is not a small parameter.
\end{remark}

The stability in the norm $\normdg{\cdot}_S$ can again be achieved through an inf-sup condition. In \ref{append:pf:stable} we provide proofs of Lemma \ref{lemma:supg:preinfsup} and Theorem \ref{second:stability}.

\begin{lemma}
\label{lemma:supg:preinfsup}
In the hypotheses of \eqref{convection:dominate:assum}, there exist two positive constants $\alpha_2$, $\alpha_3$ both independent of $h$, $\nu$, $\bm{b}$ and $c$, such that
\begin{equation}\label{supg:preinfsup}
\sup_{w_h\in V_h^k}\frac{\mathcal A_h(v_h, w_h)}{\normdg{w_h}}\geqslant
\alpha_2\|v_h\|_b-
\alpha_3\normdg{v_h},\quad\forall v_h\in V_h^k.
\end{equation}
\end{lemma}

\begin{theorem}
\label{second:stability}
There exists a positive constant $\alpha_4$ such that for $h<h_0$,
\begin{equation}\label{second:infsup}
\sup_{w_h\in V_h^k}\frac{\mathcal A_h(v_h, w_h)}{\normdg{w_h}}\geqslant
\alpha_4\normdg{v_h}_S,\quad\forall v_h\in V_h^k.
\end{equation}
\end{theorem}

\section{A priori error estimates}
\label{prior:error}

In \ref{append:pf:4} we prove the proofs of all lemmas and theorems in this section.

\subsection{Error estimates in the DG energy norms}
We shall provide a priori error estimates for the discrete DG variational problem \eqref{vfd} in the norm $\normdg{\cdot}$ and $\normdg{\cdot}_S$, respectively.

\begin{theorem}\label{priori:dg}
Let $u$ and $u_h$ be the solutions of \eqref{vfc} and \eqref{vfd}, respectively. In the hypotheses of (A4), there exists a positive constant $C_1$ independent of $h$, $\nu$, $\bm{b}$ and $c$, such that
\begin{equation}\label{priori:dg:ineq}
\normdg{u-\mathcal P^k u_h}
\leqslant C_1\big(\nu^{\frac{1}{2}}
+\|\bm{b}\|_{0,\infty,\Omega}^{\frac{1}{2}}
h^{\frac{1}{2}}+\|r\|_{0,\infty,\Omega}^{\frac{1}{2}}
h\big)h^k|u|_{k+1,\Omega}.
\end{equation}
\end{theorem}

In a convection-dominated regime, owing to Theorem \ref{second:stability} and via the same steps of proving Theorem \ref{priori:dg} (see \ref{append:pf:th41}), the following theorem provides the optimal error estimate in the norm $\normdg{\cdot}_S$. So we shall omit the details of its proof.

\begin{theorem}\label{priori:supg}
Let $u$ and $u_h$ be the solutions of \eqref{vfc} and \eqref{vfd}, respectively. In the hypotheses of (A4) and \eqref{convection:dominate:assum}, there exists a positive constant $C_2$ independent of $h$, $\nu$, $\bm{b}$ and $c$, such that
\begin{equation}\label{priori:supg:ineq}
\normdg{u-\mathcal P^k u_h}_S
\leqslant C_2\big(\nu^{\frac{1}{2}}
+\|\bm{b}\|_{0,\infty,\Omega}^{\frac{1}{2}}
h^{\frac{1}{2}}+\|r\|_{0,\infty,\Omega}^{\frac{1}{2}}
h\big)h^k|u|_{k+1,\Omega}.
\end{equation}
\end{theorem}

\subsection{Error estimates in the $L^2$-norm}

First, we stress that the error estimates shown below in the $L^2$-norm are available only for the divergence-free velocity field $\bm{b}$, which is up to our best efforts. Let $u$ and $u_h$ be the solutions of \eqref{vfc} and \eqref{vfd}, respectively. Then, based on \cite{BdaVeiga2016729}, we define the auxiliary problem with respect to the model problem \eqref{p1}--\eqref{p2} as: to find $\psi\in H^2(\Omega)$ such that
\begin{equation}\label{var:aux}
\begin{cases}
-\nu\Delta \psi+\bm{b}\cdot\nabla\psi+c\psi=u-\mathcal P^ku_h,&\text{in}~\Omega,\\
\psi=0,&\text{on}~\partial\Omega.
\end{cases}
\end{equation}

\begin{remark}
Using the fact that $\bm{b}$ is divergence free, we have
\[
\nabla\cdot(\bm{b}\psi)=(\nabla\cdot\bm{b})\psi+\bm{b}\cdot
\nabla \psi=\bm{b}\cdot
\nabla \psi.
\]
\end{remark}

\begin{lemma}\label{auxiliary:regularity}
Let $\psi$ be the solution of \eqref{var:aux}. If the domain $\Omega$ is convex, then the solution $\psi$ satisfies the following estimate:
\begin{equation}\label{var:stable}
\|\psi\|_{2,\Omega}\leqslant C_3(\nu^{-2}\|\bm{b}\|_{0,\infty,\Omega}
+\nu^{-2}\|r\|_{0,\infty,\Omega}+\nu^{-1})\|u-\mathcal P^ku_h\|_{0,\Omega},
\end{equation}
where $C_3$ is a positive constant independent of $\nu$, $\bm{b}$, and $c$.
\end{lemma}

\begin{theorem}\label{priori:l2}
Let $u$ and $u_h$ be the solutions of \eqref{vfc} and \eqref{vfd}, respectively. In the hypotheses of (A4), there exist two different positive constants $C_4$ and $C_5$ such that
\begin{equation}\label{L2result}
\|u-\mathcal P^ku_h\|_{0,\Omega}\leqslant |u|_{k+1,\Omega}
\begin{cases}
C_4h^{k+1}&\text{if diffusion dominates},\\
C_5h^{k+\frac{1}{2}}&\text{if convection dominates}.
\end{cases}
\end{equation}
\end{theorem}

\begin{remark}
If the convection dominates, the norm $\normdg{\cdot}$ should be replaced by the norm $\normdg{\cdot}_S$. In addition, the optimal error estimate of Theorem \ref{priori:supg} ensures the results \eqref{L2result} are available.
\end{remark}
\begin{remark}
In spite of the suboptimal error order when the convection dominates, we shall show that the optimal $L^2$-error order can be still achieved in the numerical examples of Section \ref{numer:ex} if the exact solution is smooth enough.
\end{remark}
\begin{remark}
As mentioned above, we can just prove the error estimates in the $L^2$-norm with the divergence-free $\bm{b}$ in theory. If $\bm{b}$ is not divergence free, the auxiliary problem \eqref{var:aux} first becomes: to find $\psi\in H^2(\Omega)$ such that
\begin{equation*}
\begin{cases}
-\nu\Delta \psi+\nabla\cdot(\bm{b}\psi)+c\psi=u-\mathcal P^ku_h,&\text{in}~\Omega,\\
\psi=0,&\text{on}~\partial\Omega.
\end{cases}
\end{equation*}
Due to the asymmetry of $\mathcal A_h$, \eqref{L2} in \ref{append:pf:th44} has to be written as:
\begin{align*}
\|u-\mathcal P^k u\|_{0,\Omega}^2&=(-\nu\Delta\psi+\bm{b}\cdot\psi+
c\psi,u-\mathcal P^k u_h)_\Omega+(\psi\nabla\cdot\bm{b},u-\mathcal P^k u_h)_\Omega\\
&=\mathcal A_h(u-\mathcal P^k u_h,\psi)+(\psi\nabla\cdot\bm{b},u-\mathcal P^k u_h)_\Omega.
\end{align*}
Without any projection or approximation property, the additional term $(\psi\nabla\cdot\bm{b},u-\mathcal P^k u_h)_\Omega$ is too tricky to reach our desirable approximation results.
\end{remark}

\section{Numerical experiments}
\label{numer:ex}

In this section we present various test examples to show the rates of convergence and the performance of the reconstruction operator. All the experiments are performed on the unit square $\Omega=(0,1)^2$ and for different reconstruction order $k$ we take the penalty parameter $\sigma_e=3k(k+1)$. Numerical errors will be measured in the $L^2$-norm $\|\cdot\|$, DG-energy norm $\normdg{\cdot}$ and/or $\normdg{\cdot}_S$, respectively.

\begin{figure}[htbp]
\centering
\subfigure[]{
\includegraphics[width=4.5cm]{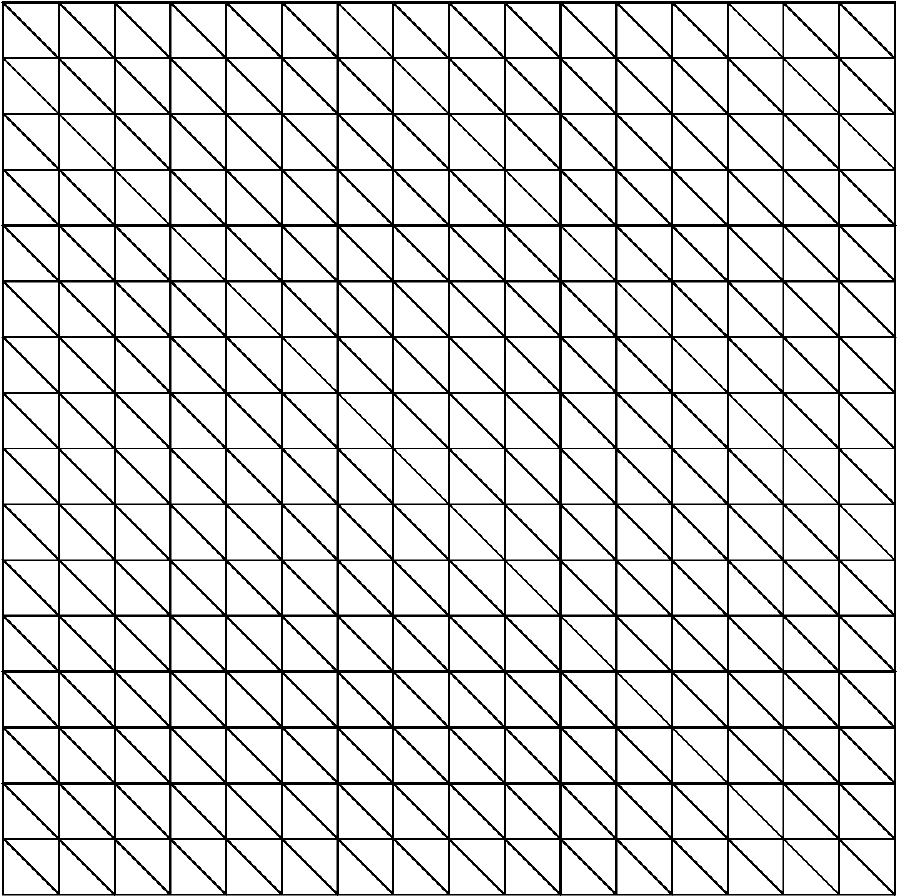}
}
\quad
\subfigure[]{
\includegraphics[width=4.5cm]{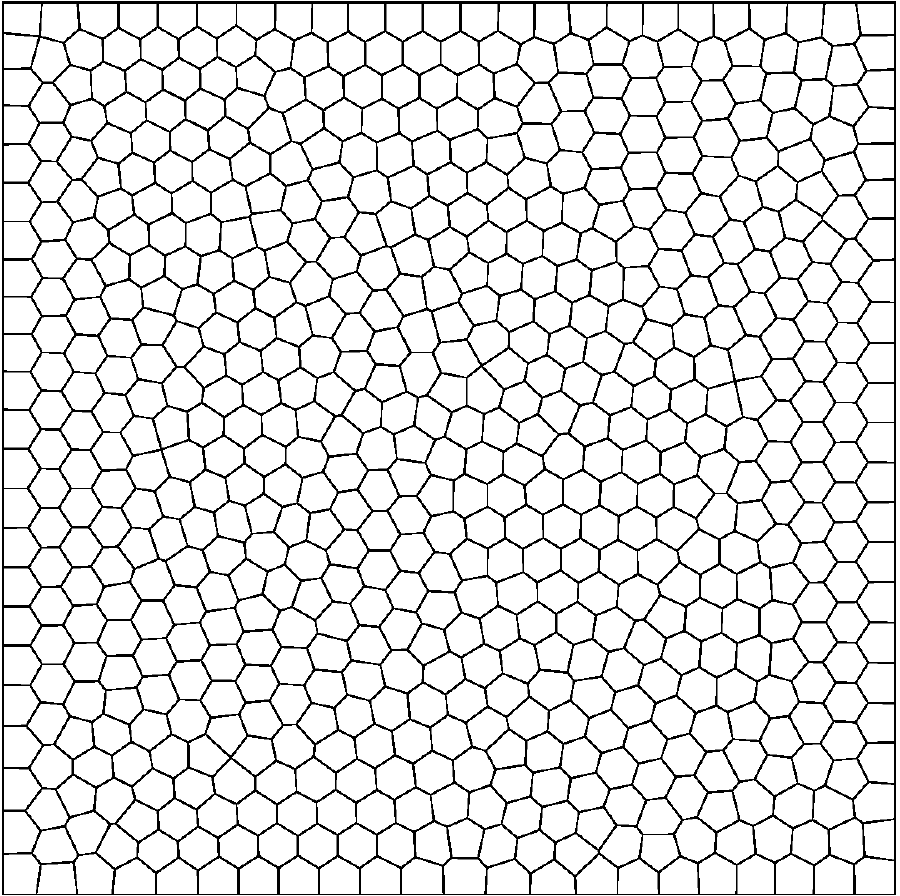}
}
\quad
\subfigure[]{
\includegraphics[width=4.5cm]{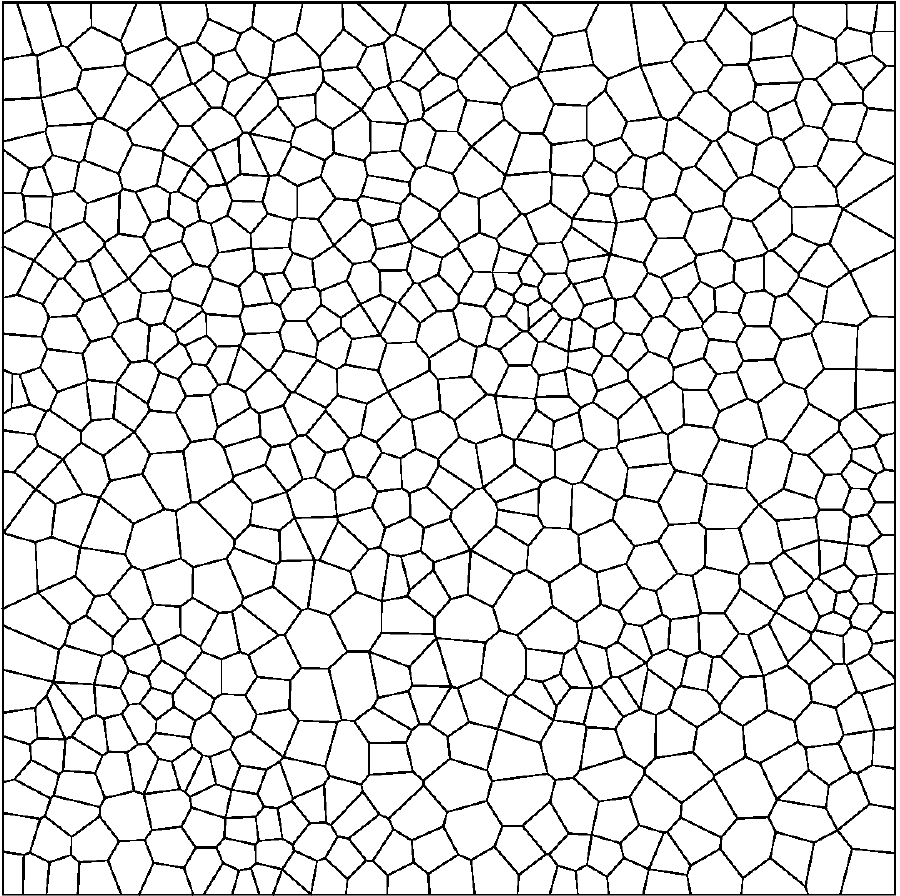}
}
\centering
\caption{The examples of three mesh families: (a) uniform shape-regular triangulation; (b) regular polygonal mesh; (c) general Voronoi mesh.}\label{fig:mesh:example}
\end{figure}

As shown in Figure \ref{fig:mesh:example}, three different meshes are employed to our experiments: uniform triangulations, regular polygonal meshes and general Voronoi meshes. Both the polygonal meshes and the general Voronoi meshes are generated by \emph{PolyMesher} \cite{Talischi2012309}. Given a partition $\mathcal T_h$, we list a group of reference values of $\#S(K)$ in Table \ref{table1} for different order $k$. A direct solver is employed to solve all the resulting linear systems.

\begin{table}[!ht]
    \centering
    \topcaption{Choices of $\#S(K)$ for different reconstruction order $k$}
    \begin{tabular*}{\hsize}{@{}@{\extracolsep{\fill}}lccc@{}}
        \toprule
        Reconstruction order $k$ & 1 & 2 & 3\\
        \hline
        $\#S(K)$ for uniform triangulations  & 4   & 7   & 11\\
        $\#S(K)$ for regular polygonal meshes   & 5   & 9   & 15\\
        $\#S(K)$ for general Voronoi meshes     & 7   & 12   & 19\\
        \bottomrule
    \end{tabular*}
    \label{table1}
\end{table}

\begin{algorithm}[!htbp]
\caption{An implementation for the convection-diffusion-reaction equation in 2D}
\label{alg:implement:method}
\begin{algorithmic}[1]
\Require
$\mathcal T_h$, $k$, $M$, $\nu$, $\bm{b}$, $c$, $f$, and $g$;
\Ensure
the numerical solution $u_h(\bm{x})$ of the model problem \eqref{p1}--\eqref{p2};
\State set $NT$ as the number of elements of $\mathcal T_h$;
\State compute the global index $N_K$ for each element $K\in\mathcal T_h$, $0\leqslant N_K\leqslant NT-1$;
\State initialize $\bm{S}$ as a $NT$ by $M$ zero matrix;
\For{$i=0$ to $NT-1$}
\State compute $\bm{S}[i,:]=(N_{E_0}, N_{E_1},
\cdots, N_{E_{M-1}})$ based on Algorithm \ref{alg:element:patch};
\EndFor
\State initialize $\bm{\Phi}$ as a $NT$ by $NT$ zero matrix;
\For{$i=0$ to $NT-1$}
\For{$j=0$ to $NT-1$}
\If{$\textsf{ismember}$($i$, $\bm{S}[j,:]$)}
\State find the index of $i$ in $\bm{S}[j,:]$, and denote it by $j_i$;
\State compute $\bm{G}$ based on \eqref{solution:wls} for the element whose global index is $j$;
\State compute $\lambda_{j_i}(\bm{x})$ defined in \eqref{re:poly:2}, and denote it by $a_i^{(j)}$;
\State set $\bm{\Phi}[i,j]=a_i^{(j)}$;
\EndIf
\EndFor
\EndFor
\For{$i=0$ to $NT-1$}
\State take all terms $\bm{\Phi}[i,j]$ from $\bm{\Phi}[i,:]$, $j=0,1,\cdots,NT-1$;
\State set the basis function $\phi_i$ as $\phi_i|_{K_j}=\bm{\Phi}[i,j]$ for the element $K_i\in\mathcal T_h$;
\EndFor
\State set the stiff matrix $\bm{M}$ as a $NT$ by $NT$ matrix with $\bm{M}[i,j]=\mathcal A_h(\phi_i,\phi_j)$, $i,j=0,1,\cdots,NT-1$;
\State set the vector $\bm{u}$ as $\bm{u}[j,]=u_j$ with $u_j$ unknown value at the barycenter of the element $K_j$, $j=0,1,\cdots,NT-1$;
\State set the load vector $\bm{F}$ as $\bm{F}[j,]=\mathcal L(\phi_j)$, $j=0,1,\cdots,NT-1$;
\State solve the linear algebraic system $\bm{M}\bm{u}=\bm{F}$ to obtain the approximate value $\widetilde{u}_i$ at the barycenters of the element $K_i$, $i=0,1,\cdots,NT-1$;
\State set the final numerical solution $\displaystyle u_h(\bm{x})=\sum_{i=0}^{NT-1}\widetilde{u}_i\phi_i(\bm{x})$.
\end{algorithmic}
\end{algorithm}

\subsection{Basic algorithm for implementing our proposed method}

Before discussing the test examples, in Algorithm \ref{alg:implement:method} we present the details of solving the model problem \eqref{p1}--\eqref{p2} with our proposed DG method. An highlight is that only one degree of freedom per element is needed to achieve arbitrary-order accuracy.

By reading each step of Algorithm \ref{alg:implement:method}, we can claim that the whole implementation is simple. Further the scale of the stiff matrix $\bm{M}$ depends only on $NT$, the number of elements of $\mathcal T_h$, and it is much smaller than the standard DG methods.

\subsection{Example 1: case of smooth solution}

\subsubsection{The convection field is not divergence free}
The variable coefficients of \eqref{p1} are given by $\bm{b}(x,y)=[x^2y+1,~xy^2+1]^\top$, $c=1$, and we vary the diffusion coefficient $\nu=1,10^{-3},10^{-9}$. The source term $f$ is chosen so that the analytical solution of \eqref{p1}, with Dirichlet boundary condition \eqref{p2}, is given by
\begin{equation}\label{exact:solution}
u(x,y)=\sin(2\pi x)\sin(2\pi y).
\end{equation}
Here the uniform triangulations and regular polygonal meshes are taken to show the convergence results. In particular, for regular polygonal meshes we use the total number of the elements to be 10, 40, 160, 640 and 2560 to simulate the effect of uniform refinement.

\begin{figure}[htbp]
\centering
\subfigure[]{
\includegraphics[width=7cm]{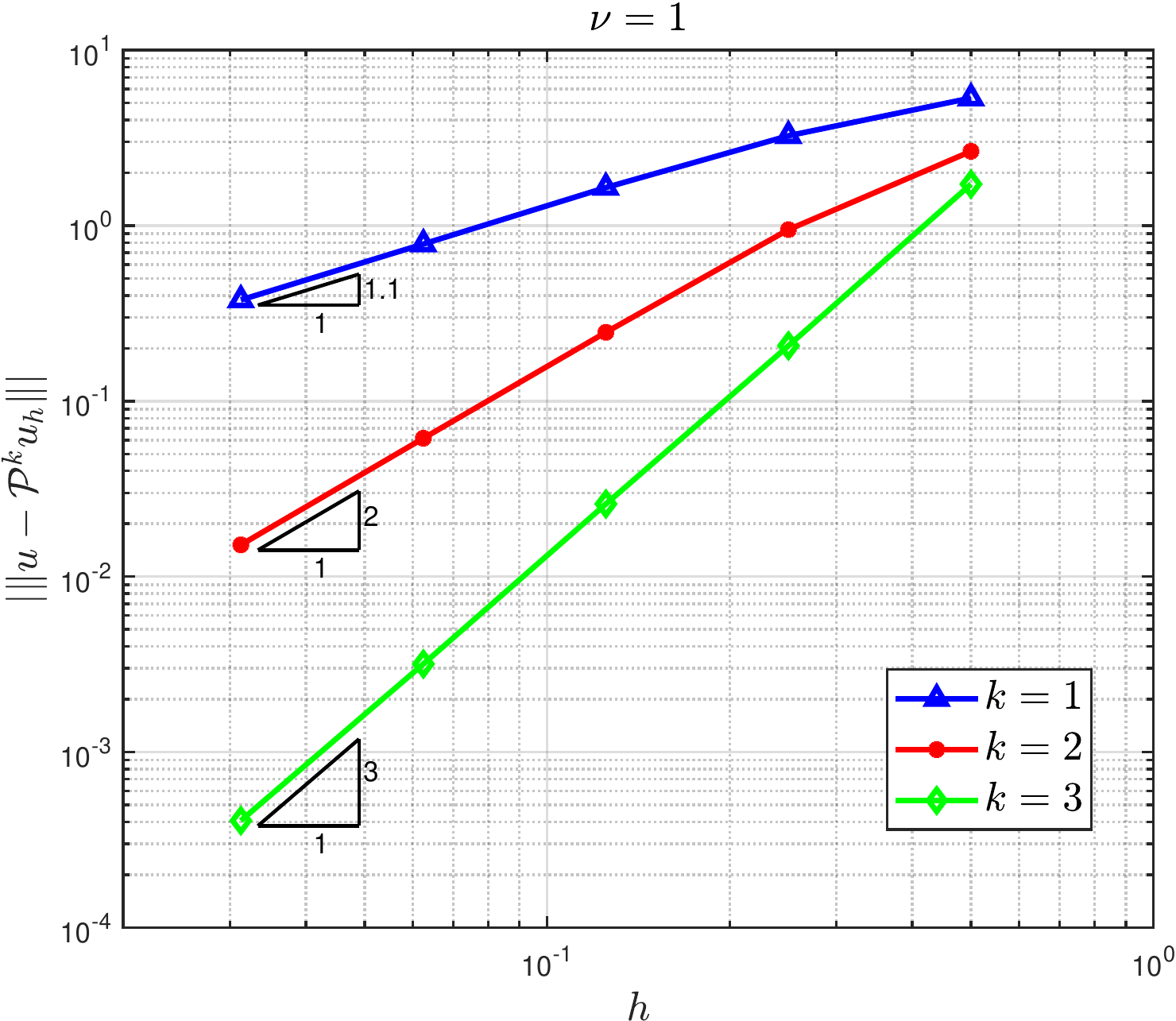}
}

\subfigure[]{
\includegraphics[width=7cm]{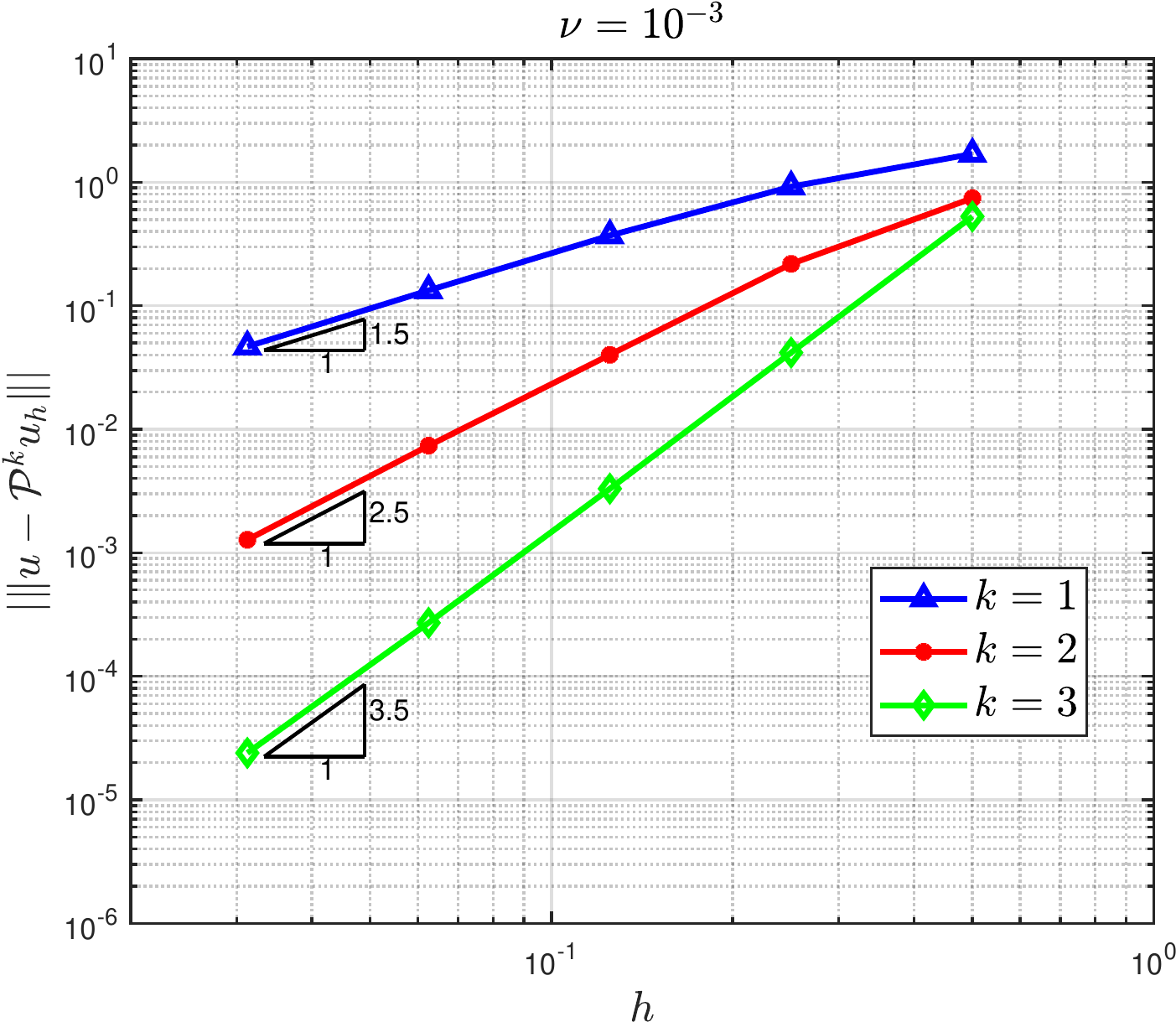}
}
\quad
\subfigure[]{
\includegraphics[width=7cm]{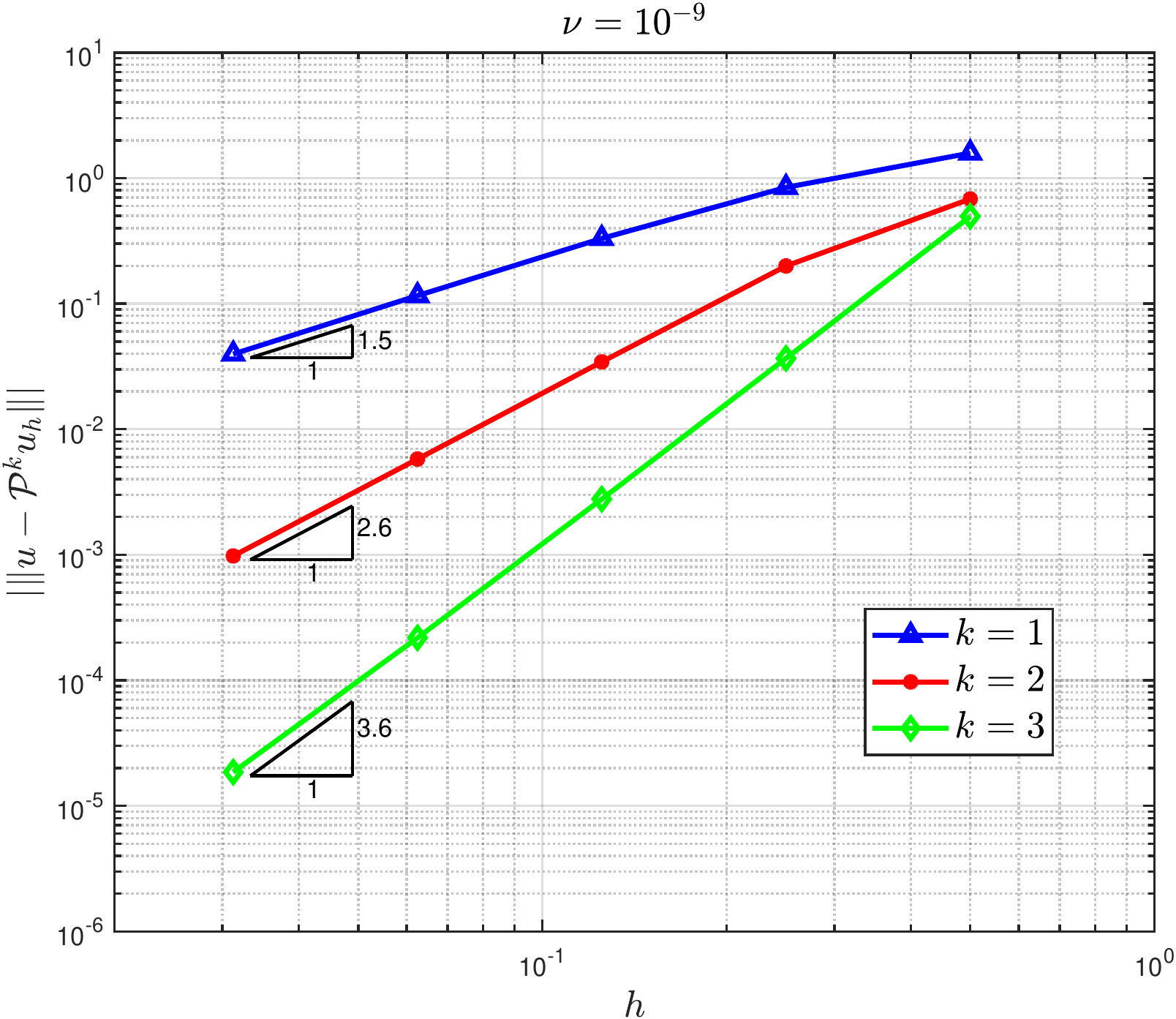}
}

\subfigure[]{
\includegraphics[width=7cm]{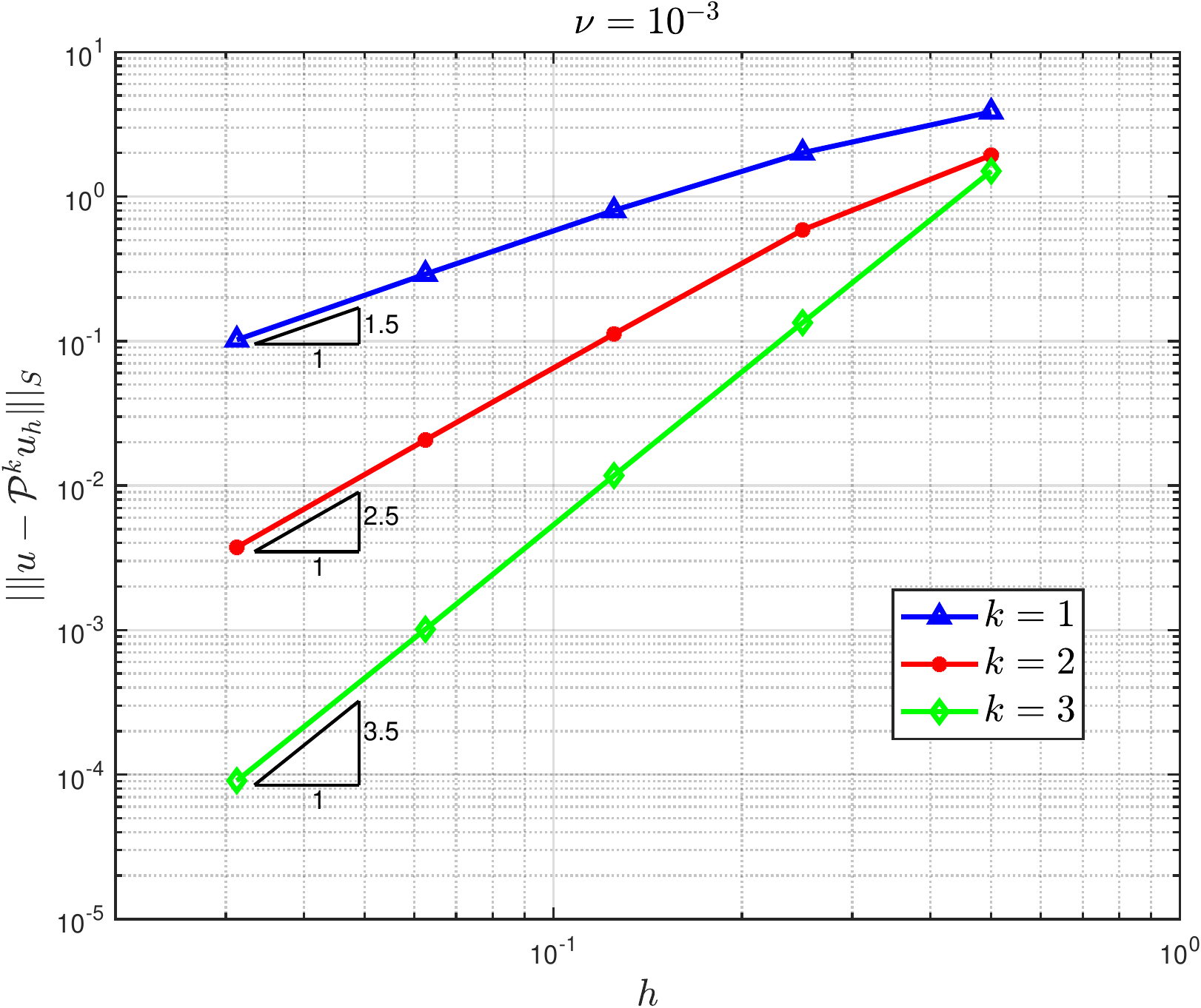}
}
\quad
\subfigure[]{
\includegraphics[width=7cm]{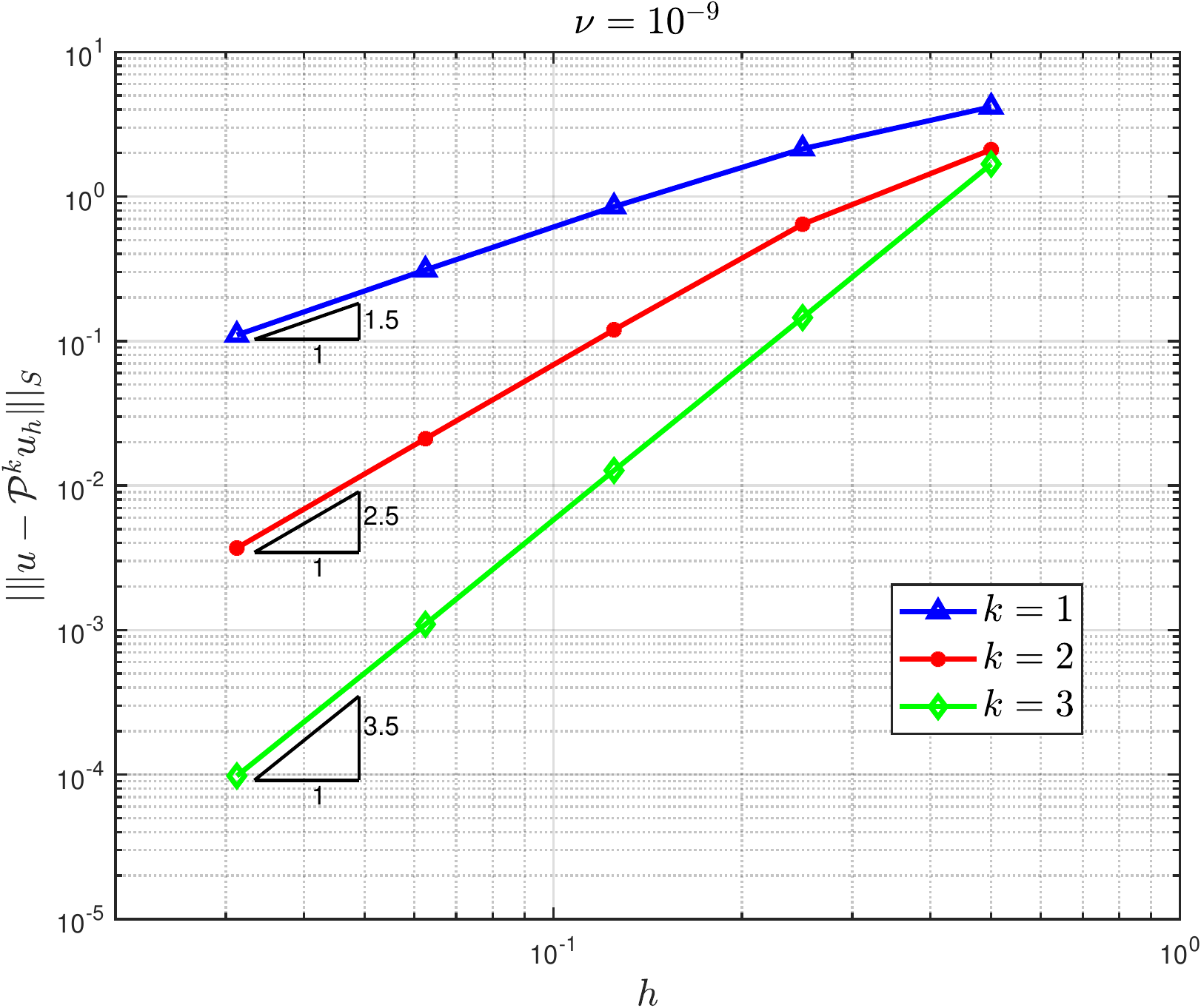}
}
\centering
\caption{\emph{Example 1.} Convergence rate results of $\normdg{u-\mathcal P^ku_h}$ and $\normdg{u-\mathcal P^ku_h}_S$ over uniform triangulations with a not divergence-free convection field: (a) $\normdg{u-\mathcal P^ku_h}$ with $\nu=1$; (b) $\normdg{u-\mathcal P^ku_h}$ with $\nu=10^{-3}$; (c) $\normdg{u-\mathcal P^ku_h}$ with $\nu=10^{-9}$; (d) $\normdg{u-\mathcal P^ku_h}_S$ with $\nu=10^{-3}$; (e) $\normdg{u-\mathcal P^ku_h}_S$ with $\nu=10^{-9}$.}\label{fig:example1:cr:l2dg:tri}
\end{figure}

\begin{figure}[htbp]
\centering
\subfigure[]{
\includegraphics[width=7cm]{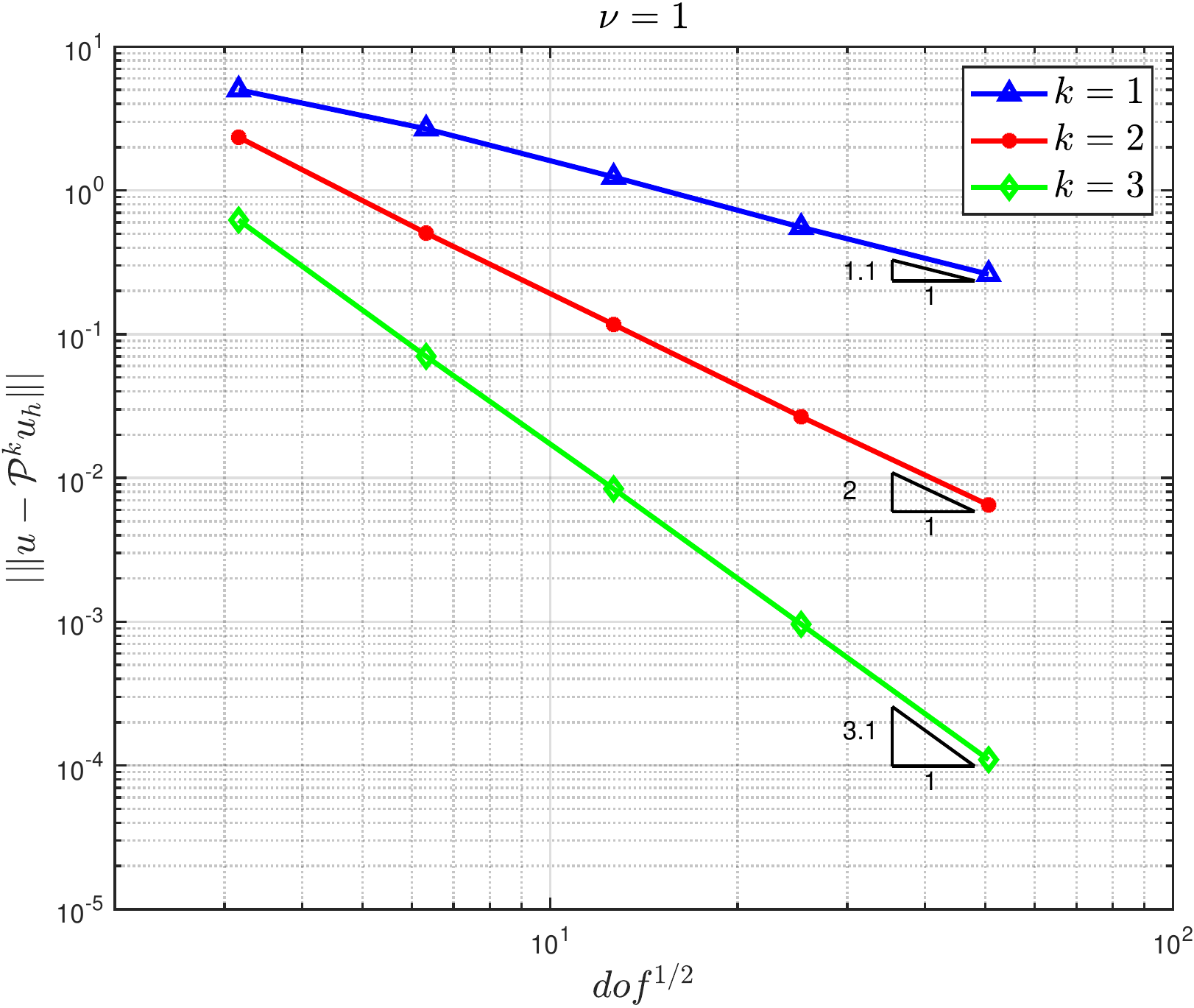}
}

\subfigure[]{
\includegraphics[width=7cm]{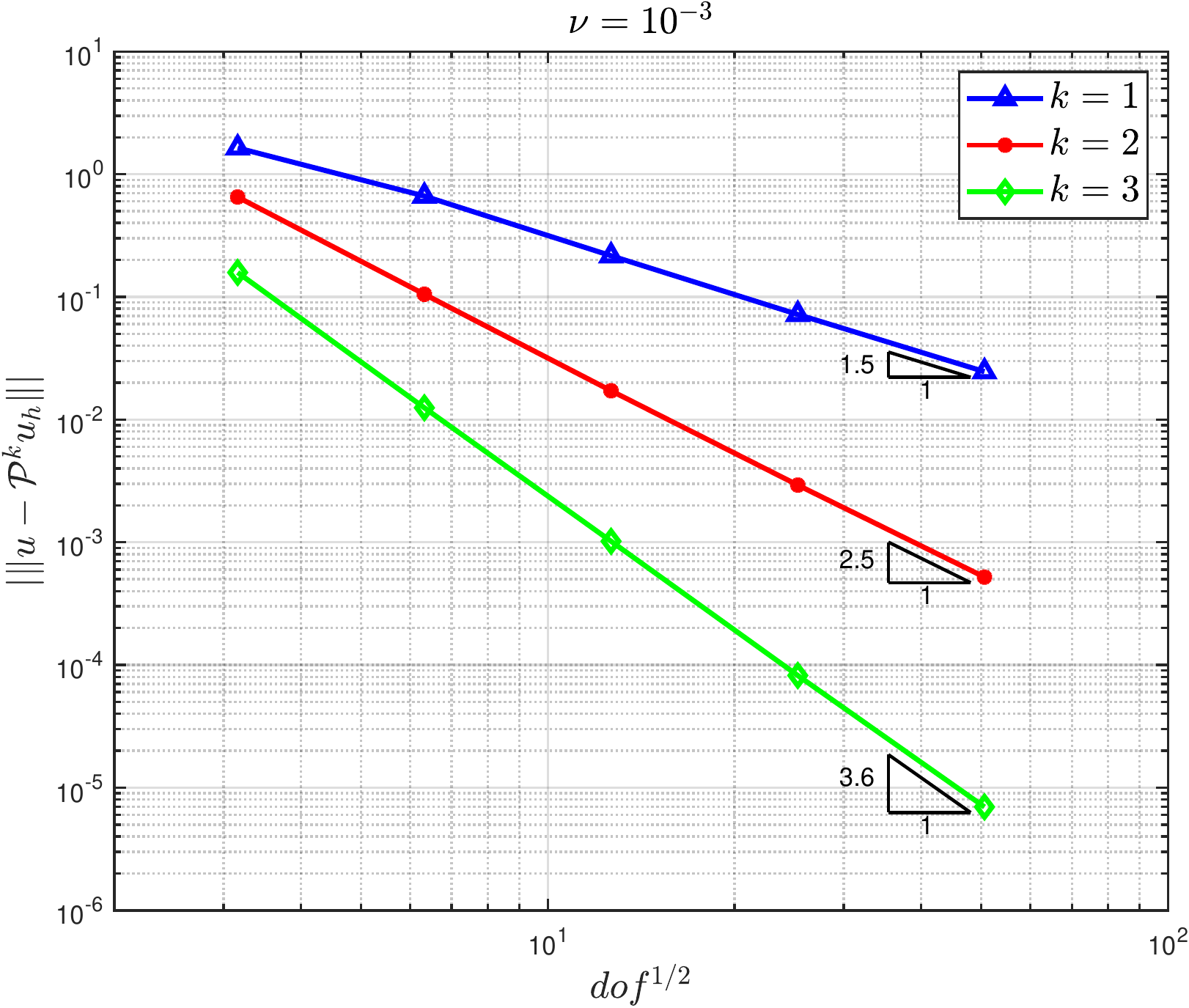}
}
\quad
\subfigure[]{
\includegraphics[width=7cm]{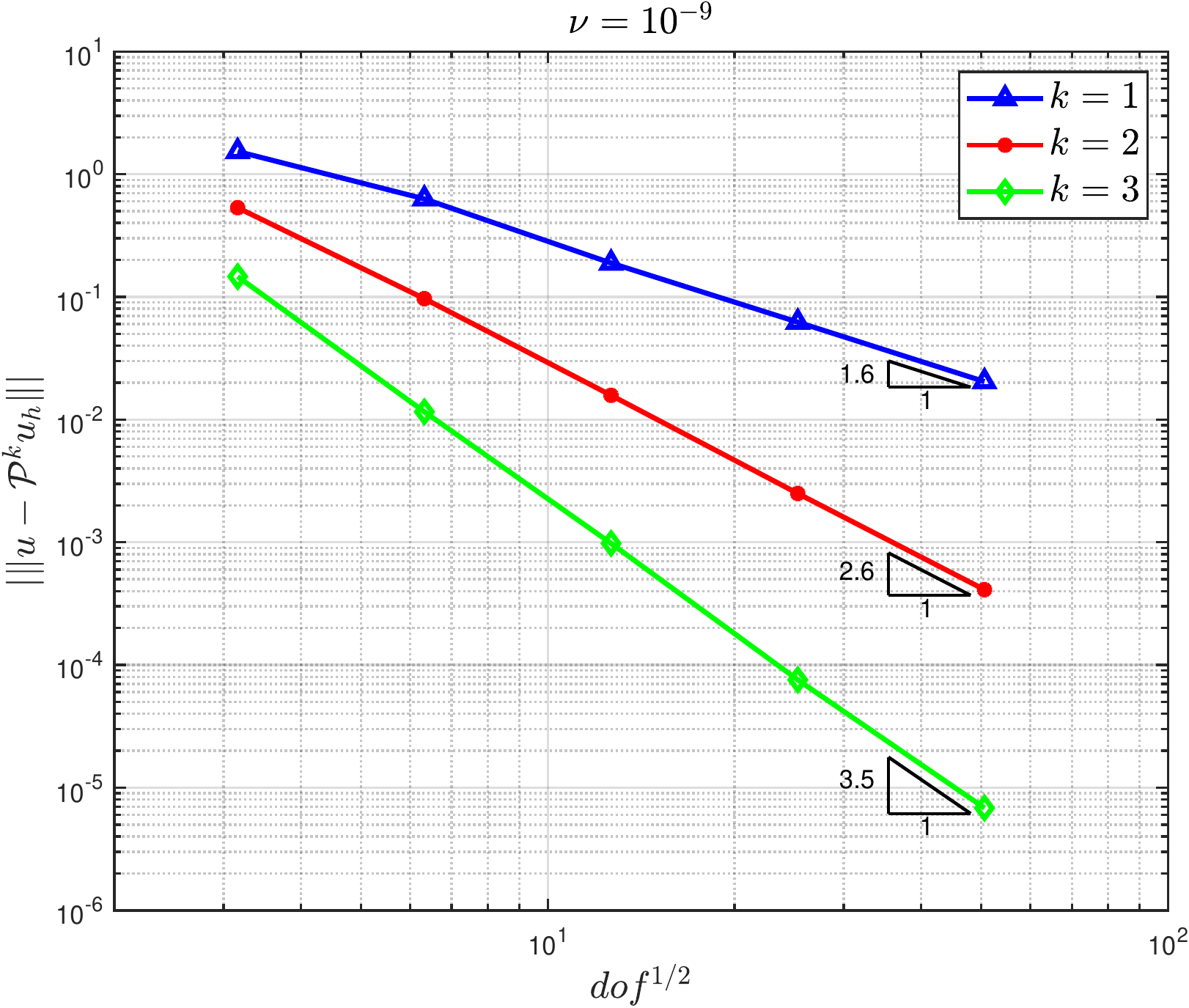}
}

\subfigure[]{
\includegraphics[width=7cm]{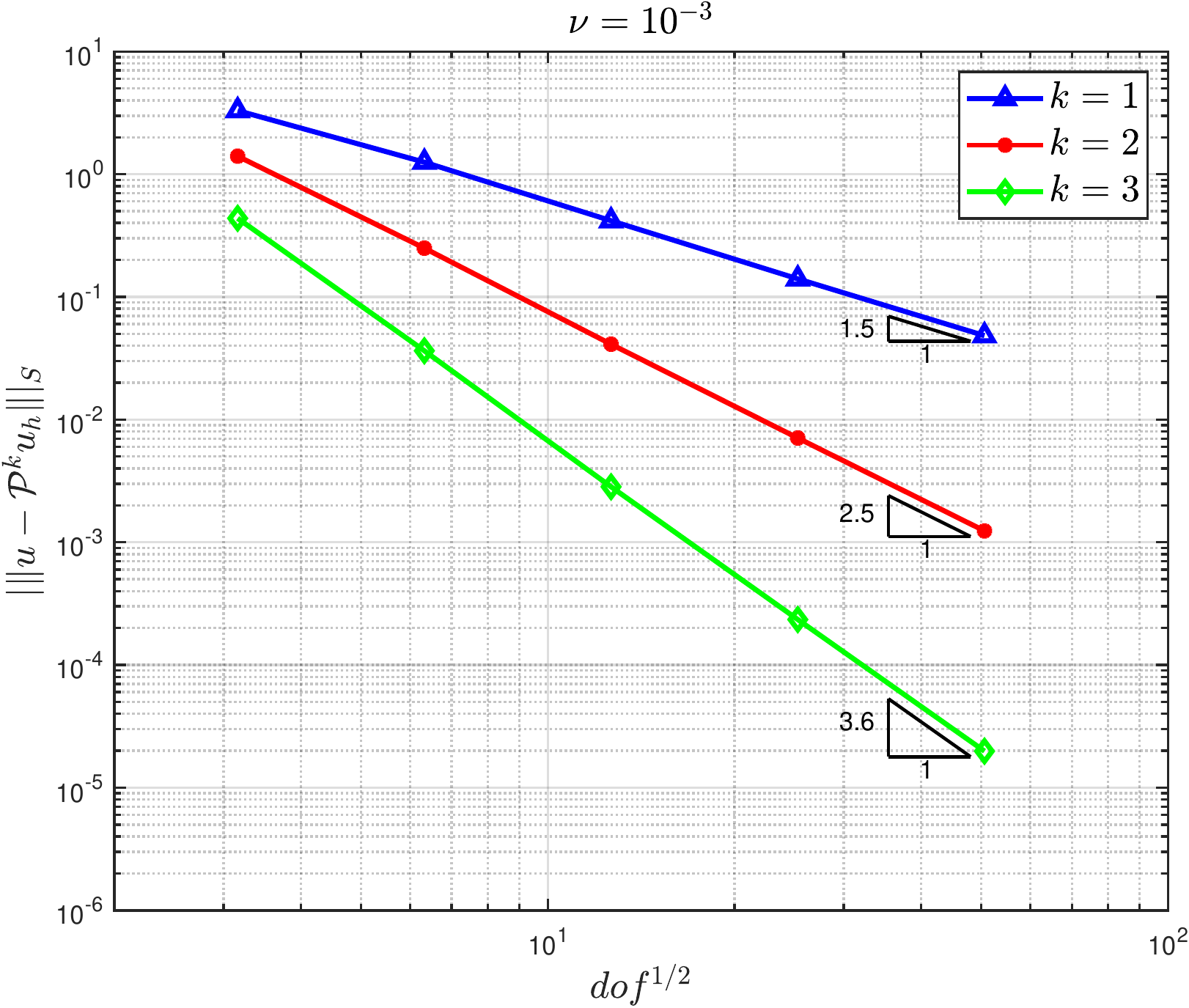}
}
\quad
\subfigure[]{
\includegraphics[width=7cm]{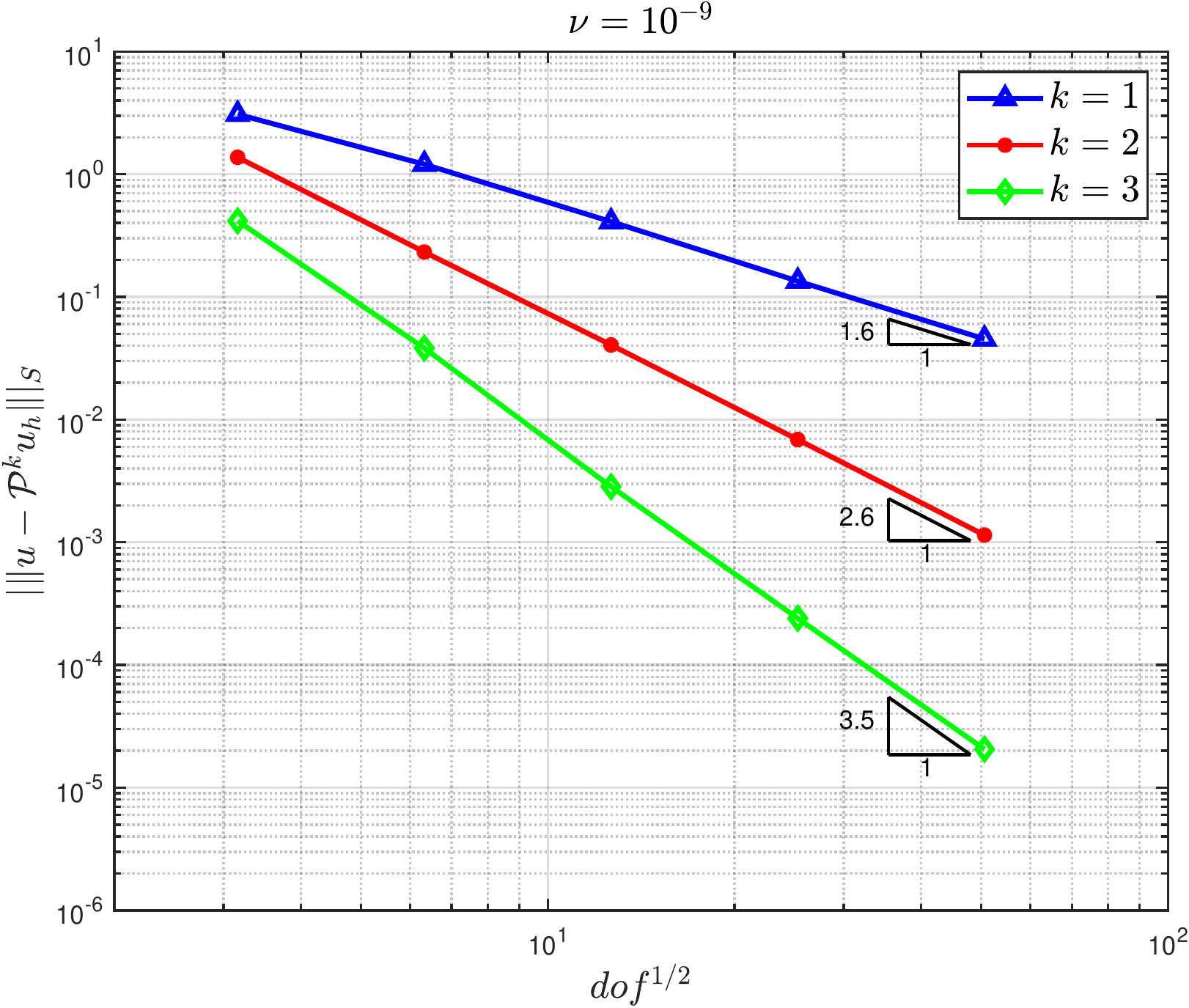}
}
\centering
\caption{\emph{Example 1.} Convergence rate results of $\normdg{u-\mathcal P^ku_h}$ and $\normdg{u-\mathcal P^ku_h}_S$ over regular polygonal meshes with a not divergence-free convection field: (a) $\normdg{u-\mathcal P^ku_h}$ with $\nu=1$; (b) $\normdg{u-\mathcal P^ku_h}$ with $\nu=10^{-3}$; (c) $\normdg{u-\mathcal P^ku_h}$ with $\nu=10^{-9}$; (d) $\normdg{u-\mathcal P^ku_h}_S$ with $\nu=10^{-3}$; (e) $\normdg{u-\mathcal P^ku_h}_S$ with $\nu=10^{-9}$.}\label{fig:example1:cr:l2dg:poly}
\end{figure}

Convergence plots of numerical errors for different reconstruction order $k$ are collected in Figure \ref{fig:example1:cr:l2dg:tri} versus uniform triangulations, and Figure \ref{fig:example1:cr:l2dg:poly} versus regular polygonal meshes. In particular, the square root of degrees of freedom ($dof^{1/2}$) replaces $h$ in Figure \ref{fig:example1:cr:l2dg:poly}.

It is observed that all optimal error orders are obtained in the norm $\normdg{\cdot}$ when $\nu=1,10^{-3},10^{-9}$, and in the norm $\normdg{\cdot}_{S}$ when $\nu=10^{-3},10^{-9}$. The above results match the statements in Theorem \ref{priori:dg} and Theorem \ref{priori:supg}. Thus, from the above results we observe that the reconstruction operator behaves well and provides good approximations when $\nabla\cdot\bm{b}\neq 0$.

\subsubsection{The convection field is divergence free}
The variable coefficients of \eqref{p1} are given by $\bm{b}(x,y)=[y,~x]^\top$, $c(x,y)=\exp(x+y)$, and we vary the diffusion coefficient $\nu=1,10^{-3},10^{-9}$. The source term $f$ is chosen so that the analytical solution of \eqref{p1}, with Dirichlet boundary condition \eqref{p2}, is given by
\begin{equation}\label{exact:solution:2}
u(x,y)=\sin(2\pi x)\sin(2\pi y)+x^5+y^5+1.
\end{equation}
Here the uniform triangulations and regular polygonal meshes are taken to show the convergence results in the $L^2$-norm $\|\cdot\|$. In particular, for regular polygonal meshes we use the total number of the elements to be 10, 40, 160, 640 and 2560 to simulate the effect of uniform refinement.

\begin{figure}[htbp]
\centering
\subfigure[]{
\includegraphics[width=7cm]{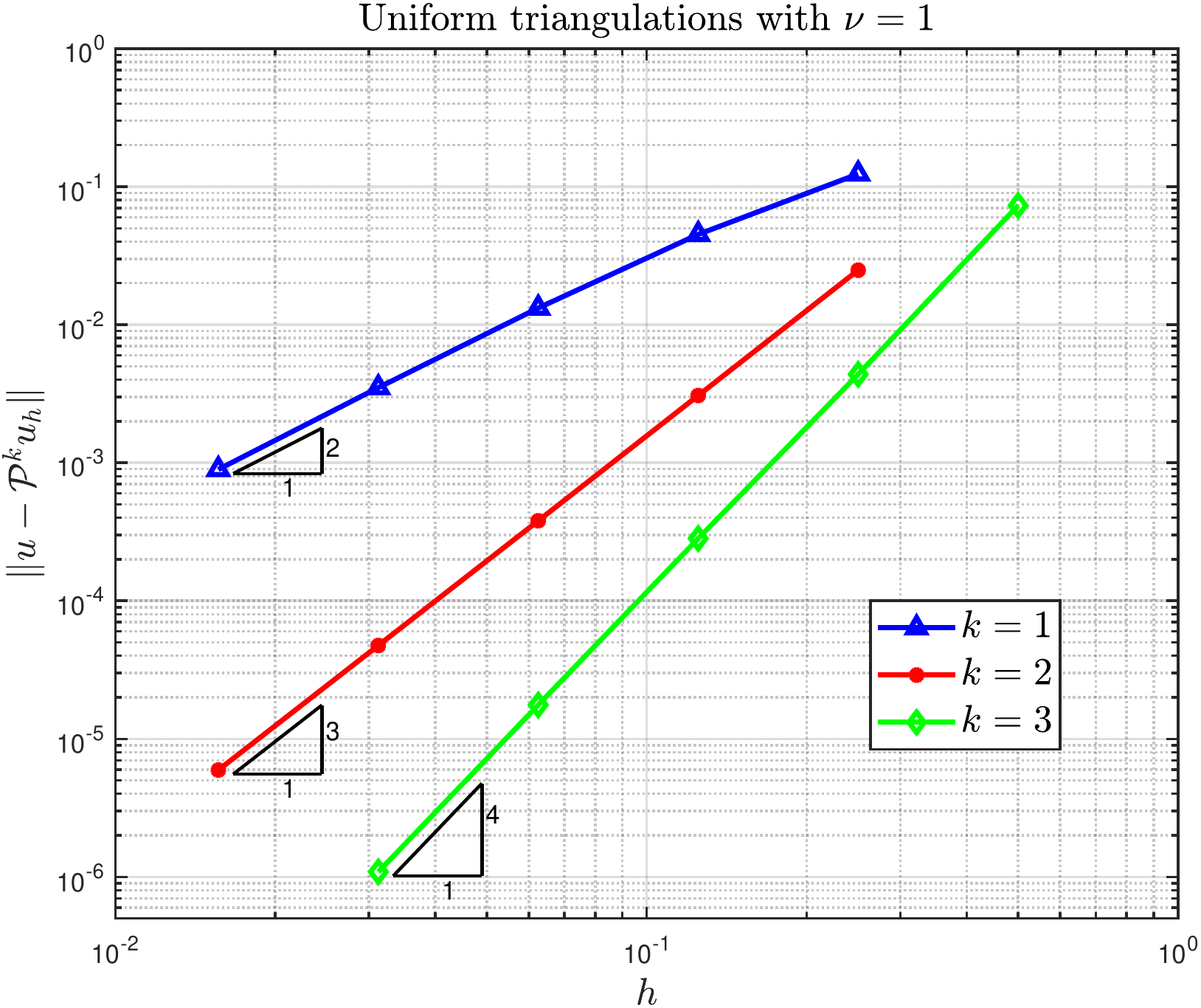}
}
\quad
\subfigure[]{
\includegraphics[width=7cm]{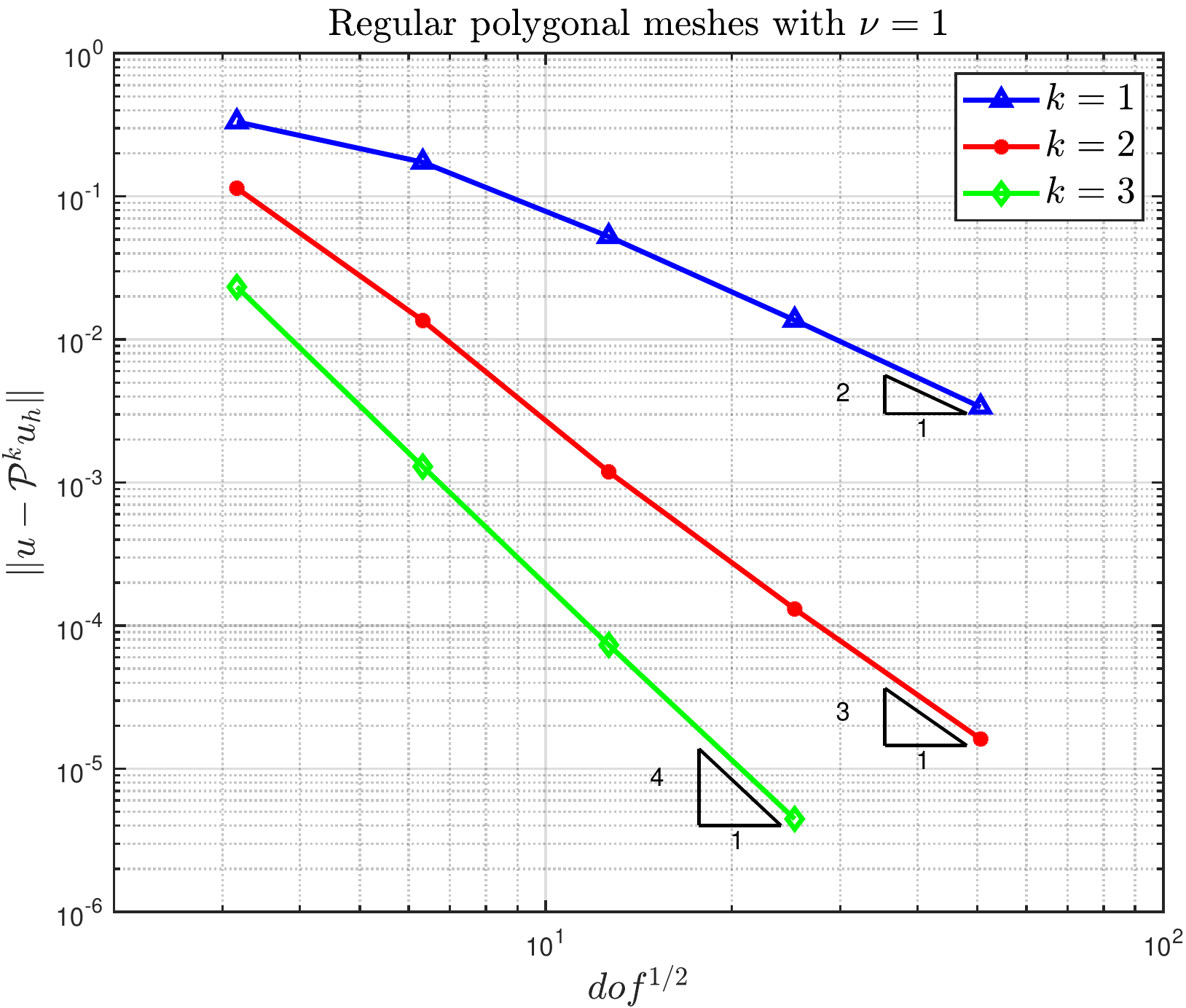}
}

\subfigure[]{
\includegraphics[width=7cm]{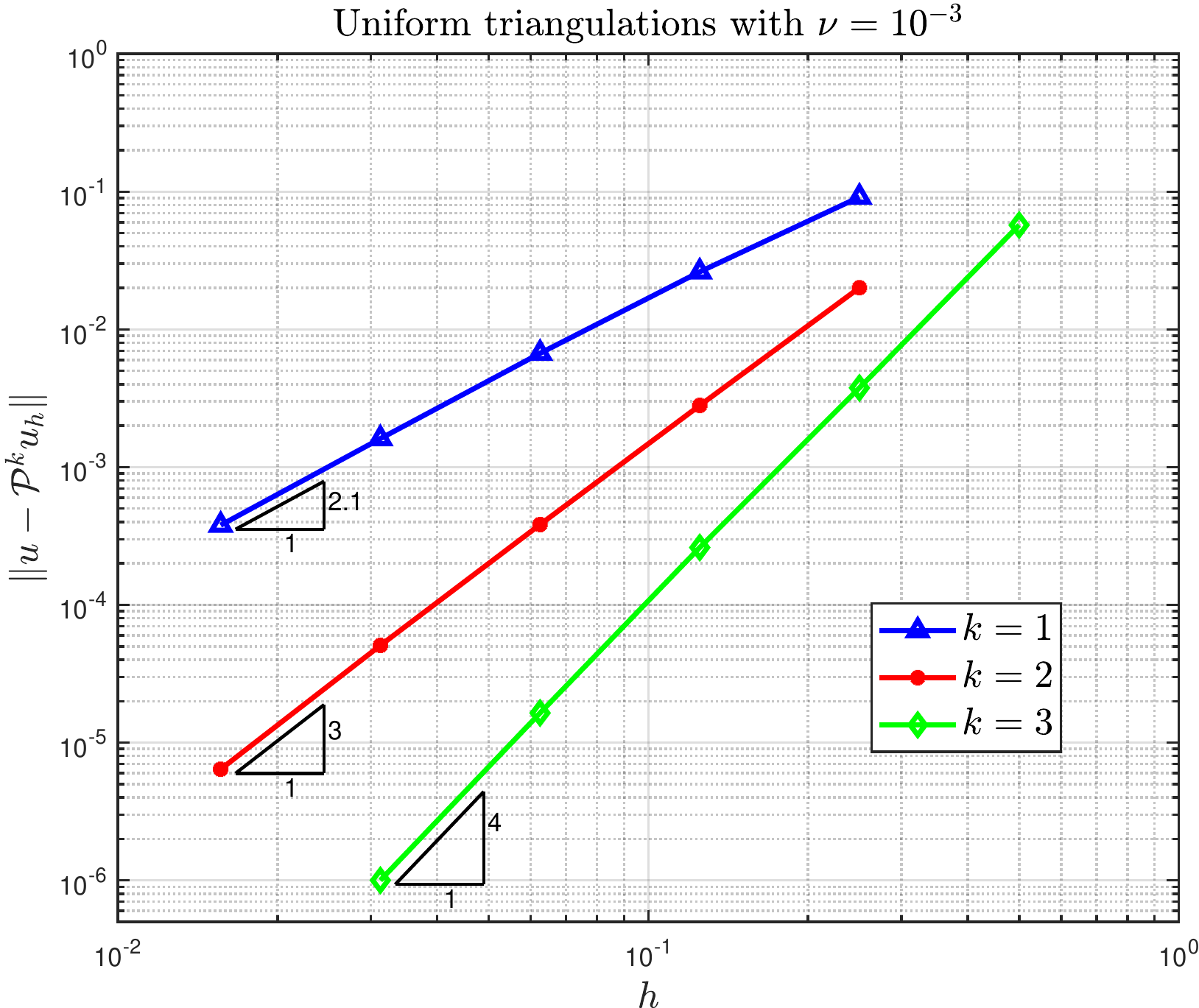}
}
\quad
\subfigure[]{
\includegraphics[width=7cm]{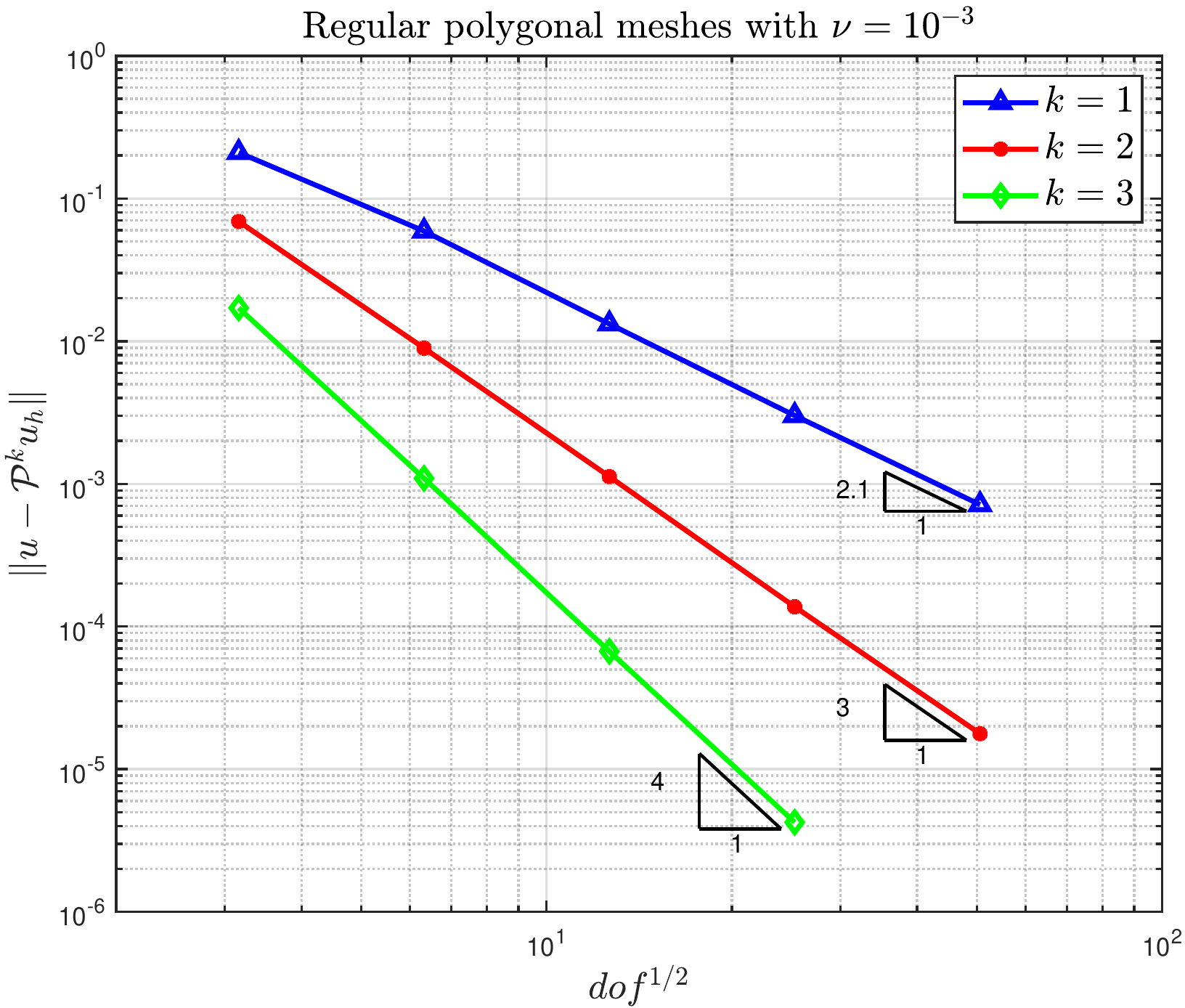}
}

\subfigure[]{
\includegraphics[width=7cm]{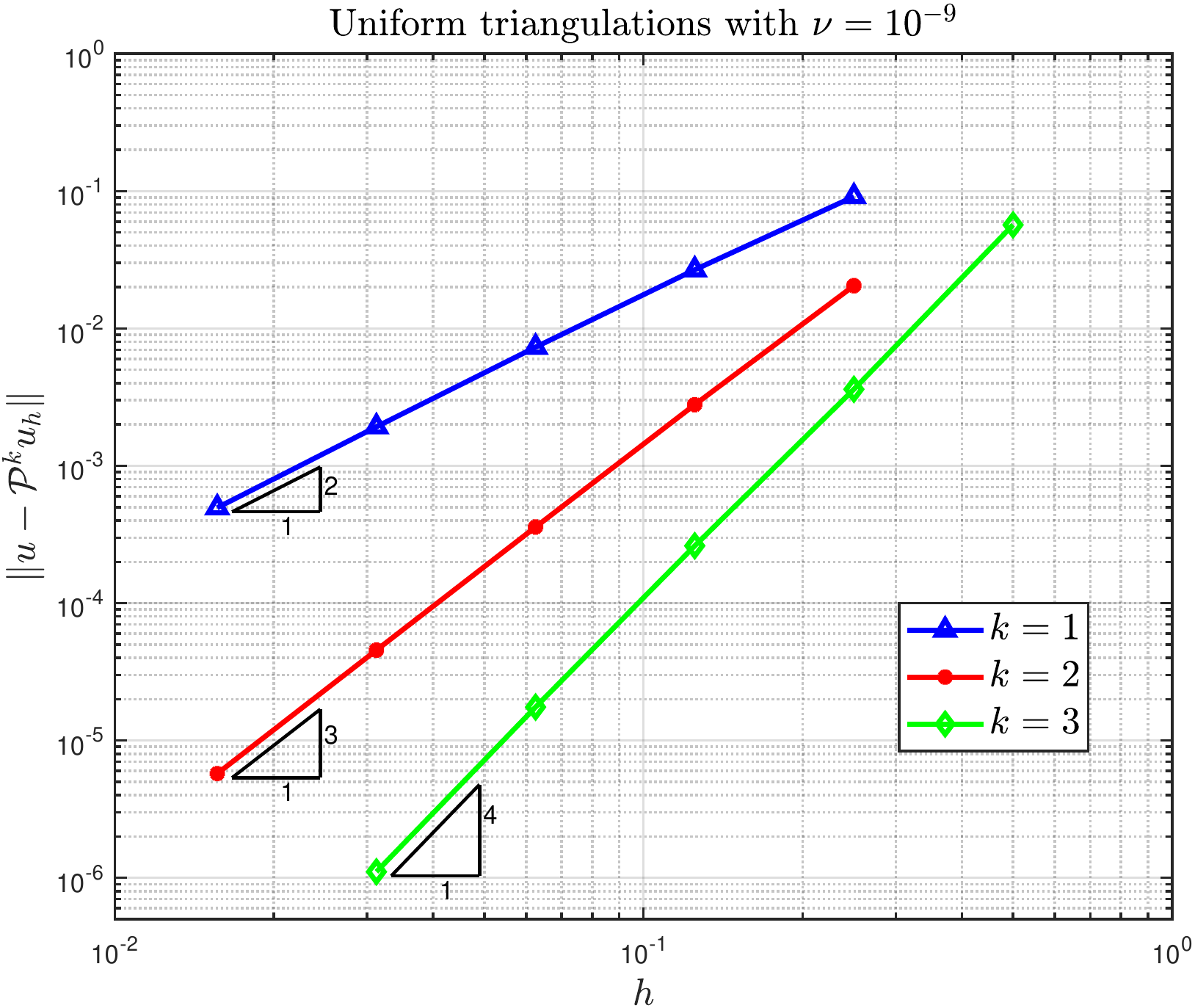}
}
\quad
\subfigure[]{
\includegraphics[width=7cm]{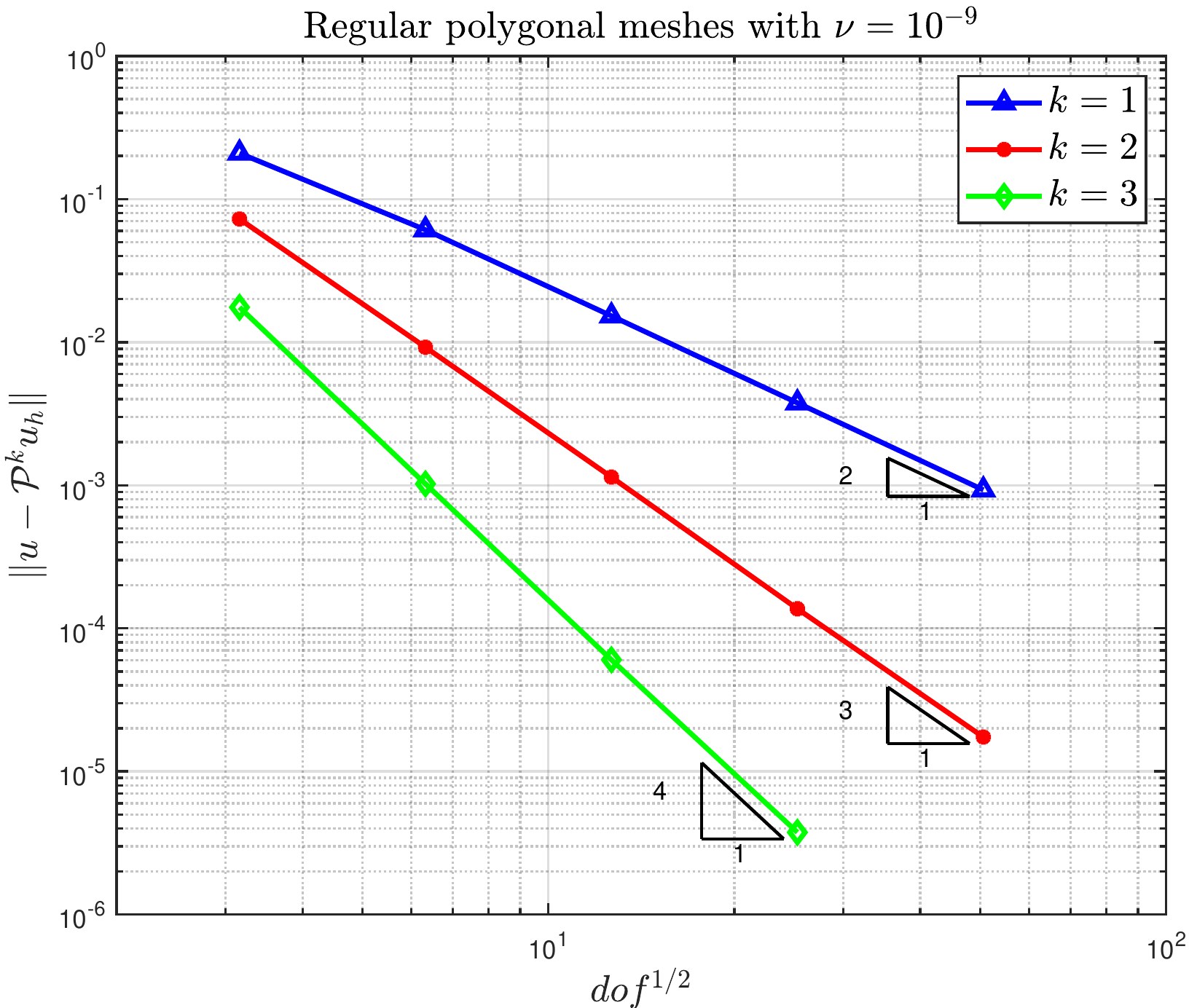}
}
\centering
\caption{\emph{Example 1.} Convergence rate results in the $L^2$-norm $\|\cdot\|$  with a divergence-free convection field for different values of $\nu$ over uniform triangulations and regular polygonal meshes: the left column results are over uniform triangulations, and the right column results are over regular polygonal meshes. }\label{fig:example1:l2:tm}
\end{figure}

Despite in Theorem \ref{priori:l2} only the theoretical sub-optimal order can be achieved in a convection-dominated regime, it is observed that in Figure \ref{fig:example1:l2:tm} the optimal error orders are achieved in the $L^2$-norm $\|\cdot\|$.

\subsection{Example 2: Interior layer with continuous boundary conditions.}
This example is taken from \cite{Lin20181482}. In this example, we examine the performance of the reconstruction operator in the occurrence of an interior layer. Let $\bm{b}=[1,0]^\top$ and $c=1$, and the source term $f$ is chosen so that the exact solution of \eqref{p1}, with Dirichlet boundary conditions, is given by
\[
u(x,y)=\frac{1}{2}x(1-x)y(1-y)
\bigg(1-\tanh\frac{l_1-x}{l_2}\bigg),
\]
where the two parameters $l_1$ and $l_2$ control the location and thickness of the interior layer, respectively.

\begin{figure}[htbp]
\centering
\subfigure[]{
\includegraphics[width=7.5cm]{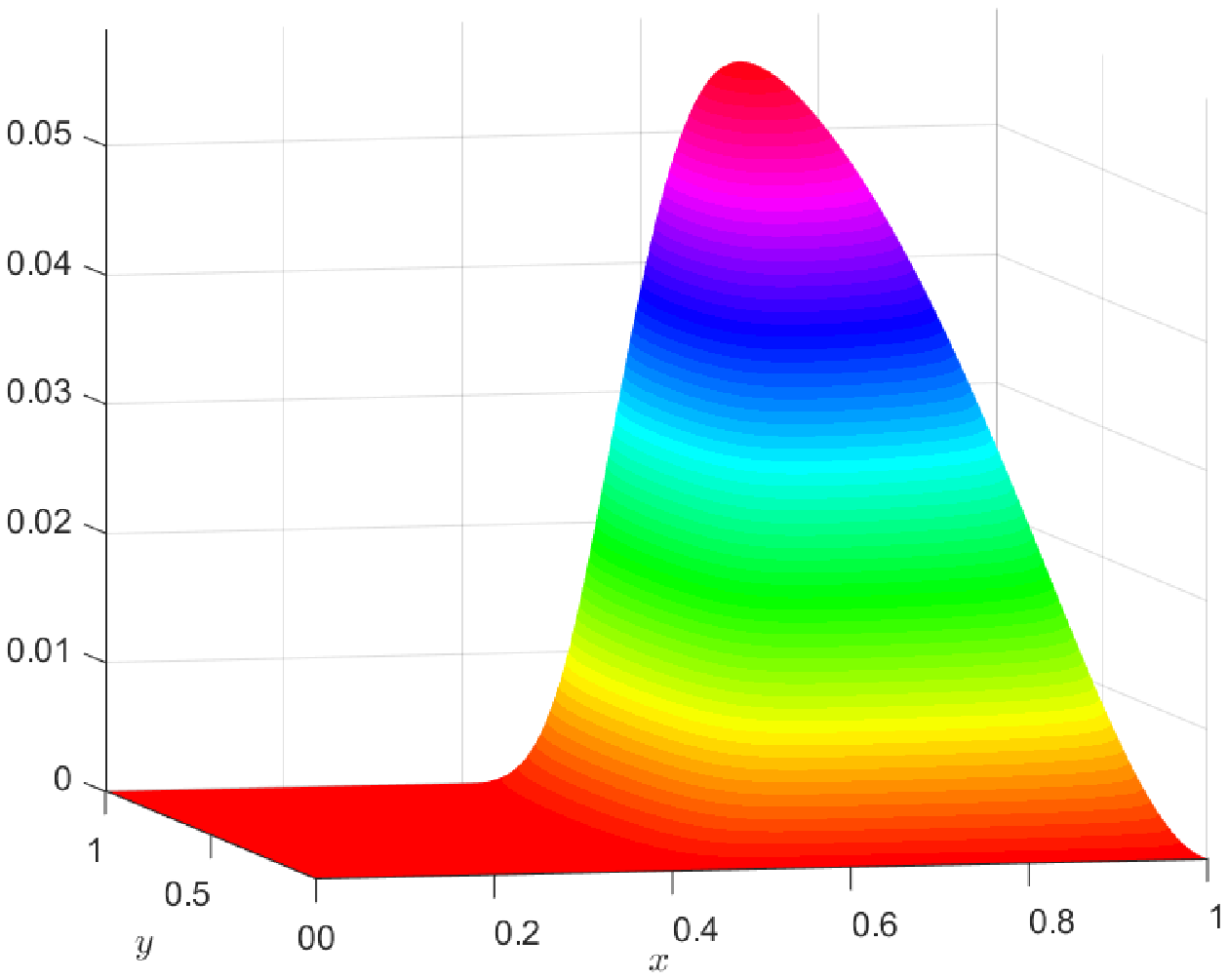}
}
\quad
\subfigure[]{
\includegraphics[width=7.5cm]{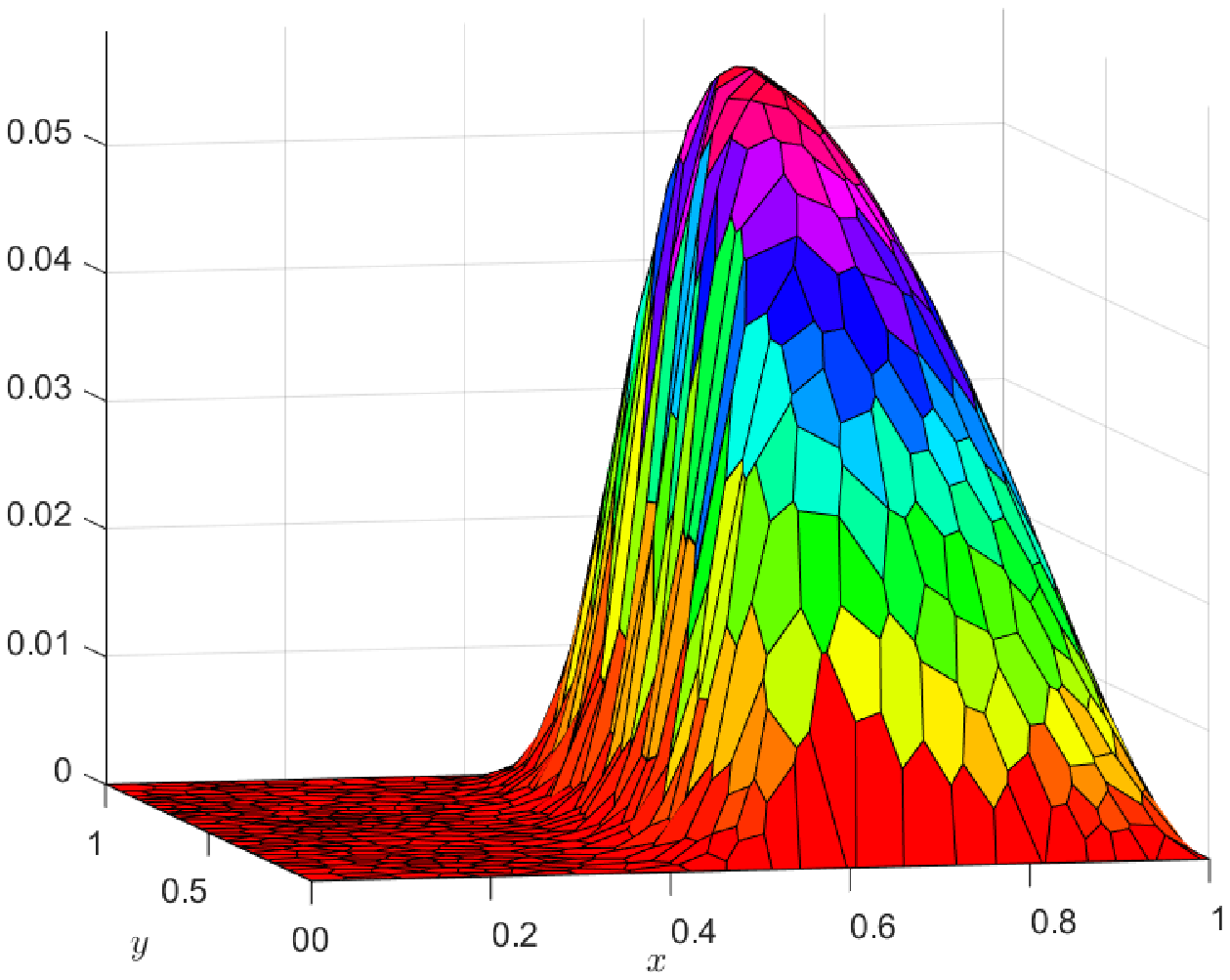}
}

\subfigure[]{
\includegraphics[width=7.5cm]{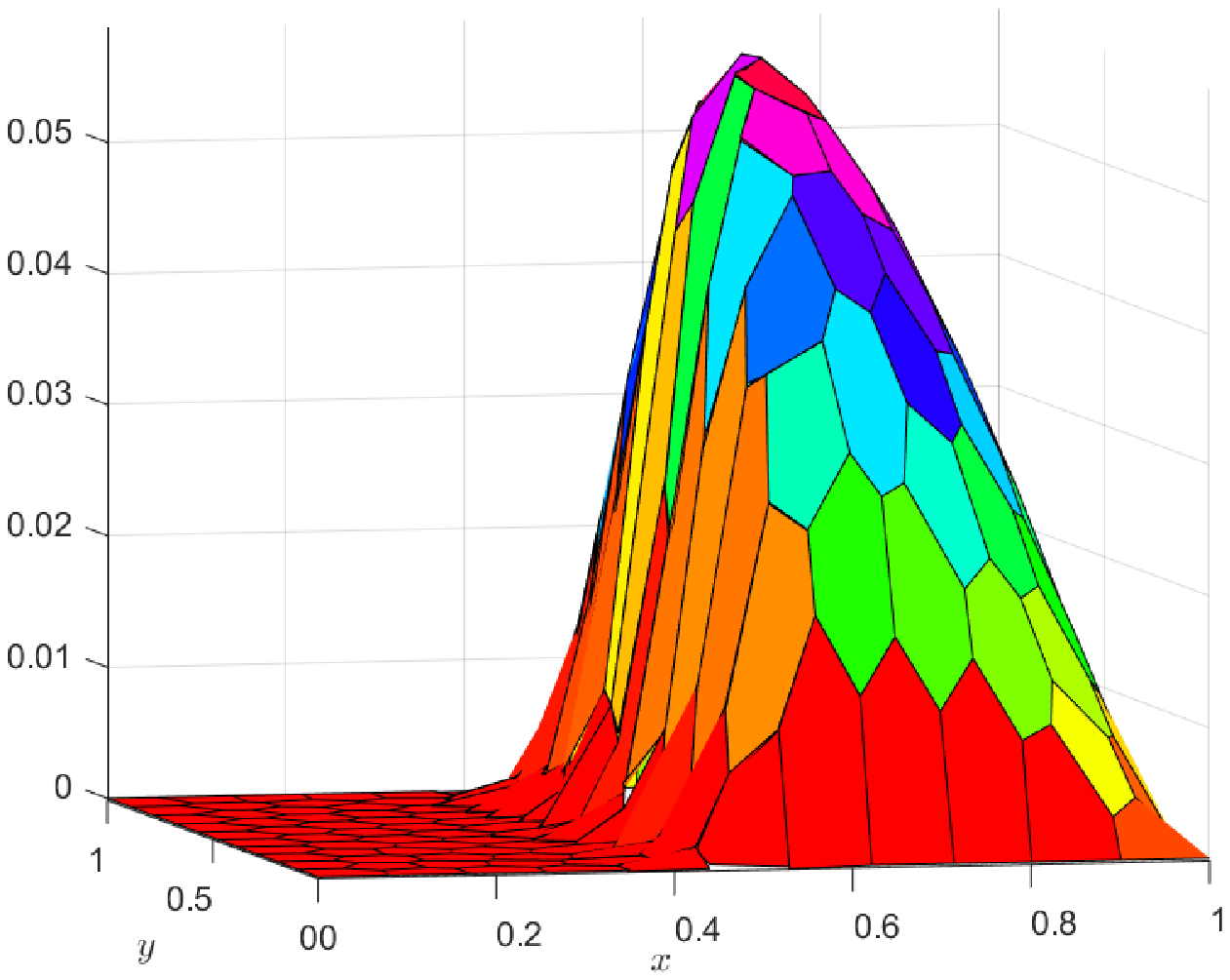}
}
\quad
\subfigure[]{
\includegraphics[width=7.5cm]{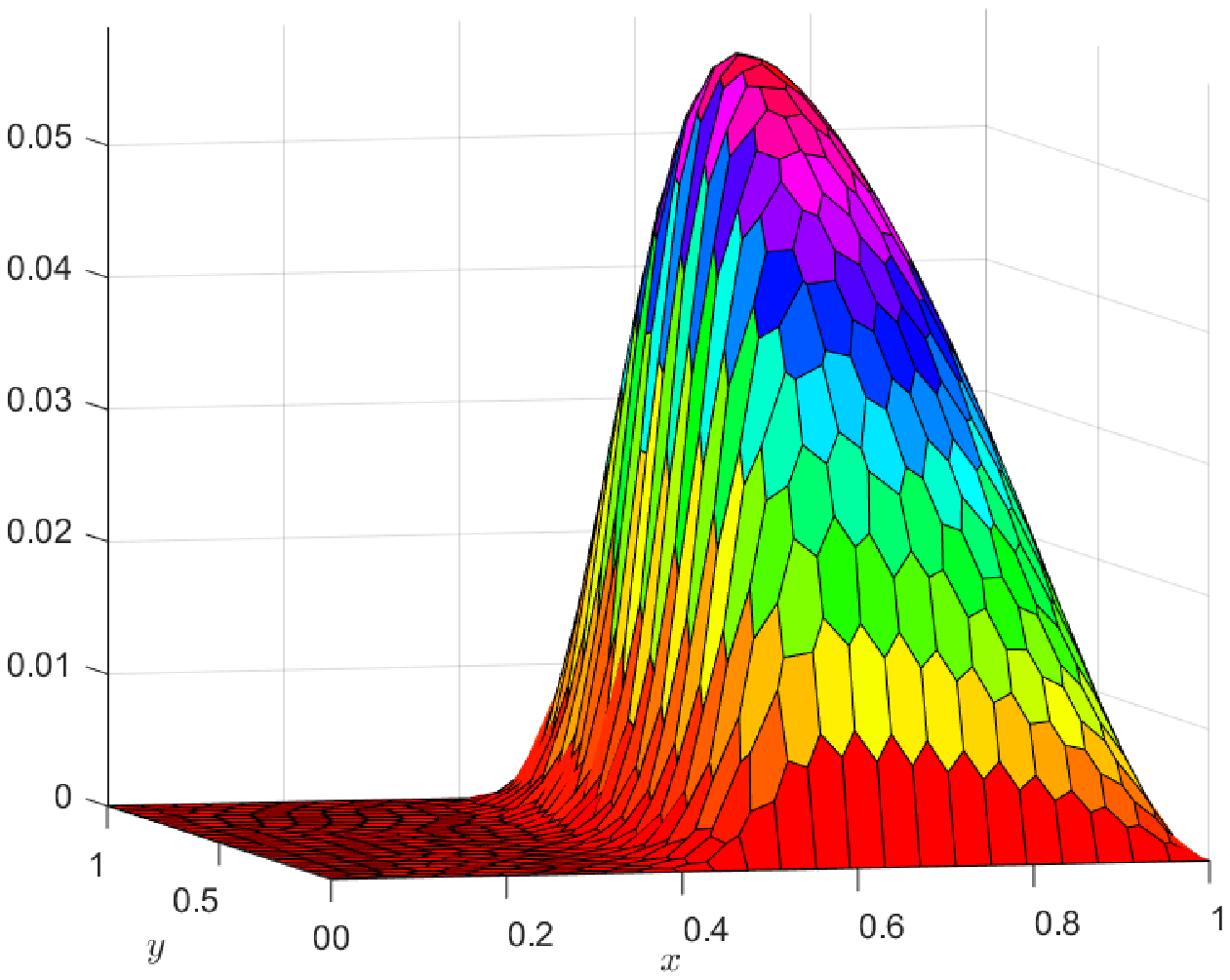}
}
\centering
\caption{\emph{Example 2.} The exact solution and numerical solutions with $k=3$, $\nu=10^{-9}$, $l_1=0.5$ and $l_2=0.05$ over regular polygonal meshes and general Voronoi meshes: (a) the exact solution; (b) the numerical solution over the general Voronoi mesh with 640 elements; (c) the numerical solution over the regular polygonal mesh with 160 elements; (d) the numerical solution over the regular polygonal mesh with 640 elements.}\label{fig:example2:solution}
\end{figure}

\begin{figure}[htbp]
\centering
\subfigure[]{
\includegraphics[width=7cm]{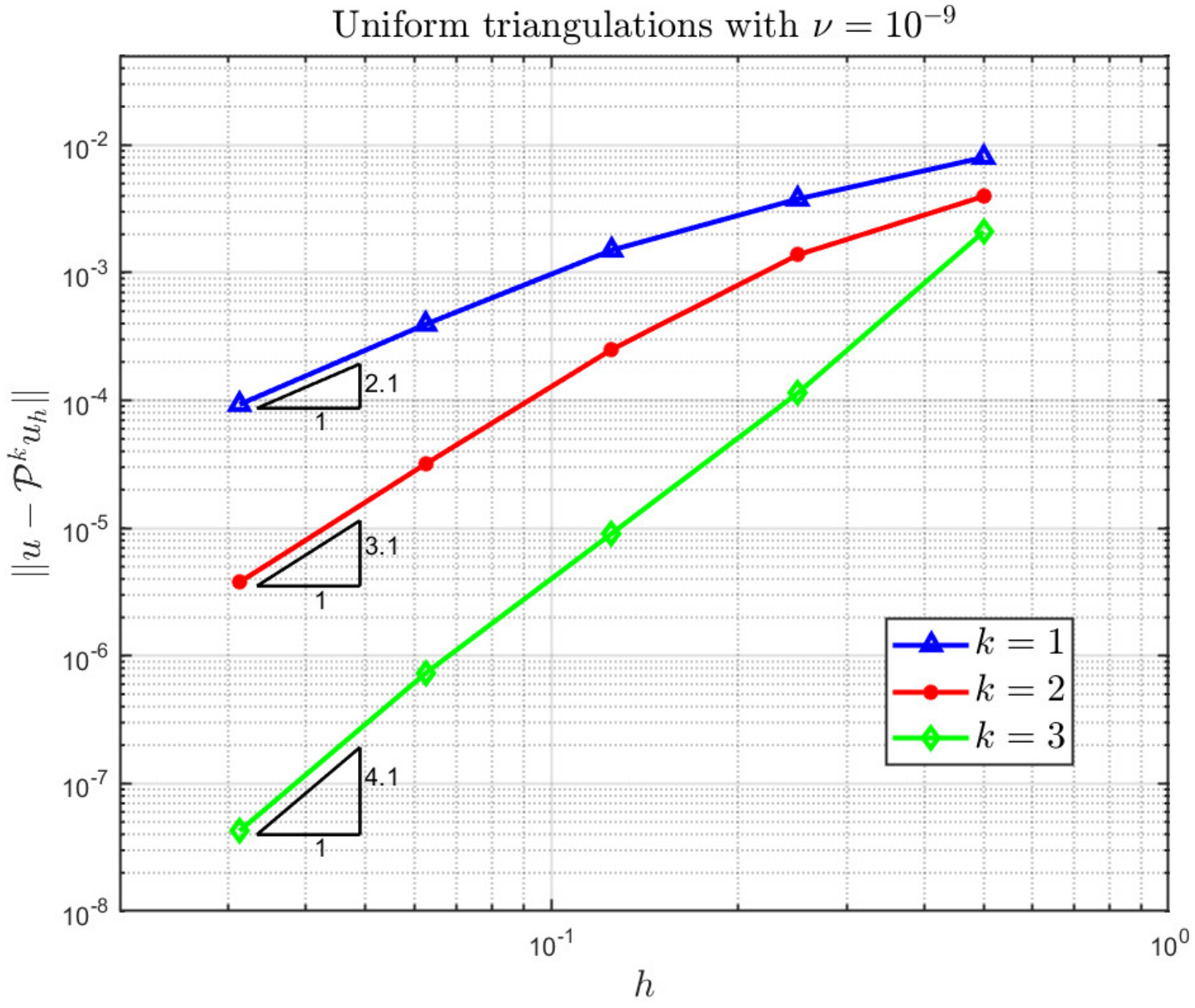}
}
\quad
\subfigure[]{
\includegraphics[width=7cm]{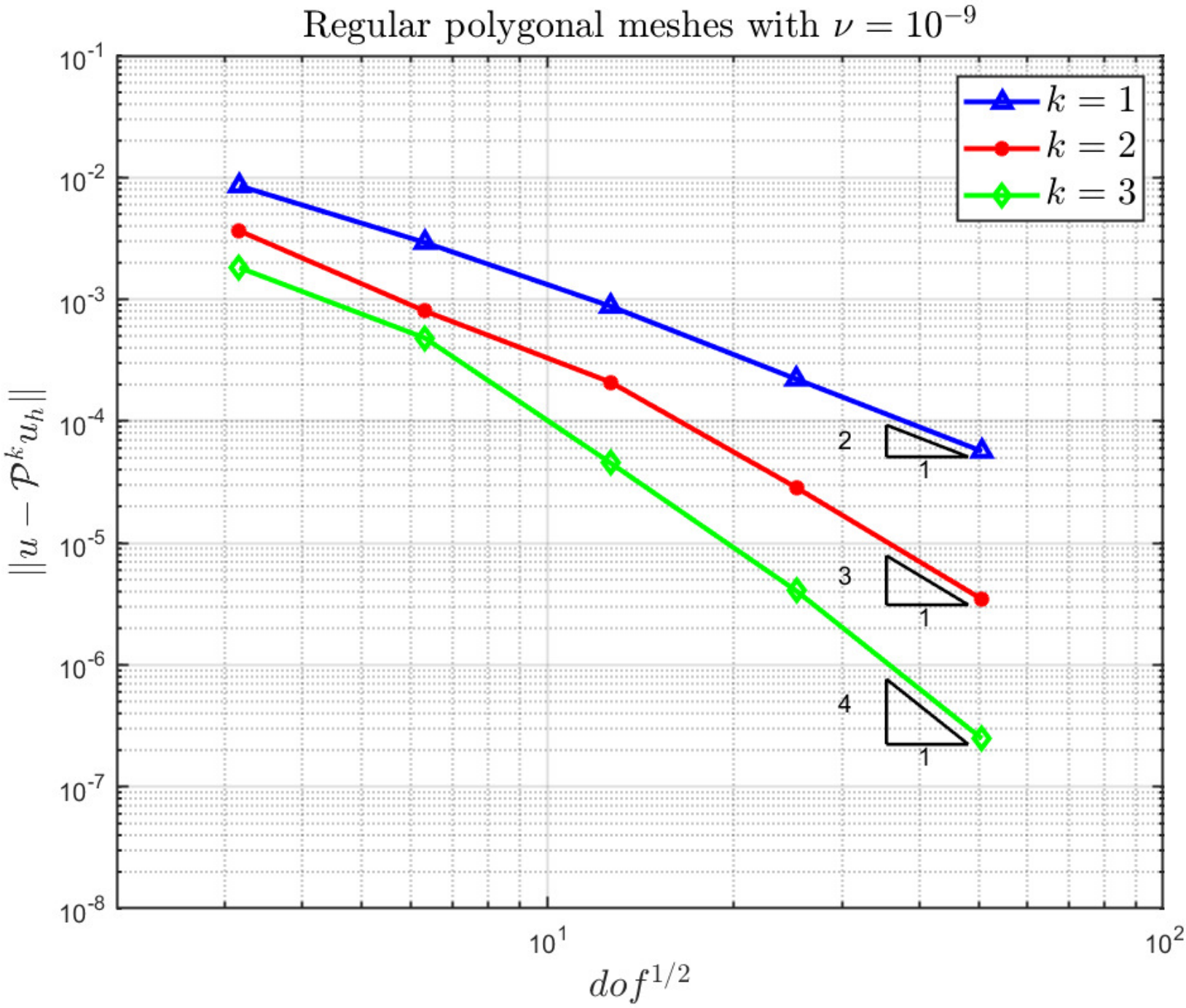}
}

\subfigure[]{
\includegraphics[width=7cm]{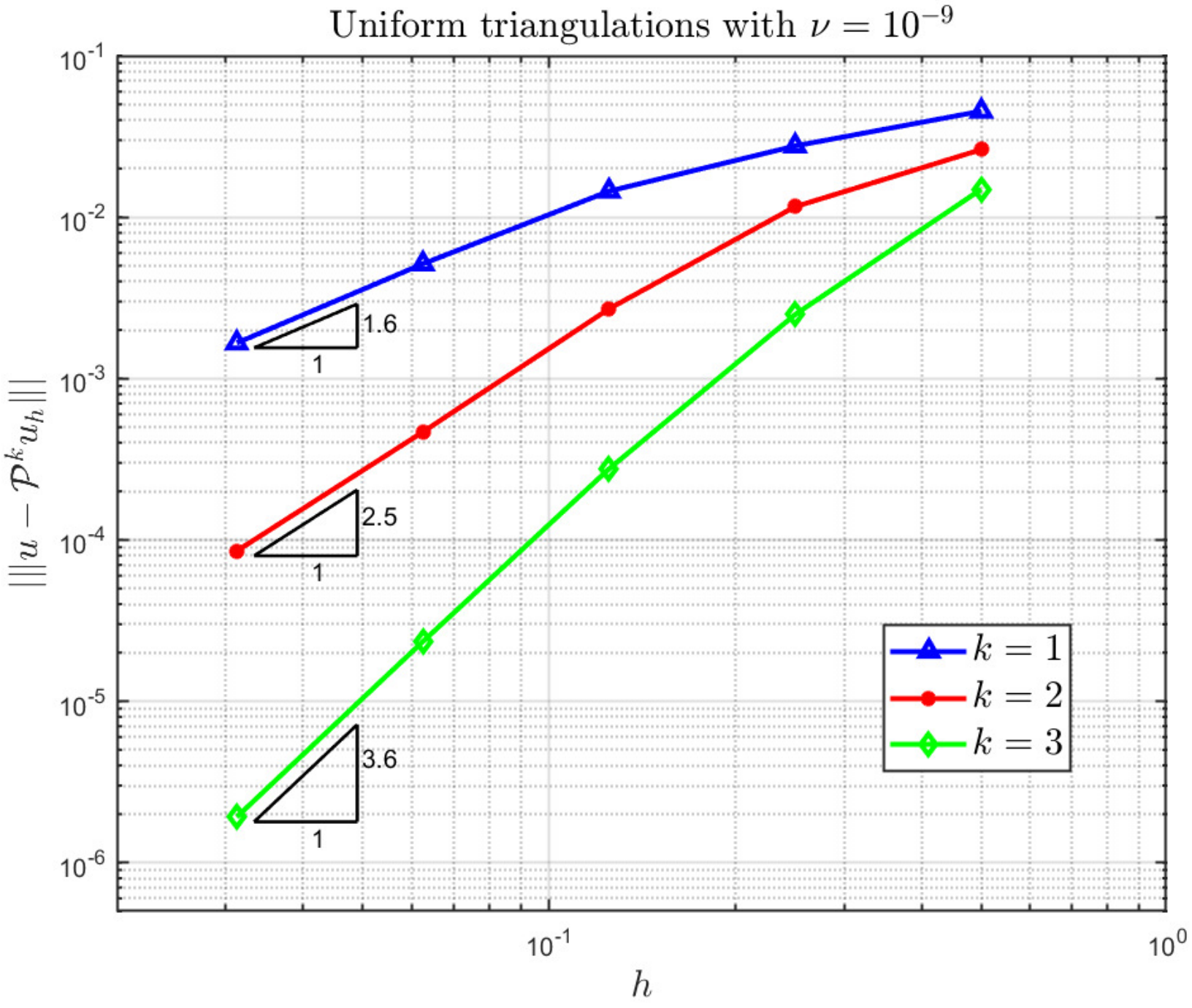}
}
\quad
\subfigure[]{
\includegraphics[width=7cm]{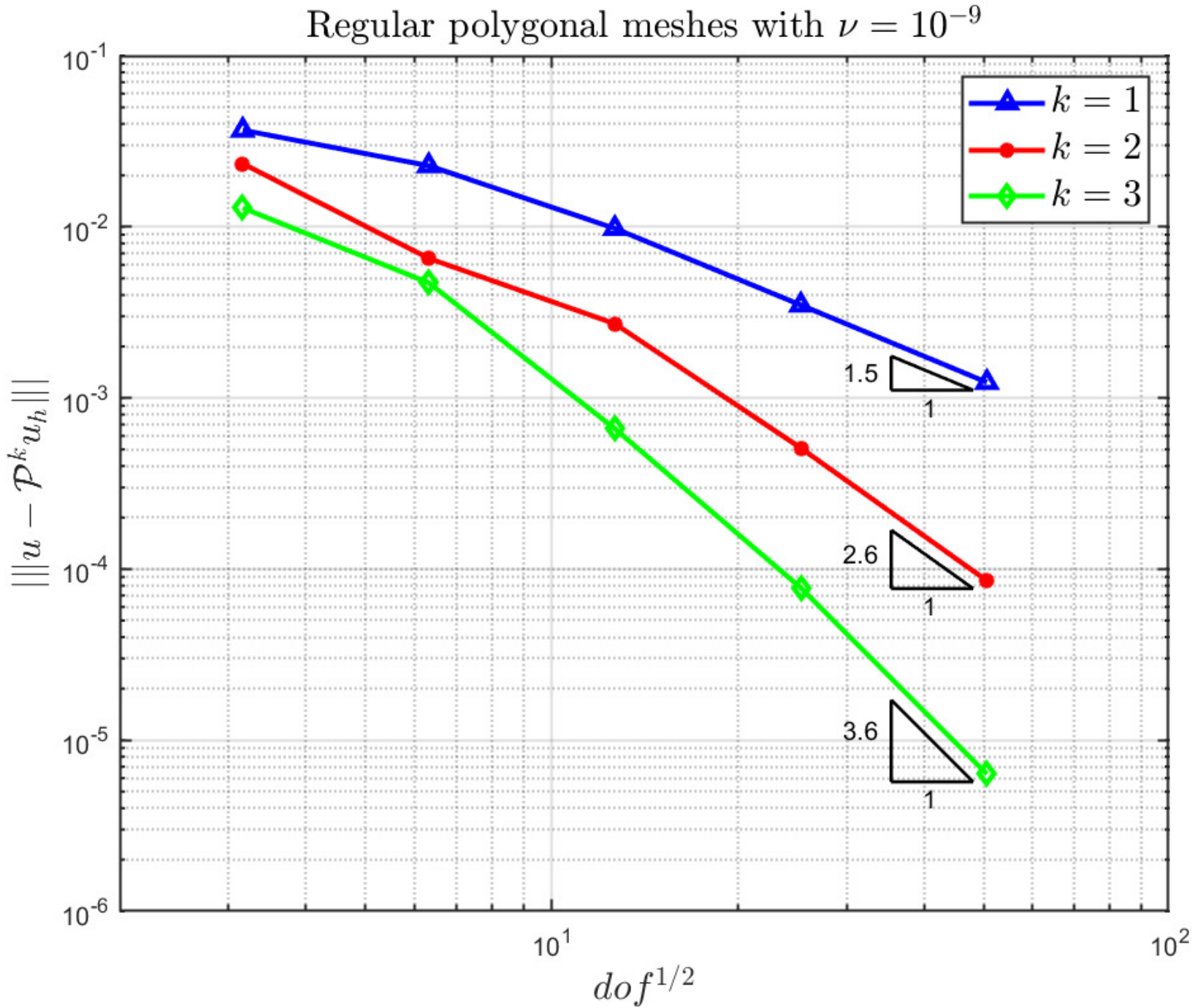}
}

\subfigure[]{
\includegraphics[width=7cm]{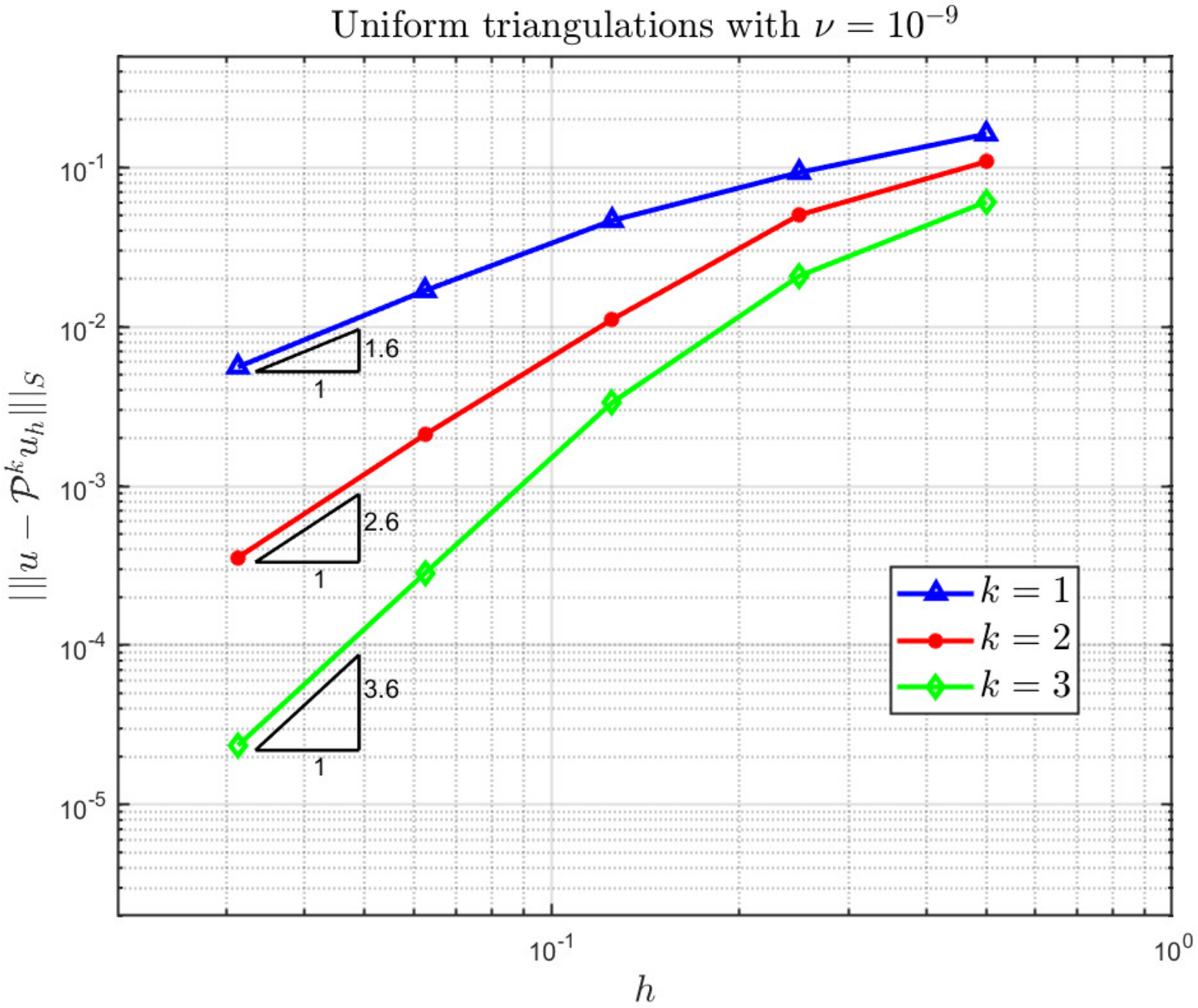}
}
\quad
\subfigure[]{
\includegraphics[width=7cm]{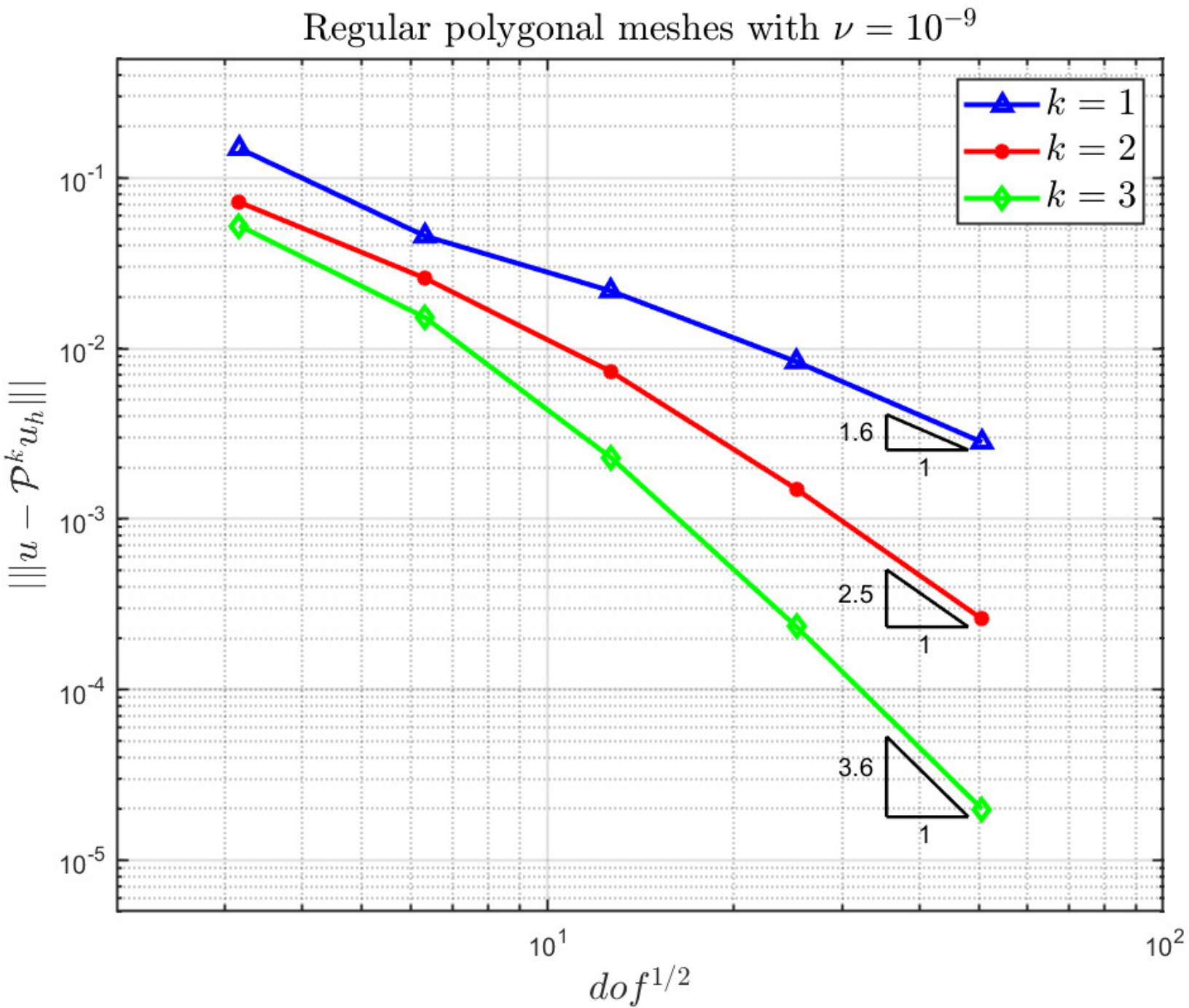}
}
\centering
\caption{\emph{Example 2.} Convergence rate results of $\|\mathcal P^ku_h-u\|$, $\normdg{\mathcal P^ku_h-u}$, and $\normdg{\mathcal P^ku_h-u}_S$ with $\nu=10^{-9}$, $l_1=0.5$ and $l_2=0.05$ over uniform triangulations and regular polygonal meshes: the left column results are over uniform triangulations, and the right column results are over regular polygonal meshes. }\label{fig:example2:cr:l2dgsupg}
\end{figure}

In Figure \ref{fig:example2:solution} we present the numerical solutions with $k=3$, $\nu=10^{-9}$, $l_1=0.5$ and $l_2=0.05$ over general Voronoi meshes and regular polygonal meshes, respectively. It is observed that the interior layer is accurately captured. In Figure \ref{fig:example2:cr:l2dgsupg}, the convergence results are also provided versus uniform triangulations and regular polygonal meshes in the norm $\|\cdot\|$, $\normdg{\cdot}$ and $\normdg{\cdot}_S$ with $\nu=10^{-9}$. We can see that all optimal error orders can be achieved.

\begin{figure}[htbp]
\centering
\subfigure[]{
\includegraphics[width=7.5cm]{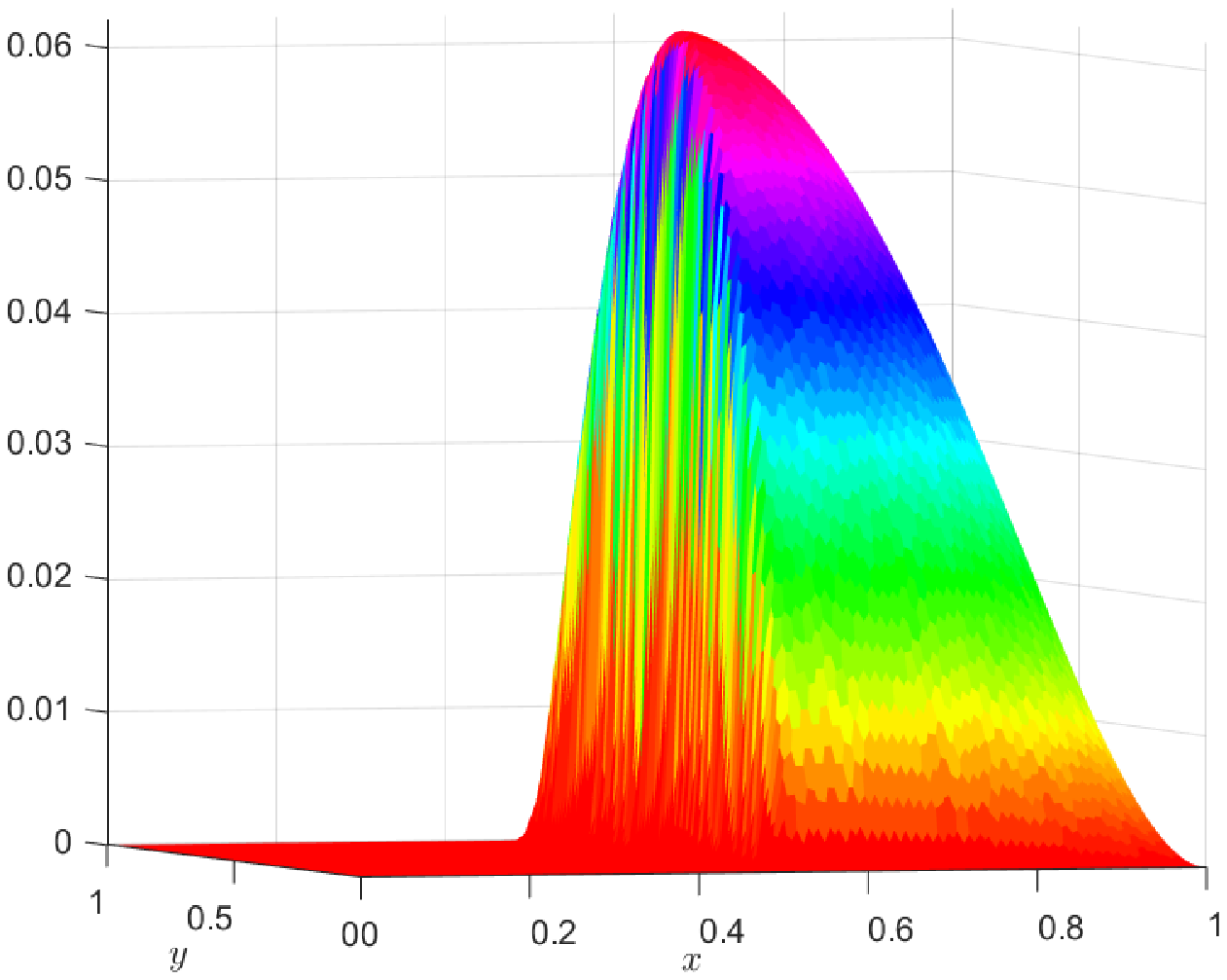}
}
\quad
\subfigure[]{
\includegraphics[width=7.5cm]{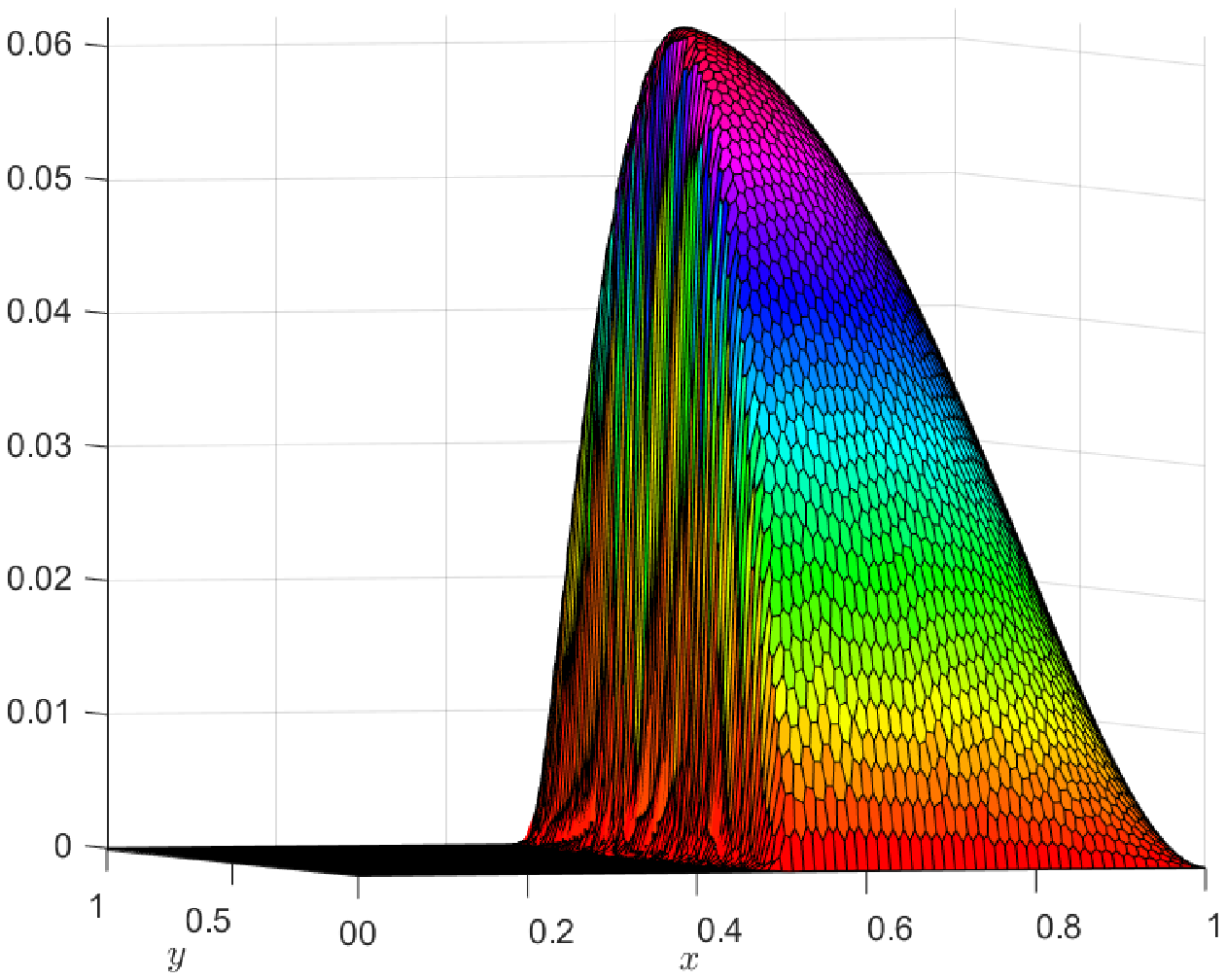}
}

\subfigure[]{
\includegraphics[width=7.5cm]{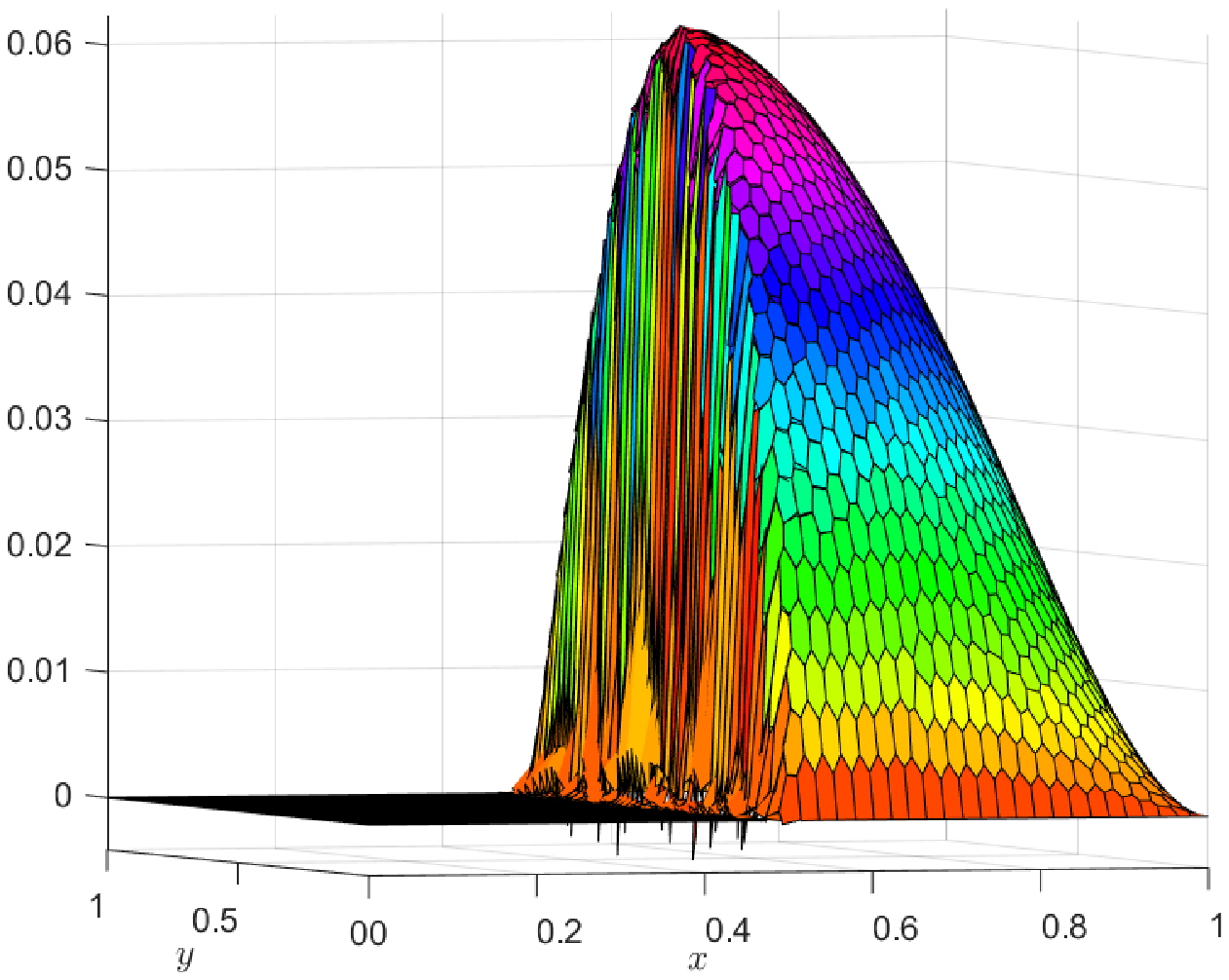}
}
\quad
\subfigure[]{
\includegraphics[width=7.5cm]{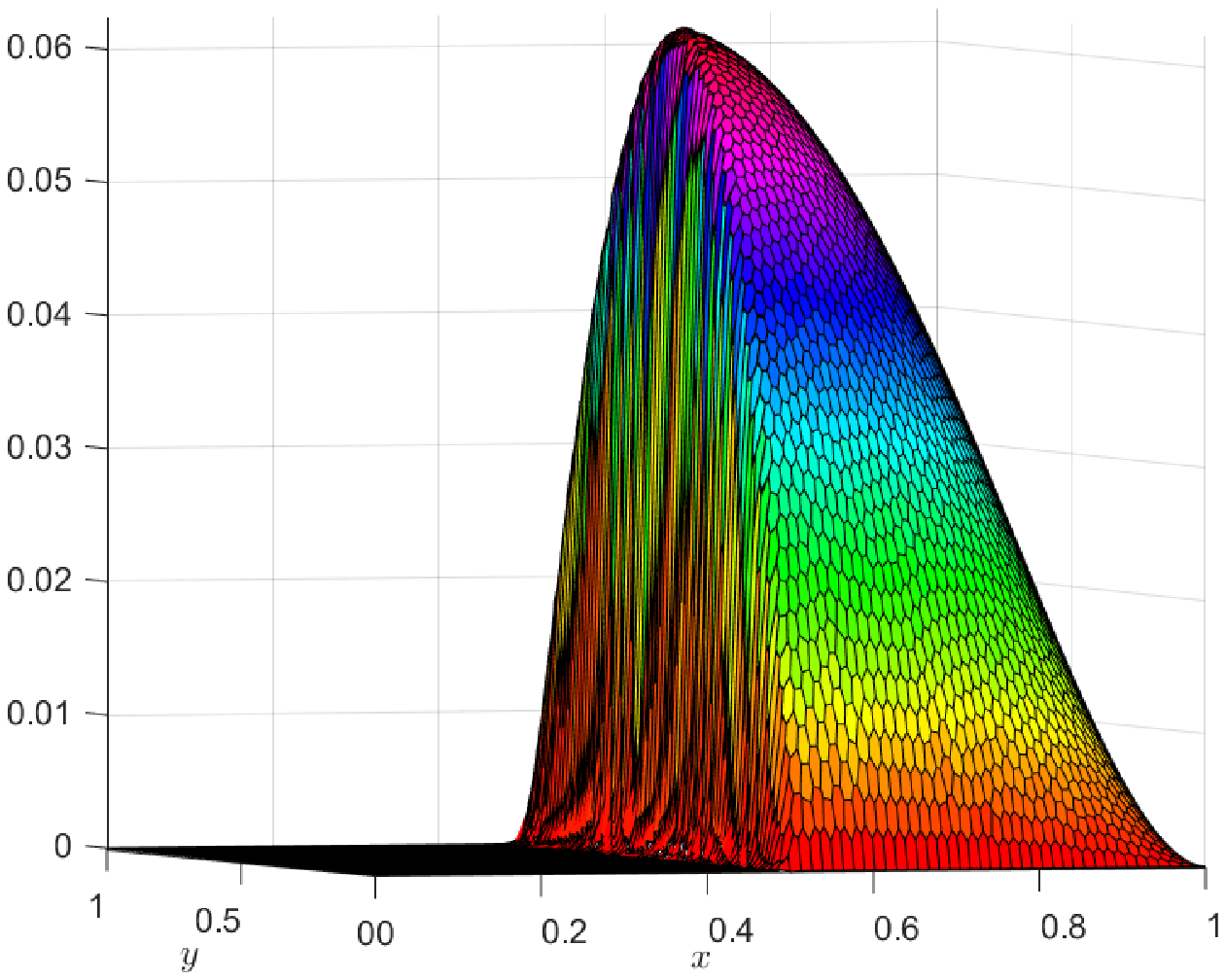}
}
\centering
\caption{\emph{Example 2.} The exact solution and numerical solutions with $k=2$, $l_1=0.5$ and $l_2=0.01$ over regular polygonal meshes: (a) the exact solution; (b) the numerical solution with 640 elements, $\nu=10^{-3}$; (c) the numerical solution with 2560 elements, $\nu=10^{-9}$; (d) the numerical solution with 6520 elements, $\nu=10^{-9}$.}\label{fig:example2:singular:solution}
\end{figure}

\begin{figure}[htbp]
\centering
\subfigure[]{
\includegraphics[width=7.5cm]{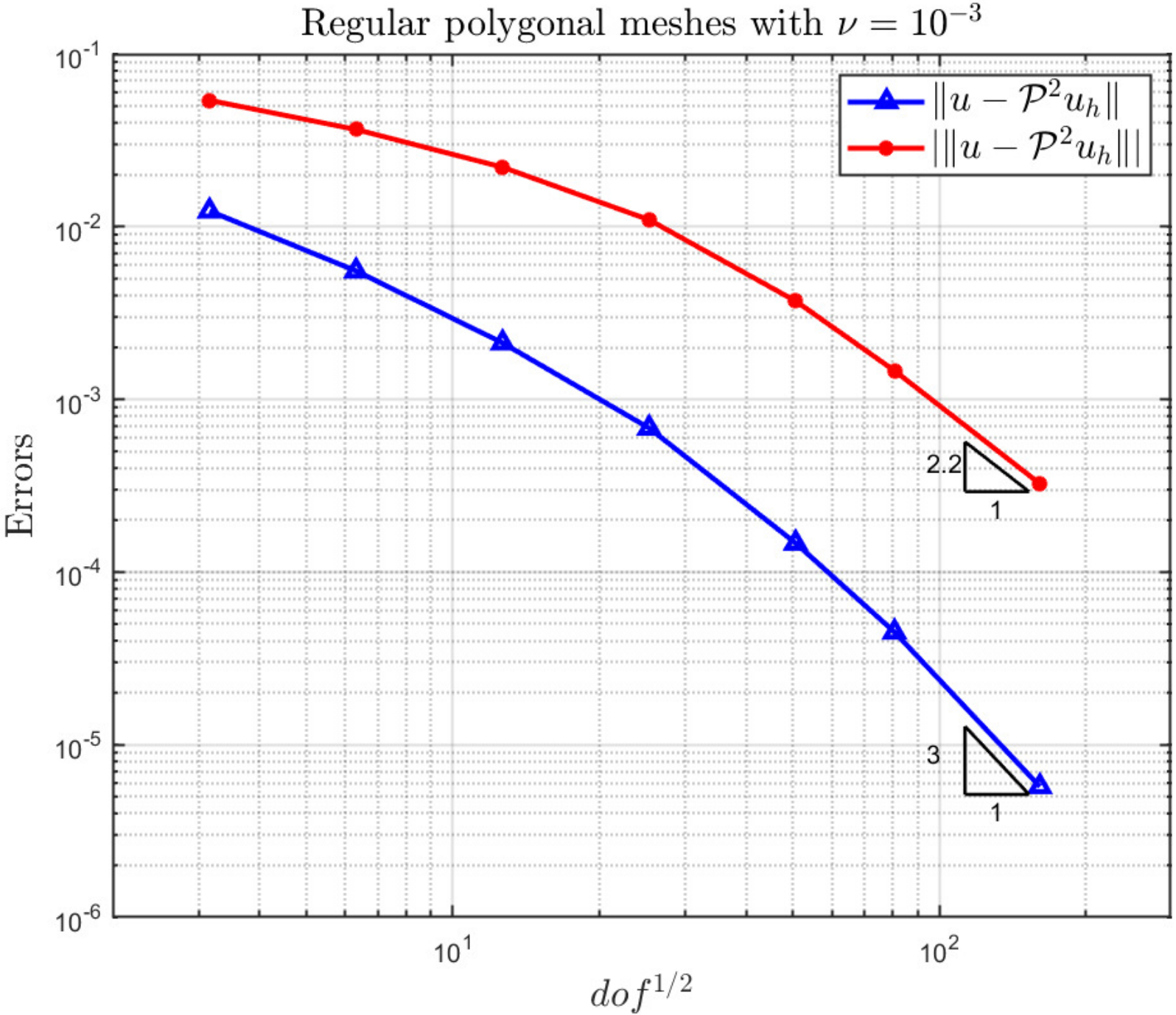}
}
\quad
\subfigure[]{
\includegraphics[width=7.5cm]{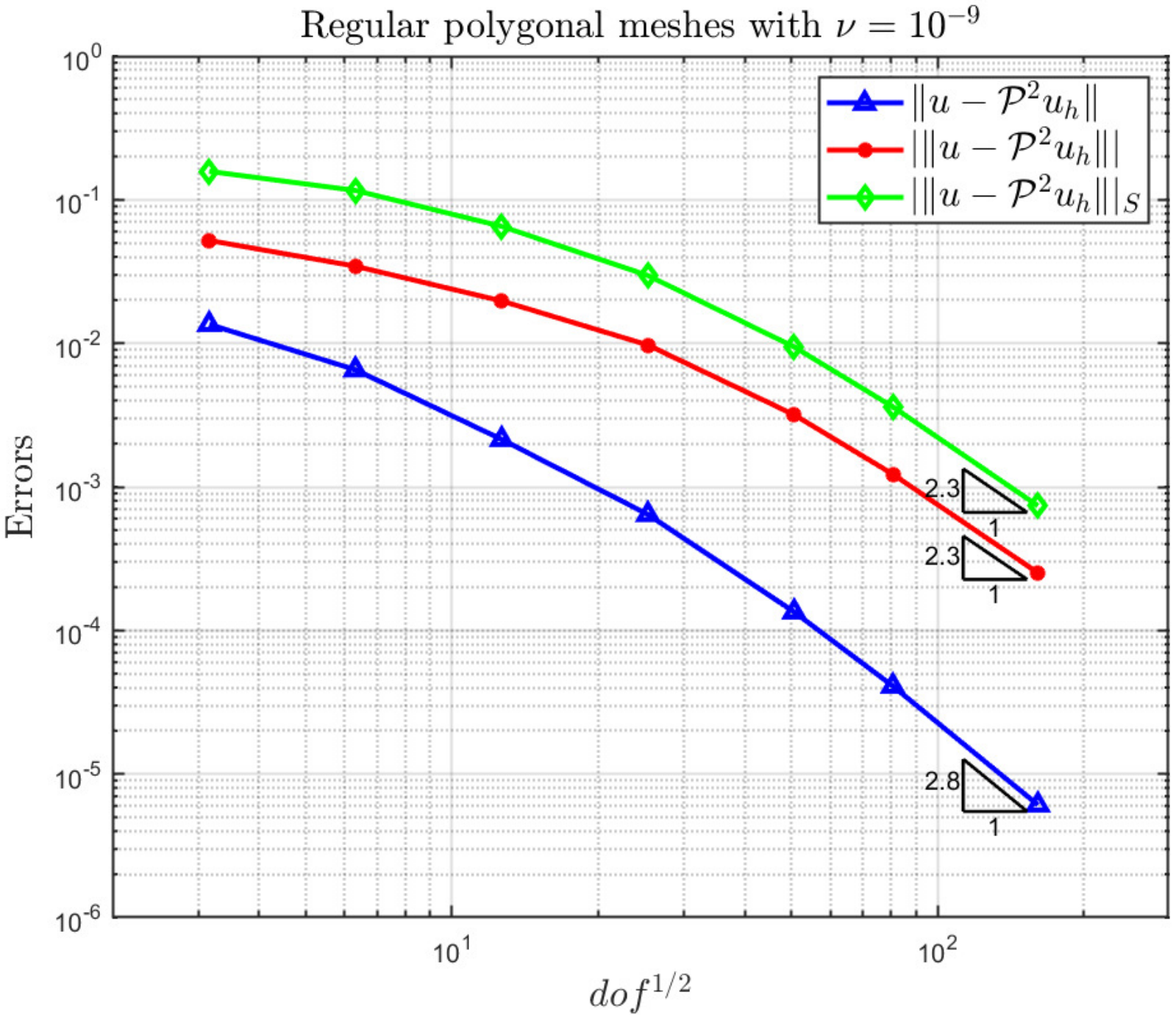}
}
\centering
\caption{\emph{Example 2.} Convergence rate results in different norms with $k=2$, $l_1=0.5$ and $l_2=0.01$ over regular polygonal meshes: (a) $\nu=10^{-3}$; (b) $\nu=10^{-9}$.}\label{fig:example2:singularcr:l2dgsupg}
\end{figure}

\begin{figure}[htbp]
\centering
\subfigure[]{
\includegraphics[width=7.5cm]{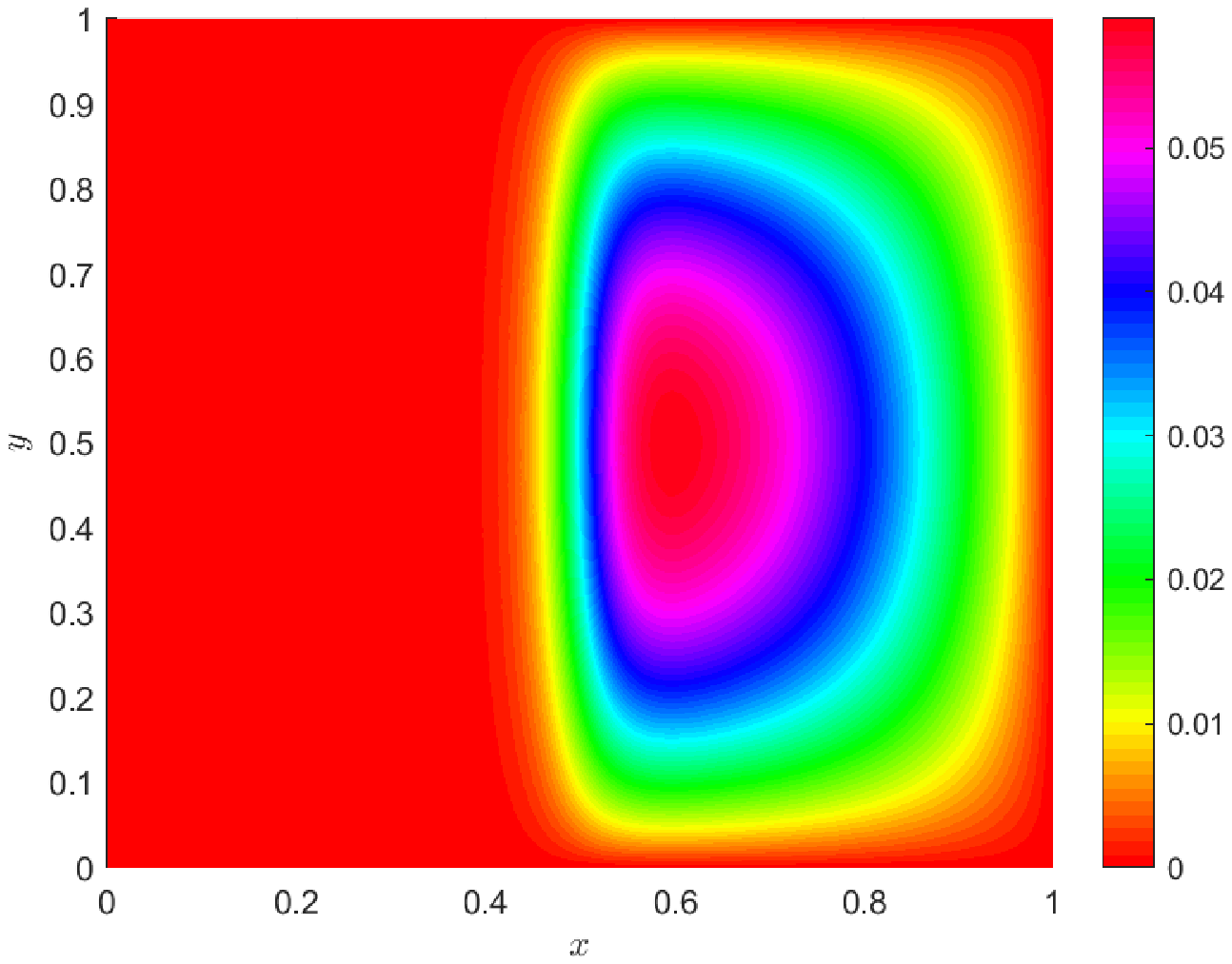}
}
\quad
\subfigure[]{
\includegraphics[width=7.5cm]{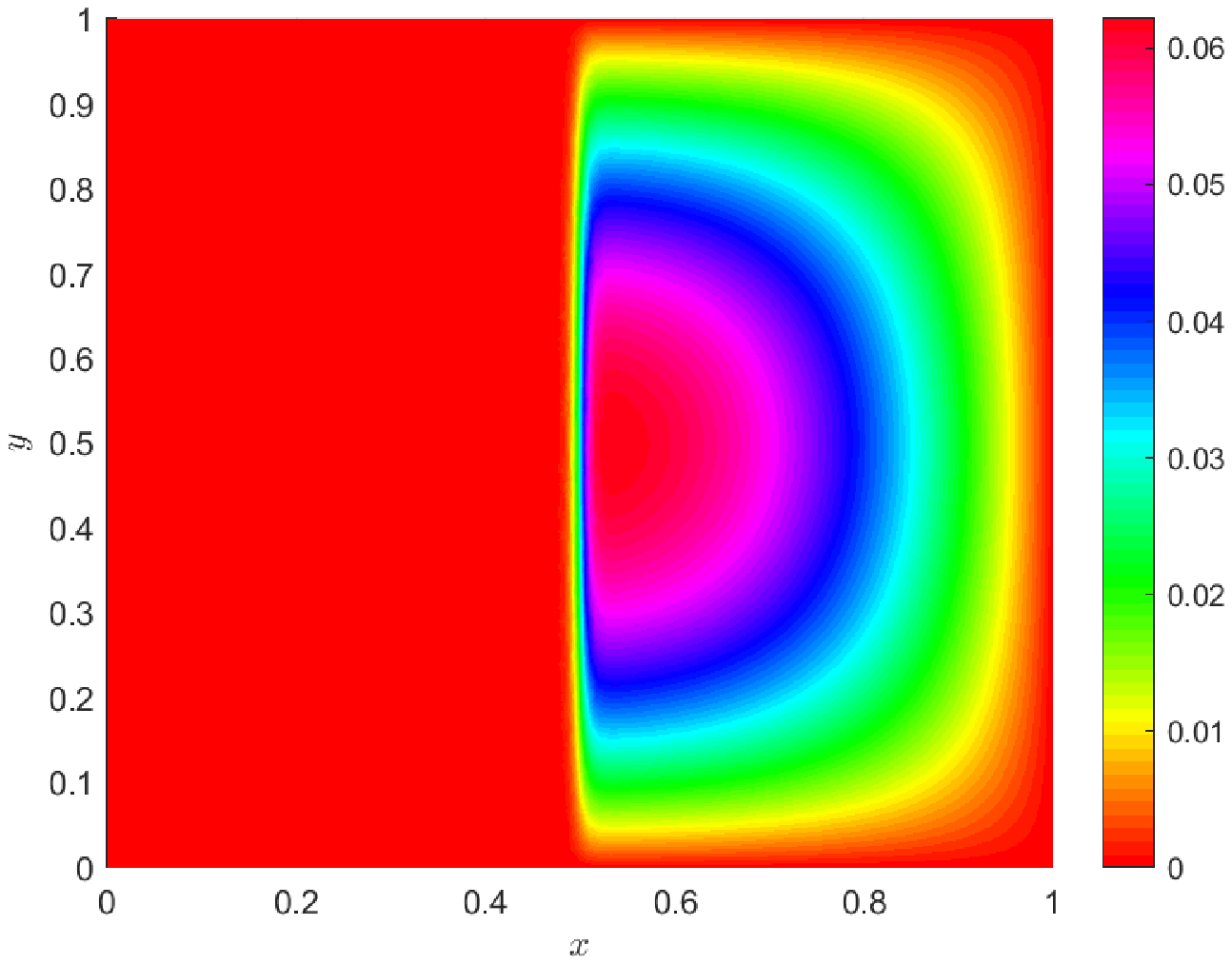}
}
\centering
\caption{\emph{Example 2.} Contour graphics of the exact solutions with $l_1=0.5$ and different value of $l_2$: (a) $l_2=0.05$; (b) $l_2=0.01$.}\label{fig:example2:solution:contours}
\end{figure}

Further we shall take the value of $l_2$ smaller and challenge behaviors of the reconstruction operator. In Figure \ref{fig:example2:singular:solution}, the exact solution of numerical results of 2nd-order reconstruction operator with $l_1=0.5$ and $l_2=0.01$ are presented over regular polygonal meshes. Clearly, the interior layer is also accurately captured in this case. Then the corresponding numerical errors measured in different norms with $\nu=10^{-3},10^{-9}$ over regular polygonal meshes are shown in Figure \ref{fig:example2:singularcr:l2dgsupg}. More specifically, two different thicknesses of the interior layer are observed in Figure \ref{fig:example2:solution:contours}, where the exact solution of \eqref{p1} is close to be singular at $x=0.5$ when $l_2=0.01$ but the other is relatively smooth when $l_2=0.05$. Therefore, the numerical errors in Figure \ref{fig:example2:singularcr:l2dgsupg} do not approach optimal orders due to the weak regularity.

\subsection{Example 3: Boundary layers.}
This example is taken from \cite{Ayuso20091391}. The data are $\bm{b}=[1,1]^\top$ and $c=0$, and we again vary $\nu$ to examine the performances of capturing boundary layers over regular polygonal meshes and general Voronoi meshes. The source term $f$ is chosen so that the exact solution of \eqref{p1}, with Dirichlet boundary conditions, is given by
\begin{equation}\label{exact:solution:bl}
u(x,y)=x+y(1-x)+\frac{\exp(-1/\nu)
-\exp(-(1-x)(1-y)/\nu)}{1-\exp(-1/\nu)}.
\end{equation}

\begin{figure}[htbp]
\centering
\subfigure[]{
\includegraphics[width=7.5cm]{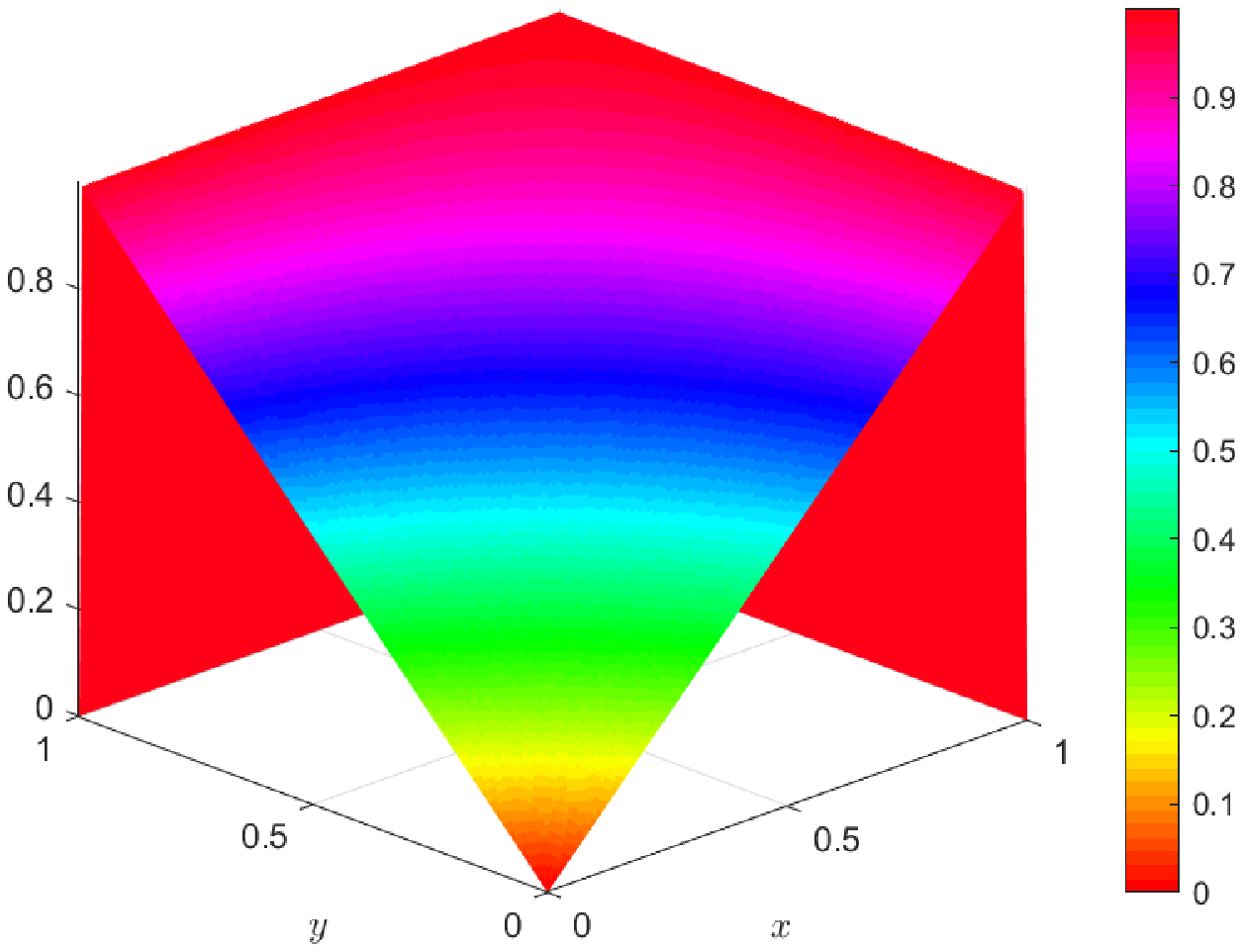}
}

\subfigure[]{
\includegraphics[width=7.5cm]{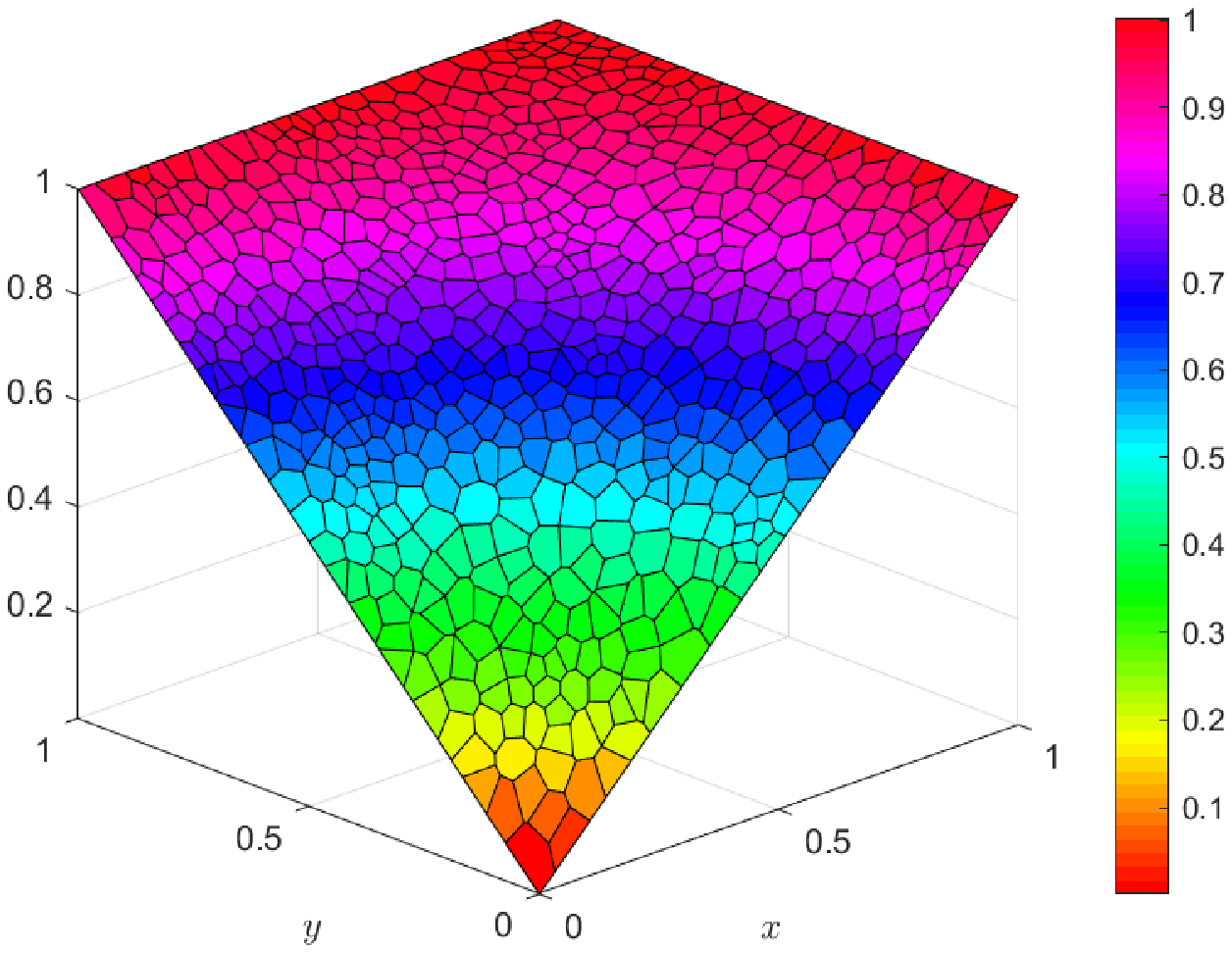}
}
\quad
\subfigure[]{
\includegraphics[width=7.5cm]{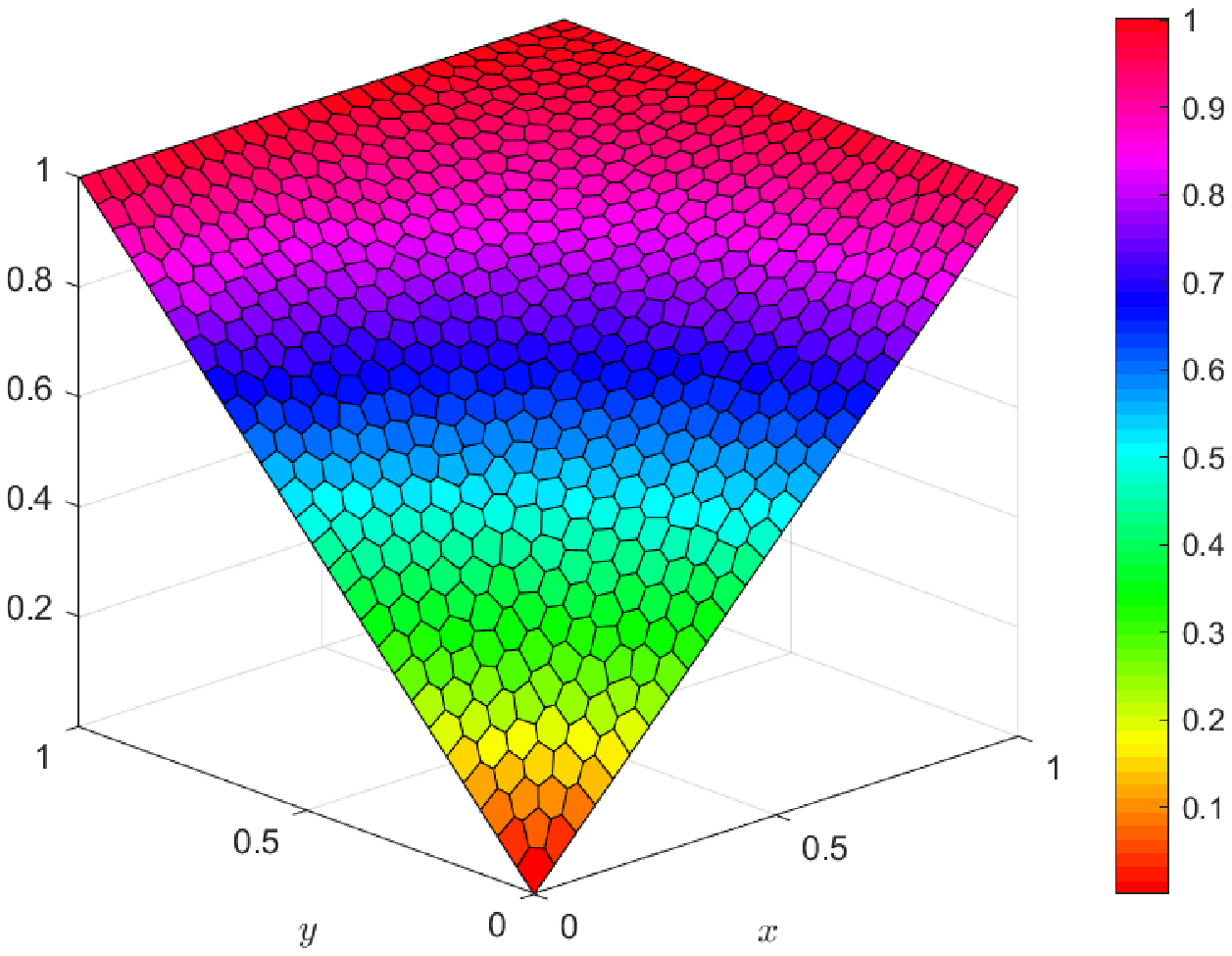}
}

\subfigure[]{
\includegraphics[width=7.5cm]{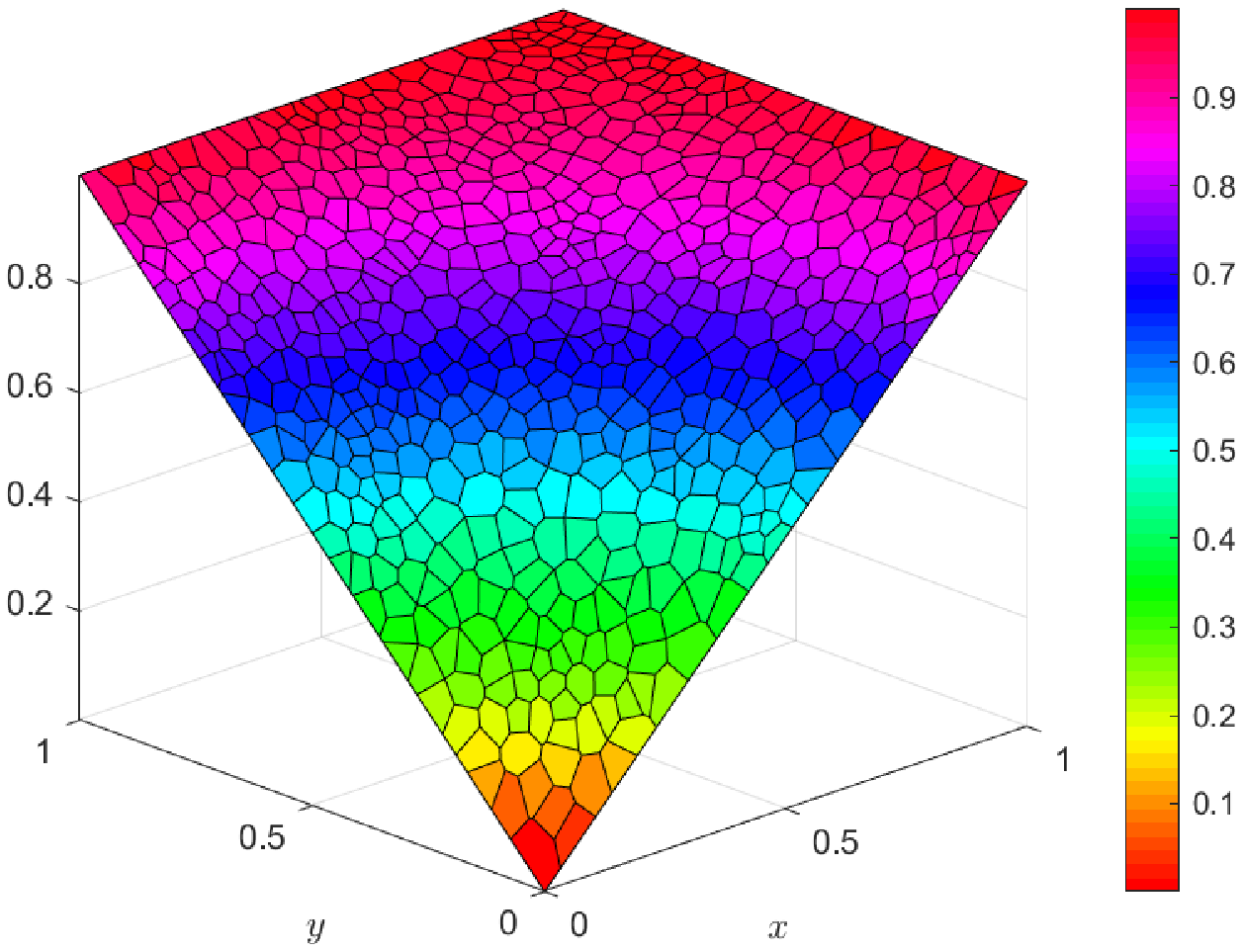}
}
\quad
\subfigure[]{
\includegraphics[width=7.5cm]{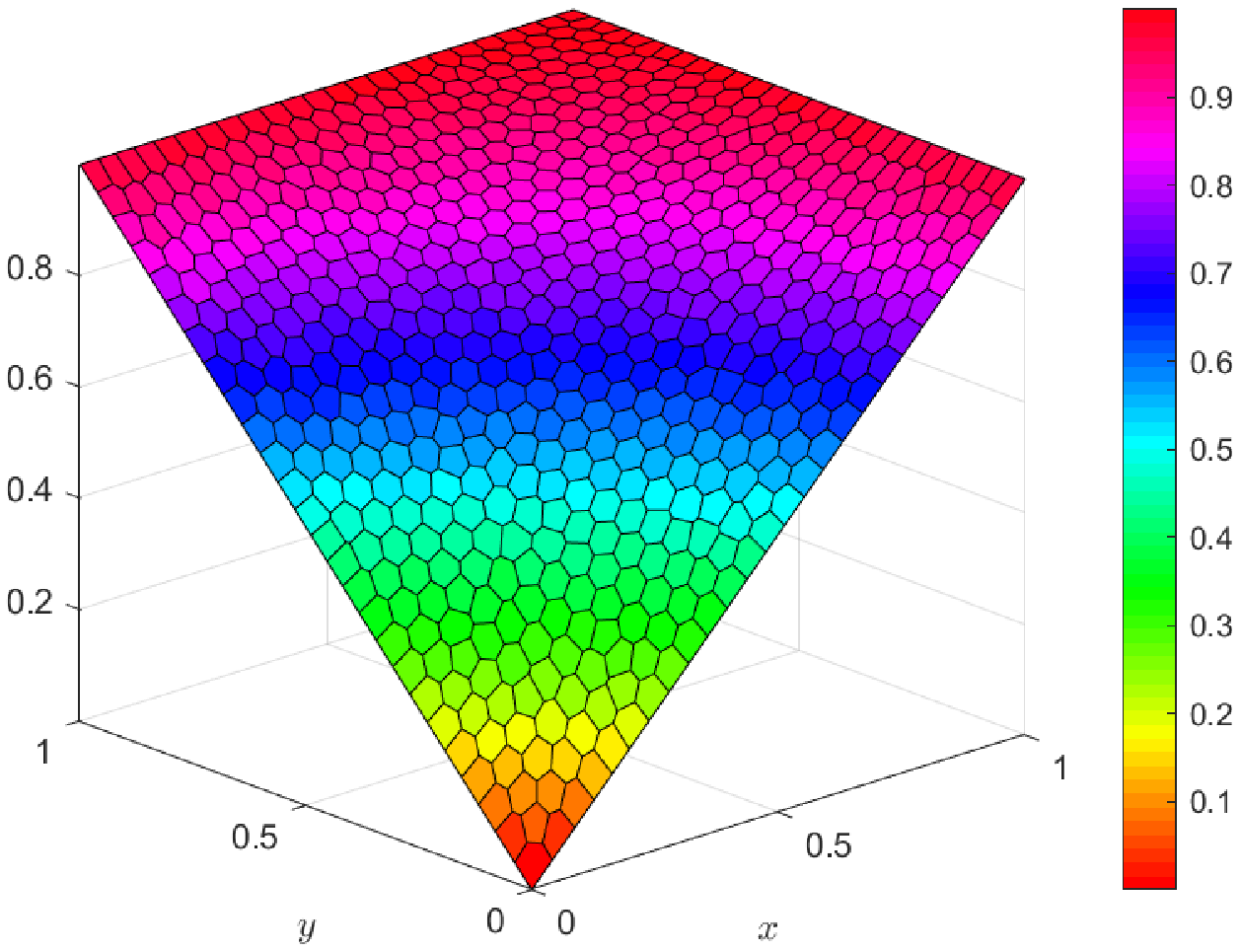}
}
\centering
\caption{\emph{Example 3.} A strongly convection-dominated case with $\nu=10^{-9}$ and 640 elements in the following meshes: (a) the exact solution; (b) the numerical solution over the general Voronoi mesh with $k=1$; (c) the numerical solution over the regular polygonal mesh with $k=1$; (d) the numerical solution over the general Voronoi mesh with $k=2$; (e) the numerical solution over the regular polygonal mesh with $k=2$.}\label{fig:example3:pV}
\end{figure}

In this example, the solution \eqref{exact:solution:bl} has boundary layers along $x=1$ and $y=1$ for $0<\nu\ll 1$. Figure \ref{fig:example3:pV} shows the exact solution and numerical solutions of different reconstruction order $k$ with $\nu=10^{-9}$. We have to admit that the reconstruction operator has a disappointing performance of capturing the so thin layer in a strongly convection-dominated regime although they have no spurious oscillations. The similar results are also shown in \cite{Ayuso20091391, Kim2019207, Lin20181482}.

\begin{figure}[htbp]
\centering
\subfigure[]{
\includegraphics[width=7.5cm]{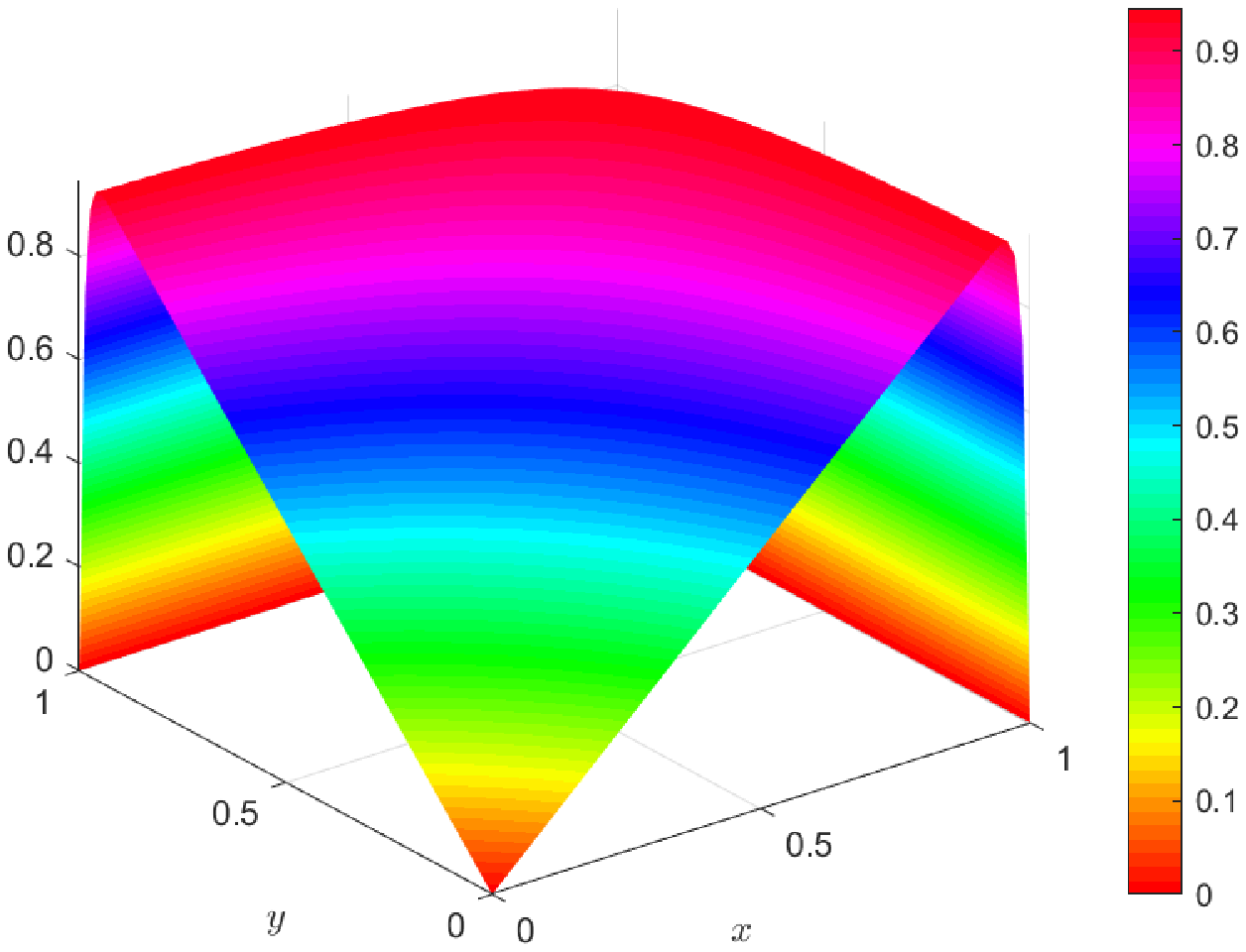}
}
\quad
\subfigure[]{
\includegraphics[width=7.5cm]{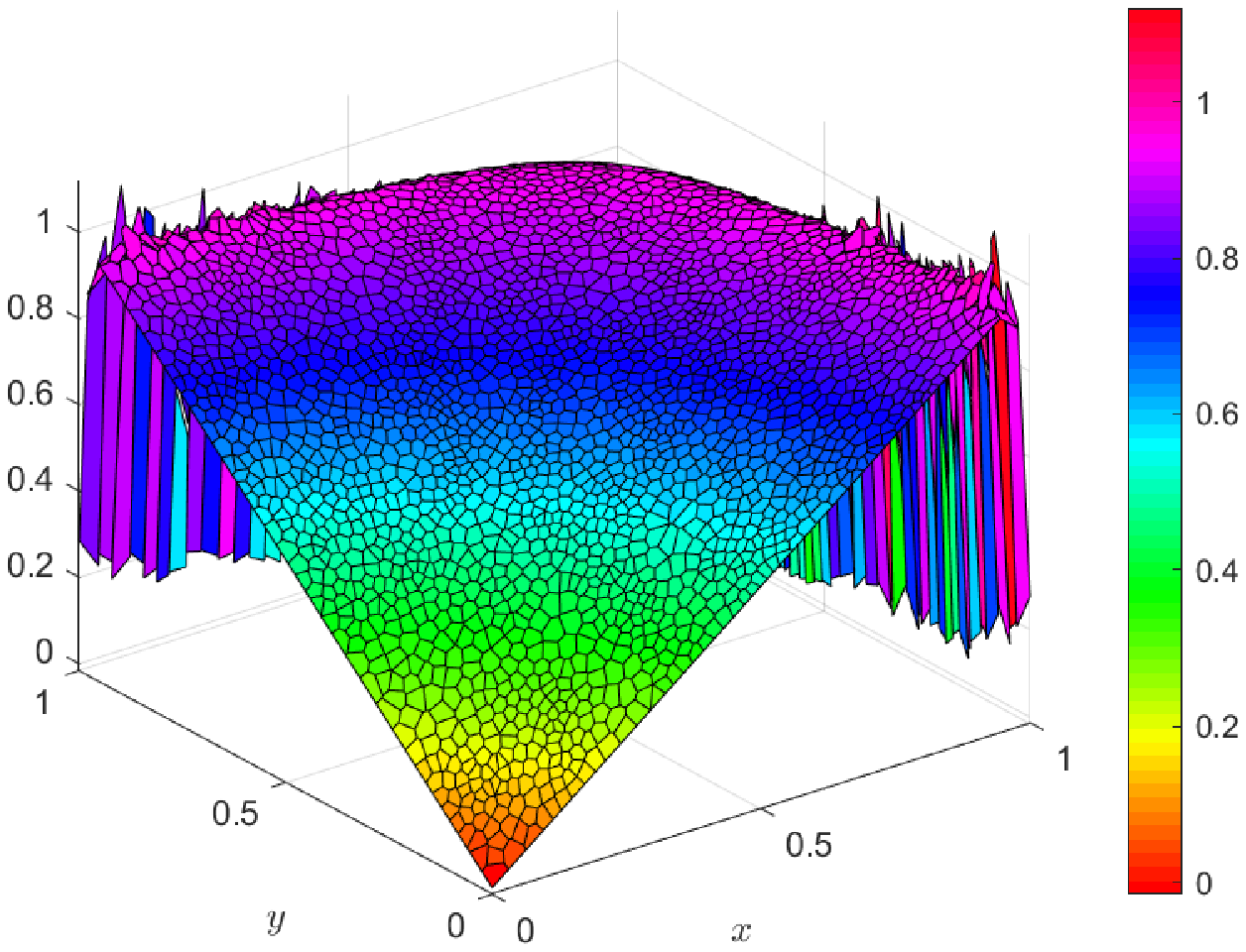}
}

\subfigure[]{
\includegraphics[width=7.5cm]{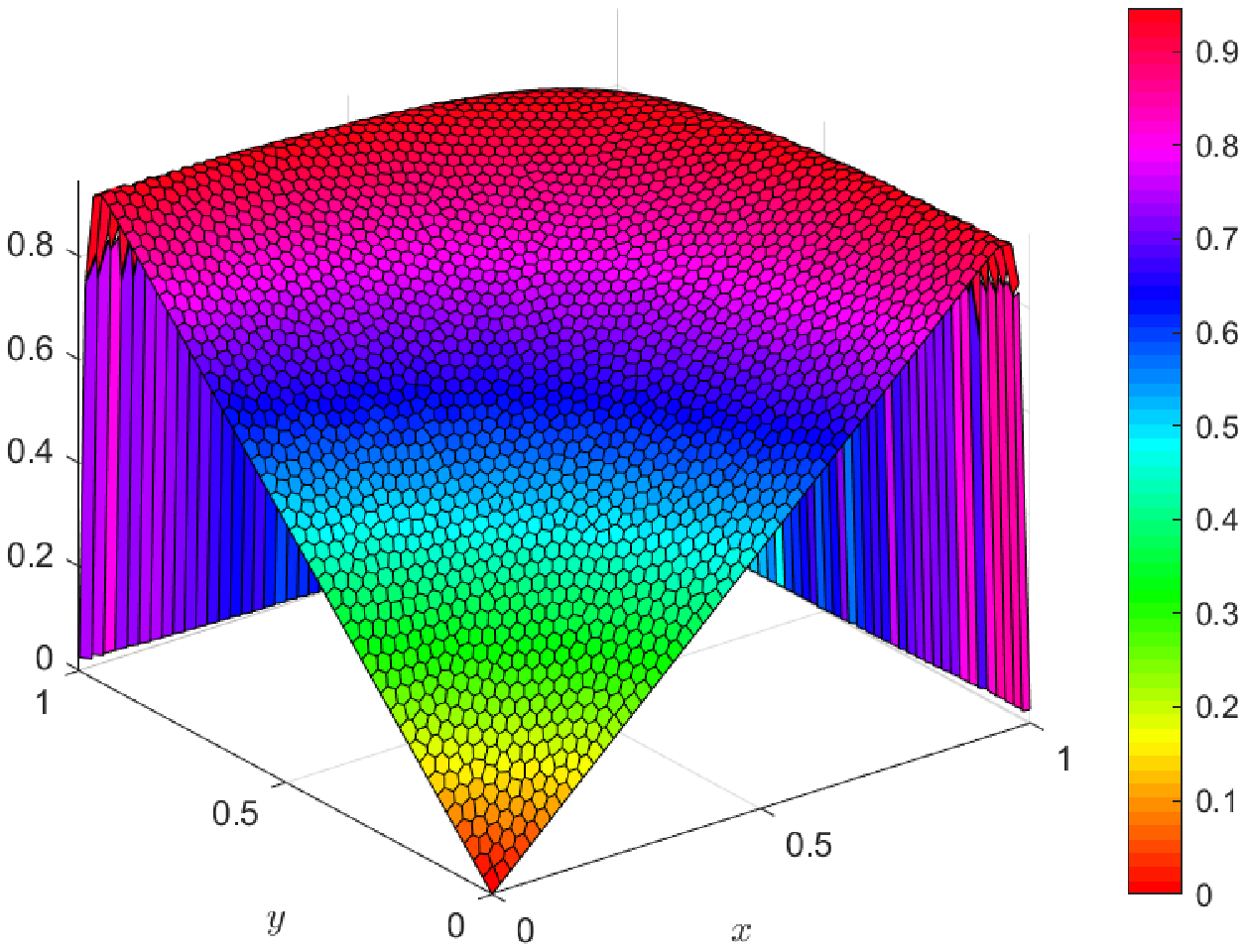}
}
\quad
\subfigure[]{
\includegraphics[width=7.5cm]{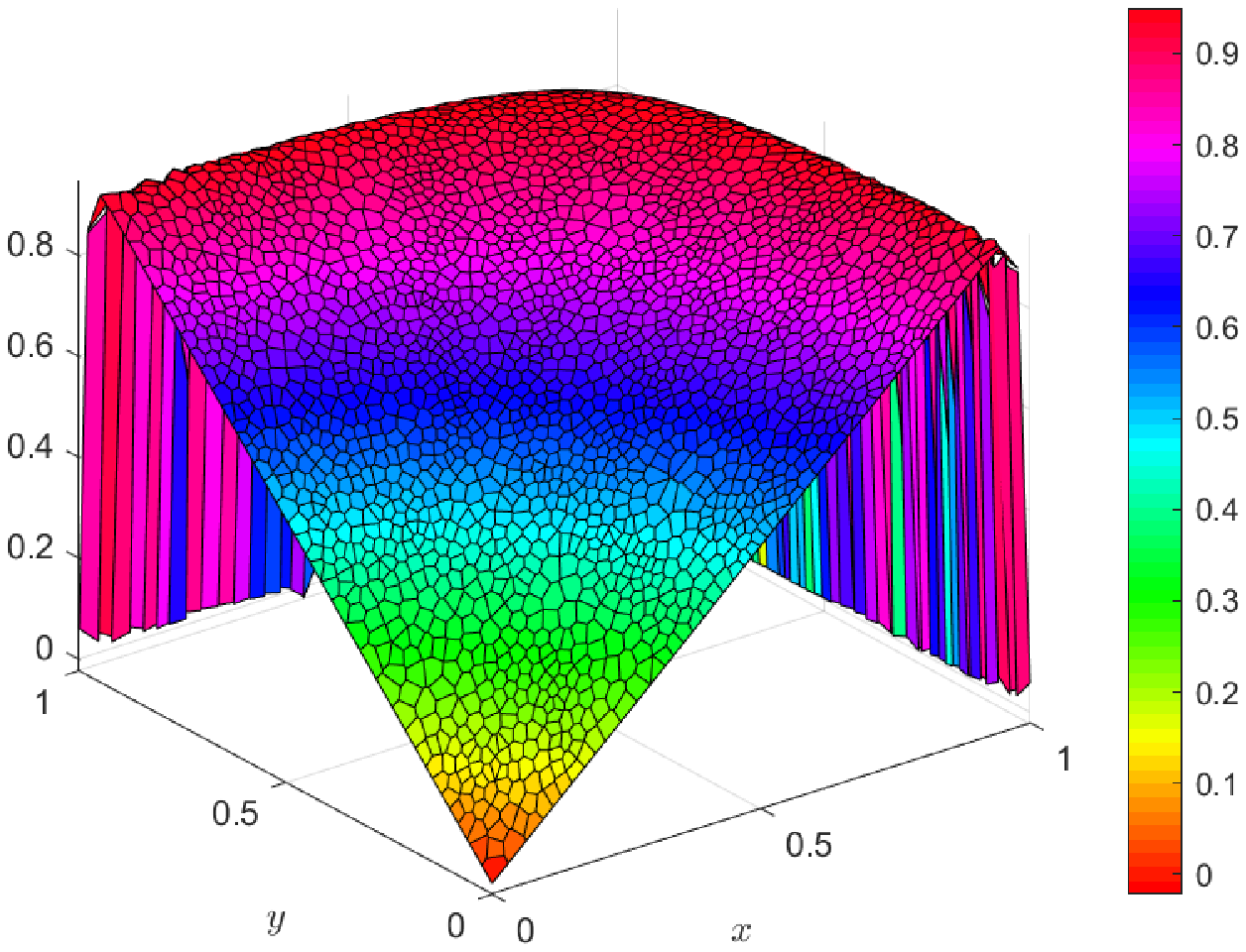}
}
\centering
\caption{\emph{Example 3.} A common case with $\nu=10^{-2}$ and 2560 elements in the following meshes: (a) the exact solution; (b) the numerical solution over the general Voronoi mesh with $k=1$; (c) the numerical solution over the regular polygonal mesh with $k=2$; (d) the numerical solution over the general Voronoi mesh with $k=2$.}\label{fig:example3:solution}
\end{figure}

\begin{figure}[htbp]
  \centering
  \includegraphics[width=0.7\textwidth]{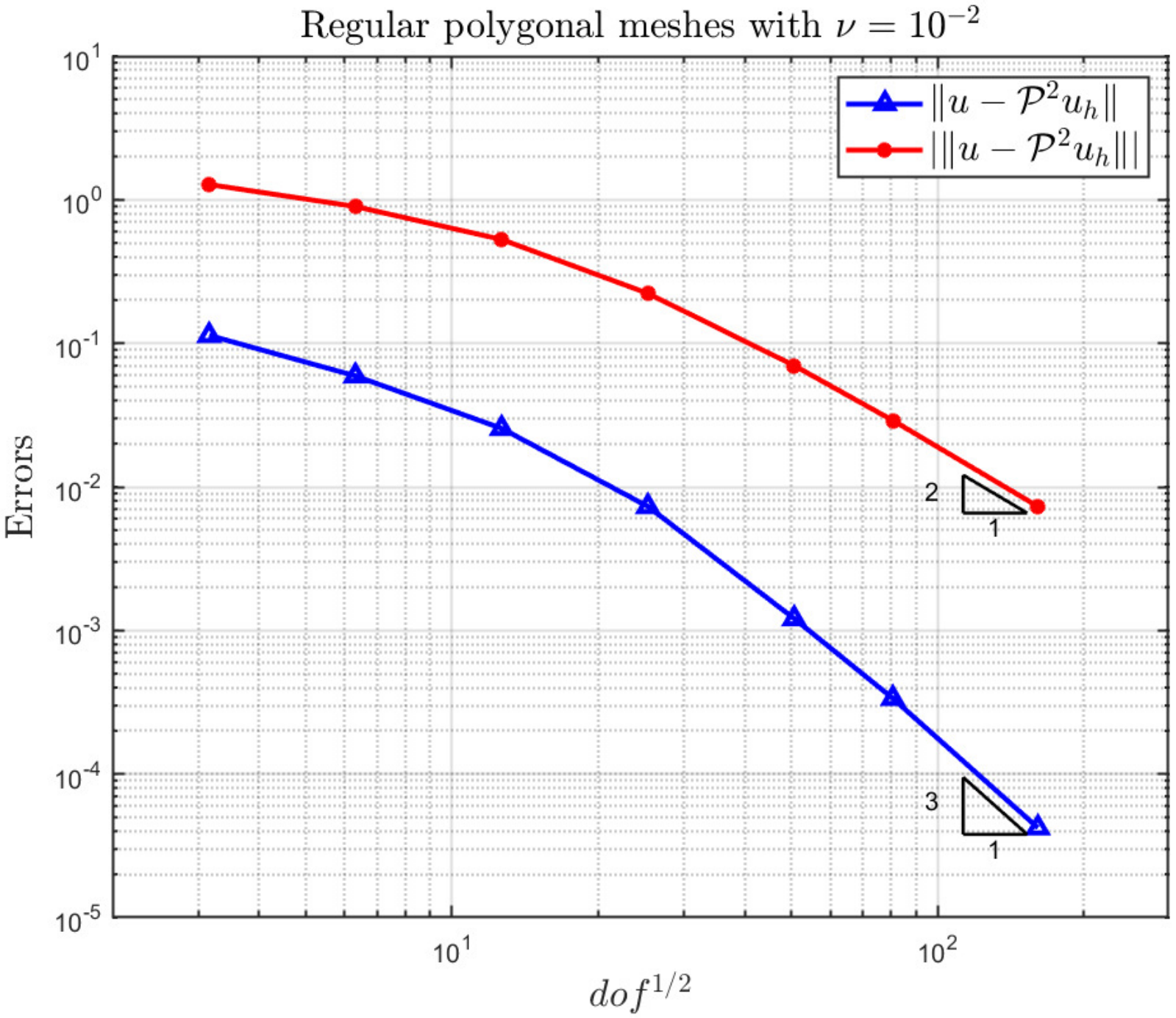}\\
  \centering
  \caption{\emph{Example 3.} Convergence rate results in the norm $\|\cdot\|$ and $\normdg{\cdot}$ with $k=2$ and $\nu=10^{-2}$ over regular polygonal meshes.}\label{fig:example3:cr:l2dg}
\end{figure}

Next, if we let $\nu=10^{-2}$, it cannot be called a ``boundary layer" problem, but we can also observe the behaviors over regular polygonal meshes and general Voronoi meshes. In Figure \ref{fig:example3:solution}, we present the exact solution and numerical solutions of different reconstruction order $k$ with $\nu=10^{-2}$. It is observed that the reconstruction operator has a good performance in both meshes. In Figure \ref{fig:example3:cr:l2dg}, for $\nu=10^{-2}$ the convergence results of 2nd-order reconstruction operator are provided versus refined polygon meshes in the norm $\|\cdot\|$ and $\normdg{\cdot}$ over global domain $\Omega$ including its boundary $\Gamma$. Clearly, the optimal error orders are obtained, which matches the results in Theorem \ref{priori:dg} and Theorem \ref{priori:l2}.

\begin{figure}[htbp]
\centering
\subfigure[]{
\includegraphics[width=7.5cm]{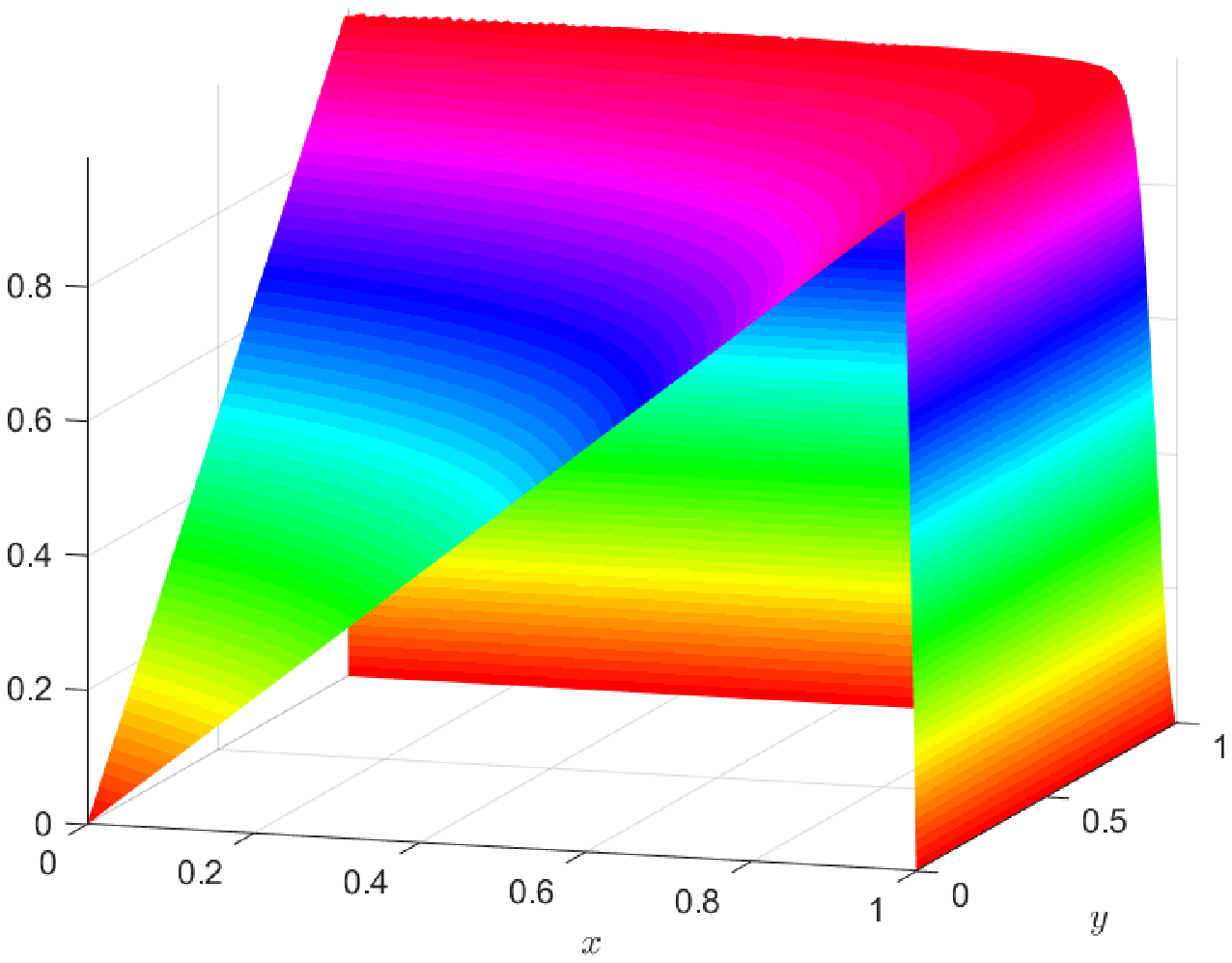}
}

\subfigure[]{
\includegraphics[width=7.5cm]{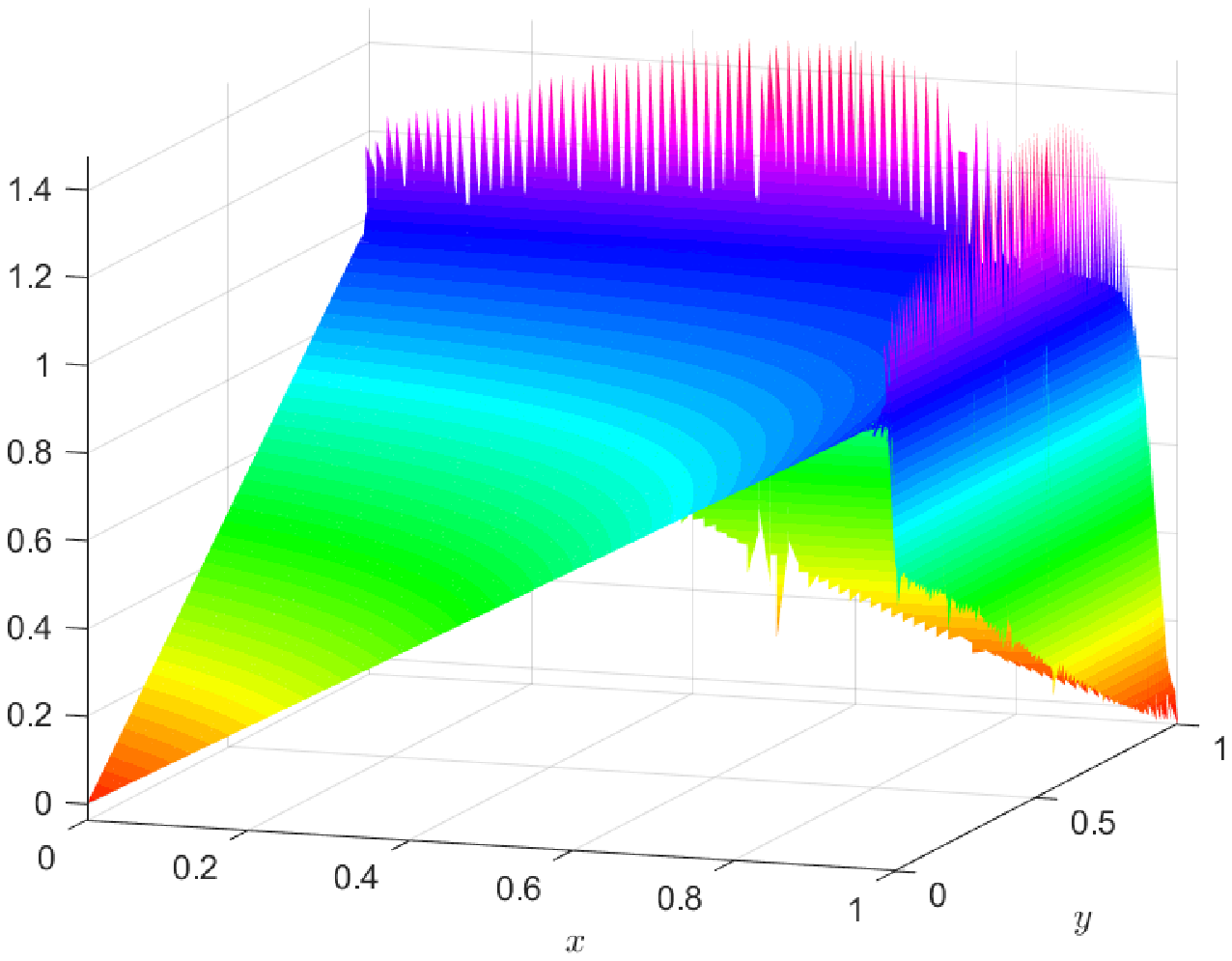}
}
\quad
\subfigure[]{
\includegraphics[width=7.5cm]{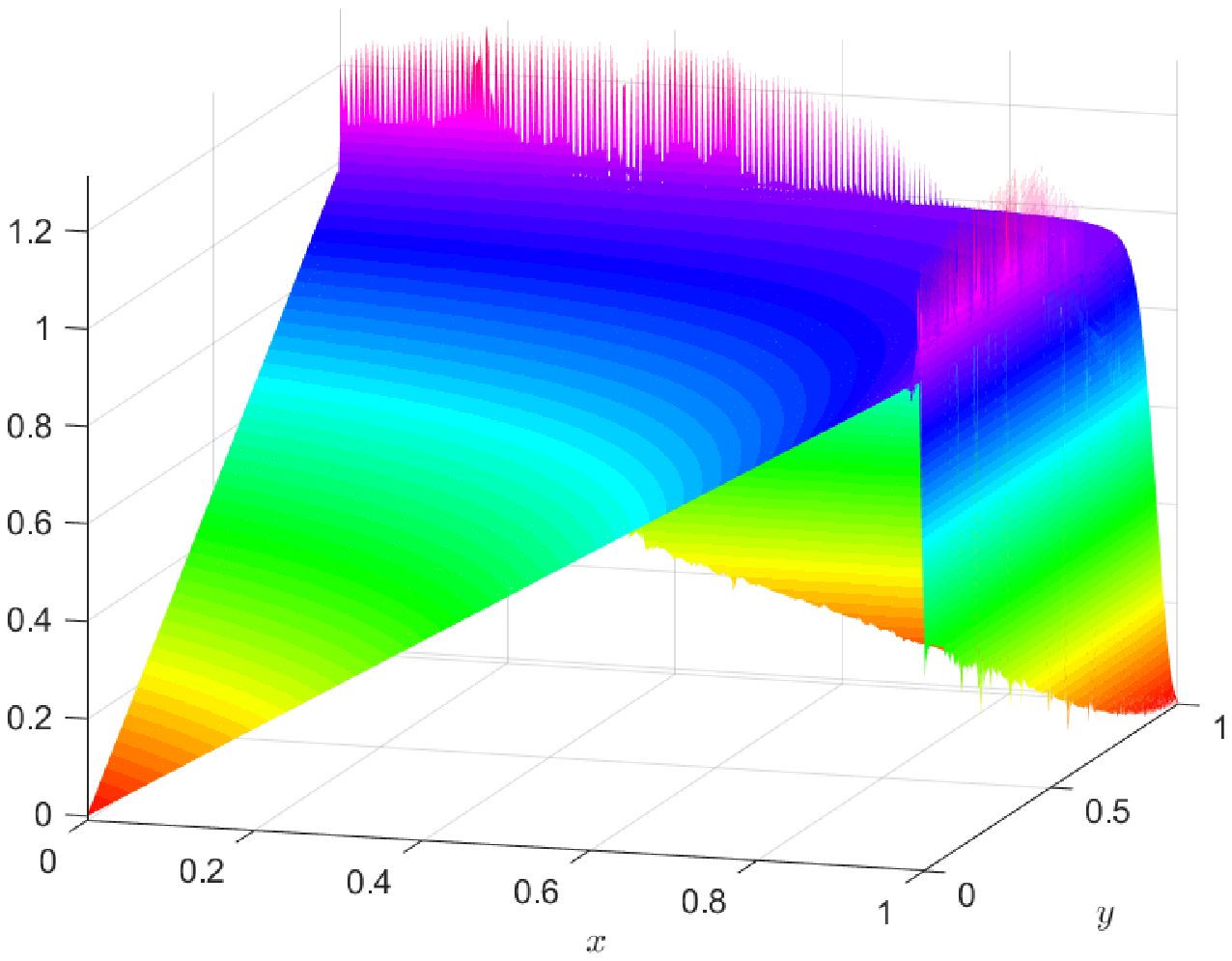}
}

\subfigure[]{
\includegraphics[width=7.5cm]{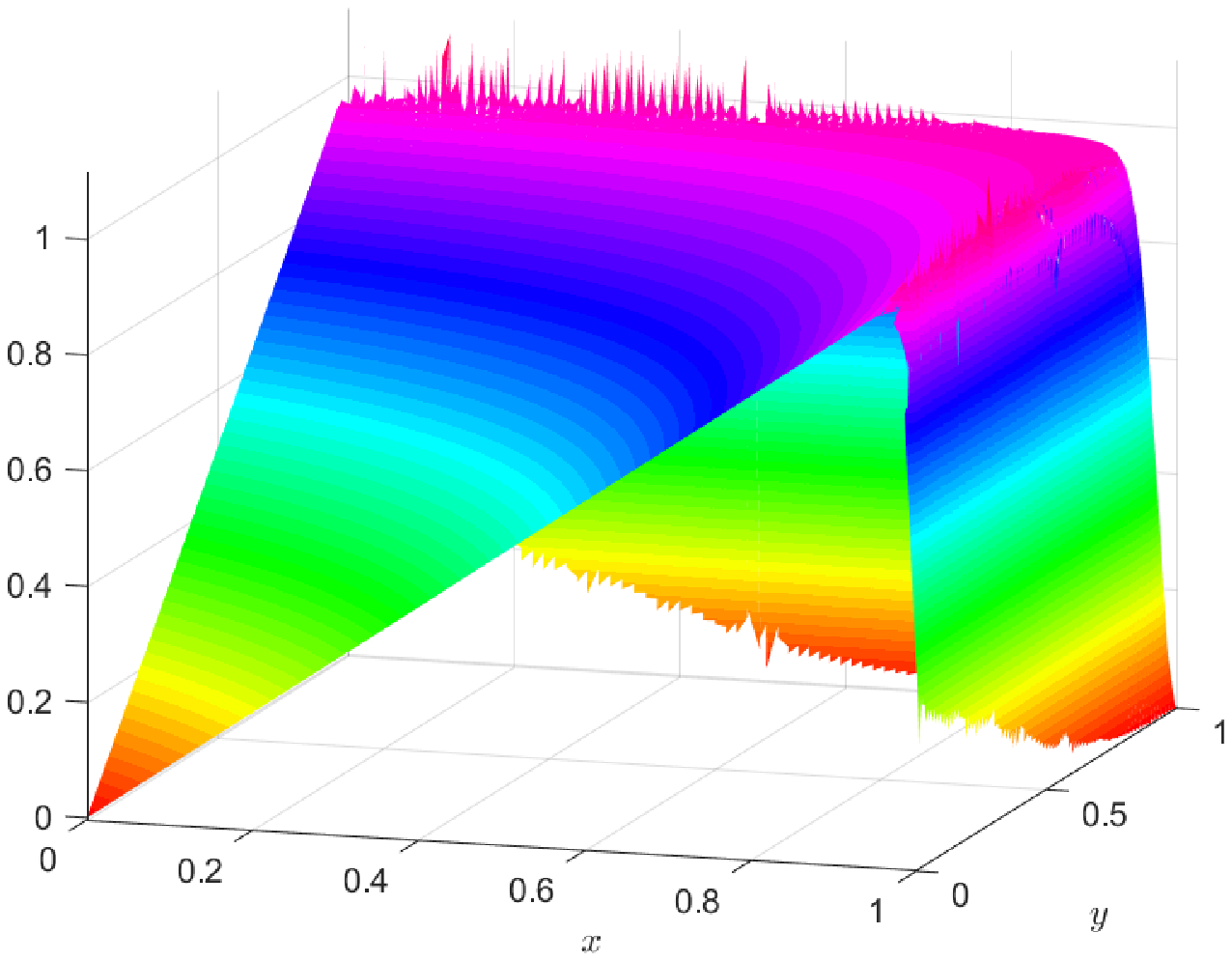}
}
\quad
\subfigure[]{
\includegraphics[width=7.5cm]{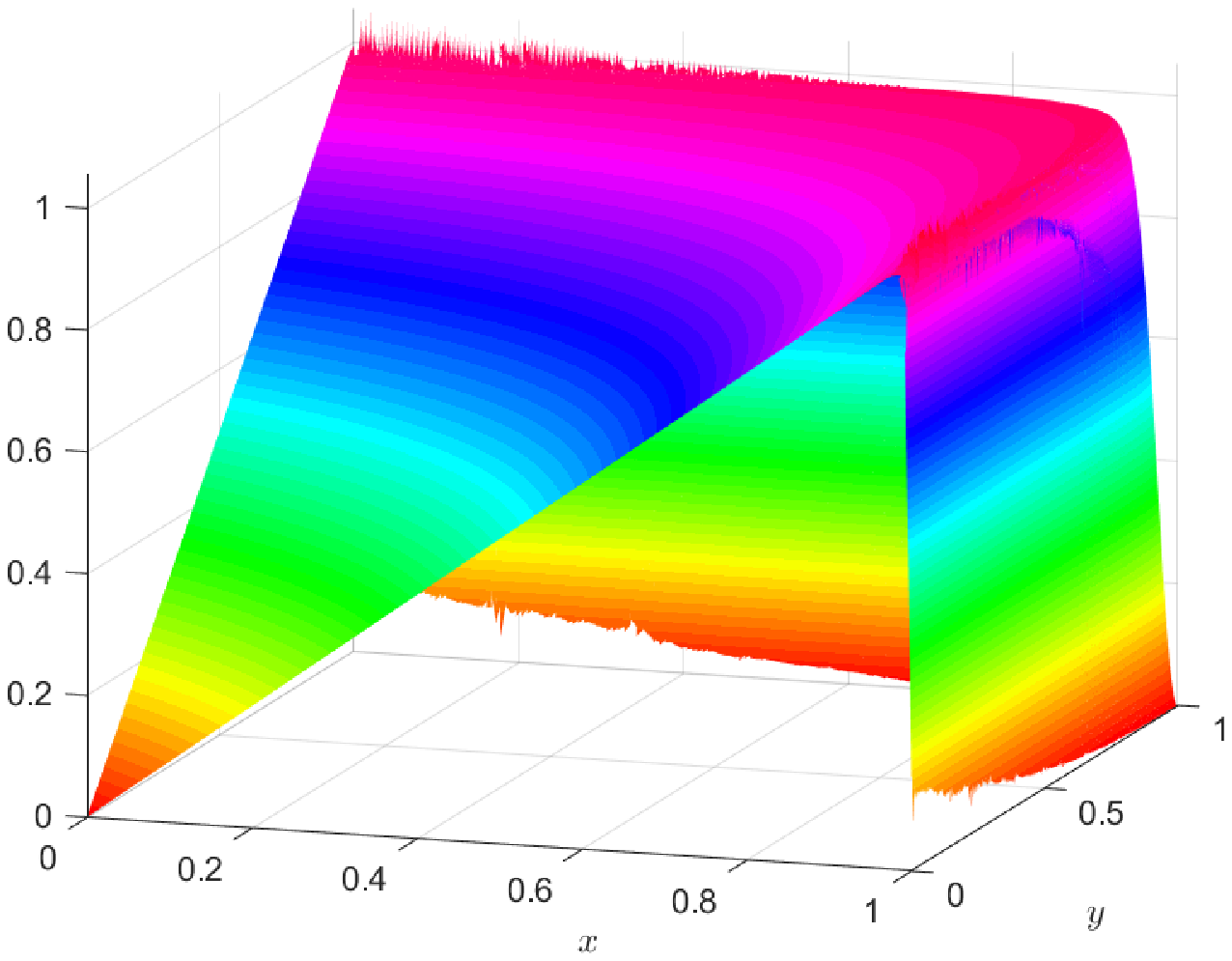}
}
\centering
\caption{\emph{Example 3.} An intermediate convection-dominated case over regular polygonal meshes with $\nu=10^{-3}$: (a) the exact solution; (b) the numerical solution with 6520 elements, $k=1$; (c) the numerical solution with 25800 elements, $k=1$; (d) the numerical solution with 6520 elements, $k=2$; (e) the numerical solution with 25800 elements, $k=2$.}\label{fig:example3:in:ea}
\end{figure}

Finally, in the intermediate regime, $\nu=10^{-3}$, we shall gradually capture the boundary layer as the higher value of $k$. In Figure \ref{fig:example3:in:ea}, for various reconstruction order $k$, we compare the numerical solutions over regular polygonal meshes with the exact solution. Further in Figure \ref{fig:example3:in:ea}, all the mesh grids keep invisible for better observations along $x=1$ and $y=1$). It can be observed that the ability of capturing the boundary layer is more powerful as the value of $k$ increases. In addition, we can also find that overshooting almost magically disappears when $k=2$.

\subsection{Example 4: Internal and boundary layers with discontinuous boundary conditions.}
Let $\bm{b}=[1/2,\sqrt{3}/2]^\top$, $c=0$ and $f=0$. The Dirichlet boundary conditions are given as follows:
\begin{equation}\label{ex4:bc:discon}
u(x,y)=g=
\begin{cases}
1&\text{on}~\{y=0,~0\leqslant x\leqslant 1\},\\
1&\text{on}~\{x=0,~y\leqslant 1/5\},\\
0&\text{elsewhere}.
\end{cases}
\end{equation}

This example is taken from \cite{Ayuso20091391, Kim2019207}. As interpreted in Example 6.5 of \cite{Kim2019207}, the exact solution is unknown, but it has both internal layer near the line $y=\sqrt{3}x+1/5$ and boundary layer along $x=1$ and $y=1$ for $0<\nu\ll 1$ due to the discontinuous Dirichlet boundary conditions given by \eqref{ex4:bc:discon}. Moreover, as the diffusion coefficient $\nu$ decreases, thickness of the internal layer becomes thinner.

\begin{figure}[htbp]
\centering
\subfigure[]{
\includegraphics[width=7.5cm]{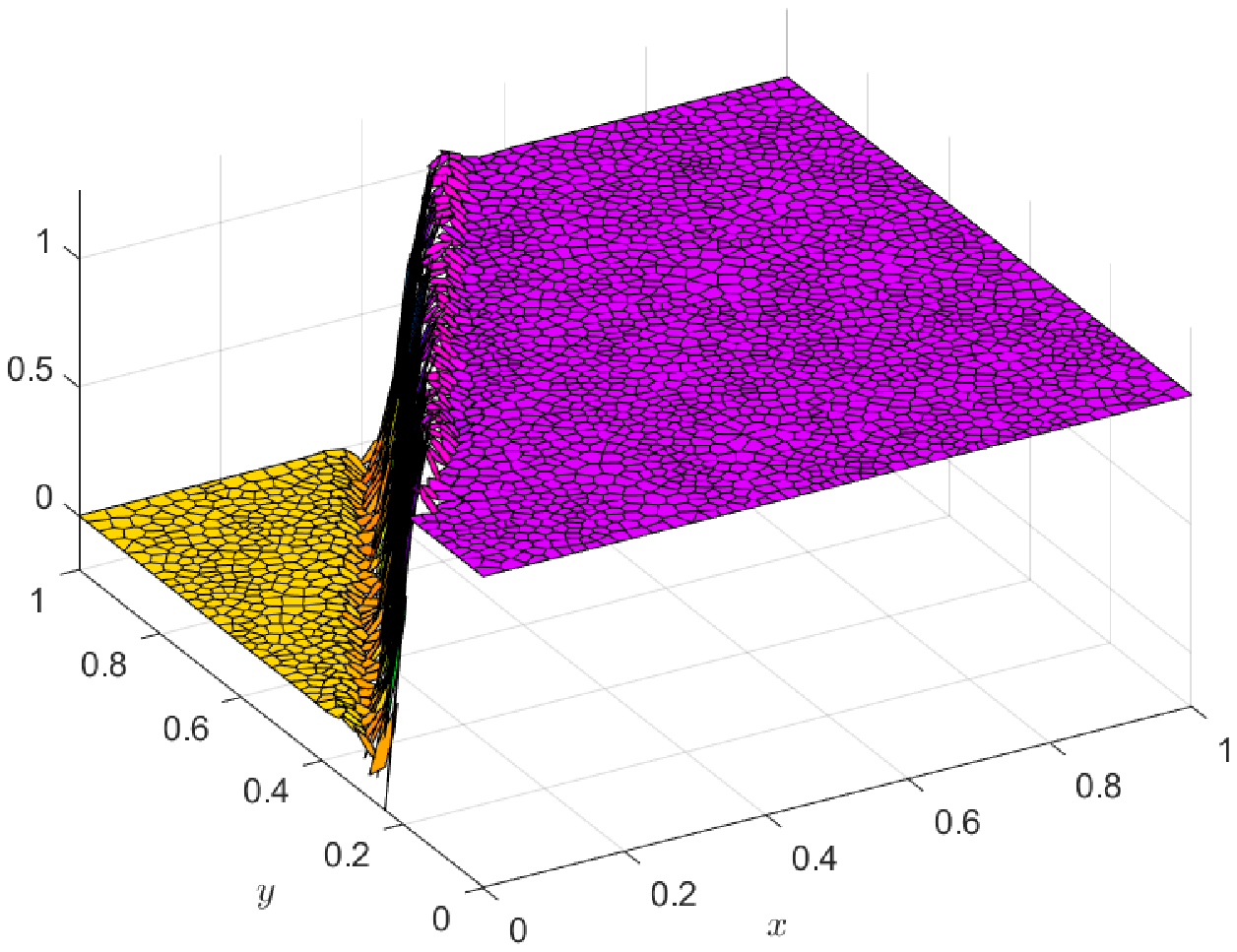}
}
\quad
\subfigure[]{
\includegraphics[width=7.5cm]{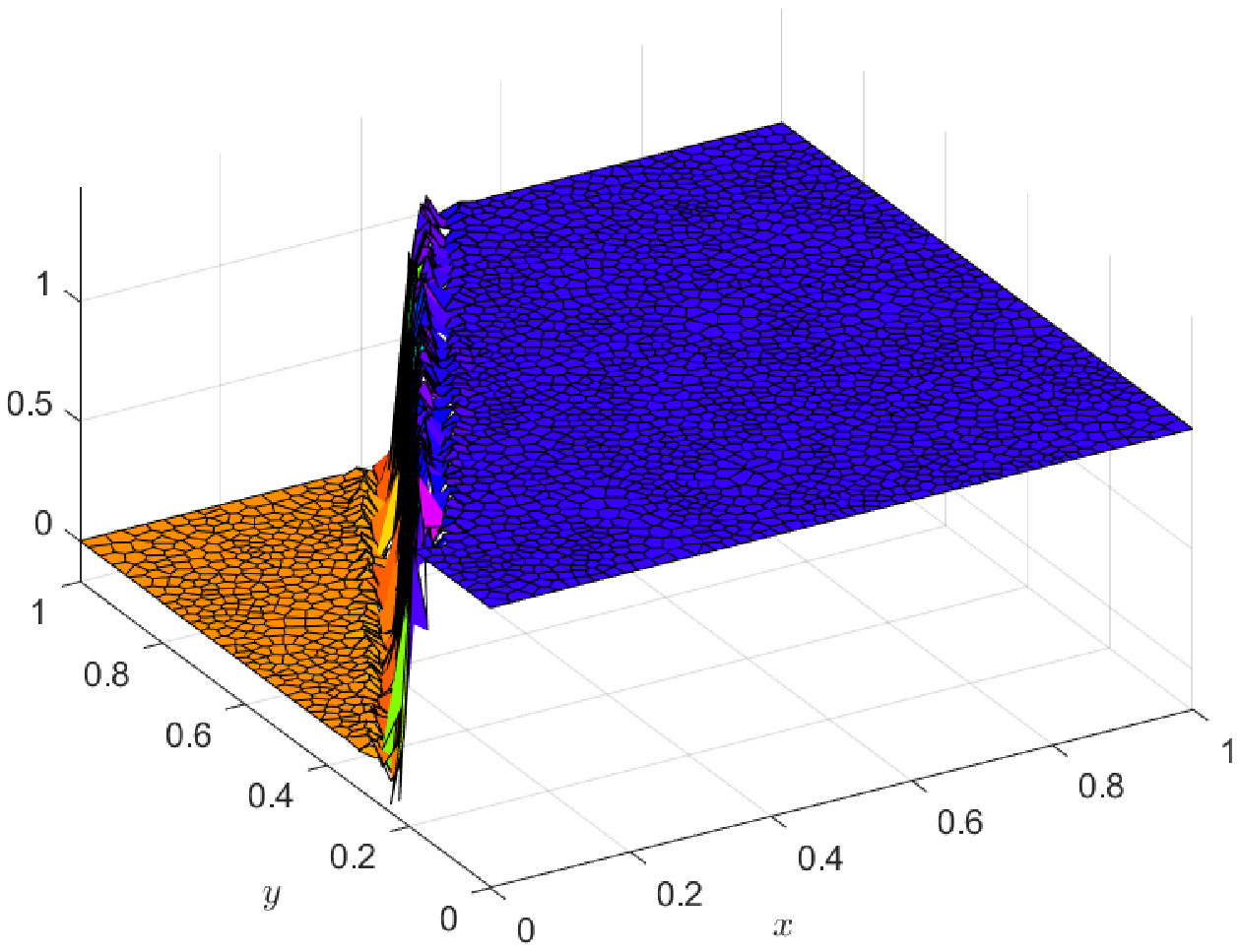}
}

\subfigure[]{
\includegraphics[width=7.5cm]{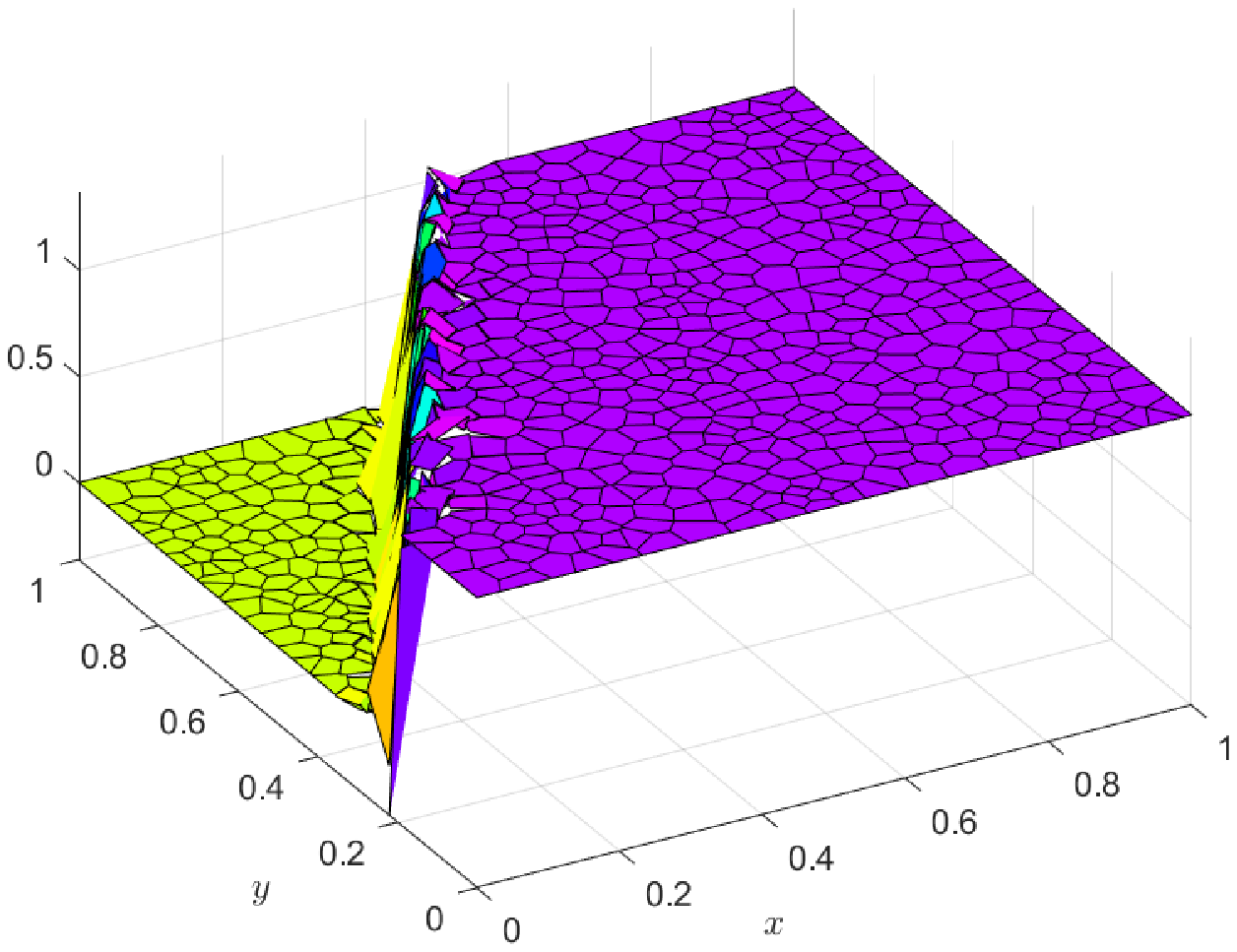}
}
\quad
\subfigure[]{
\includegraphics[width=7.5cm]{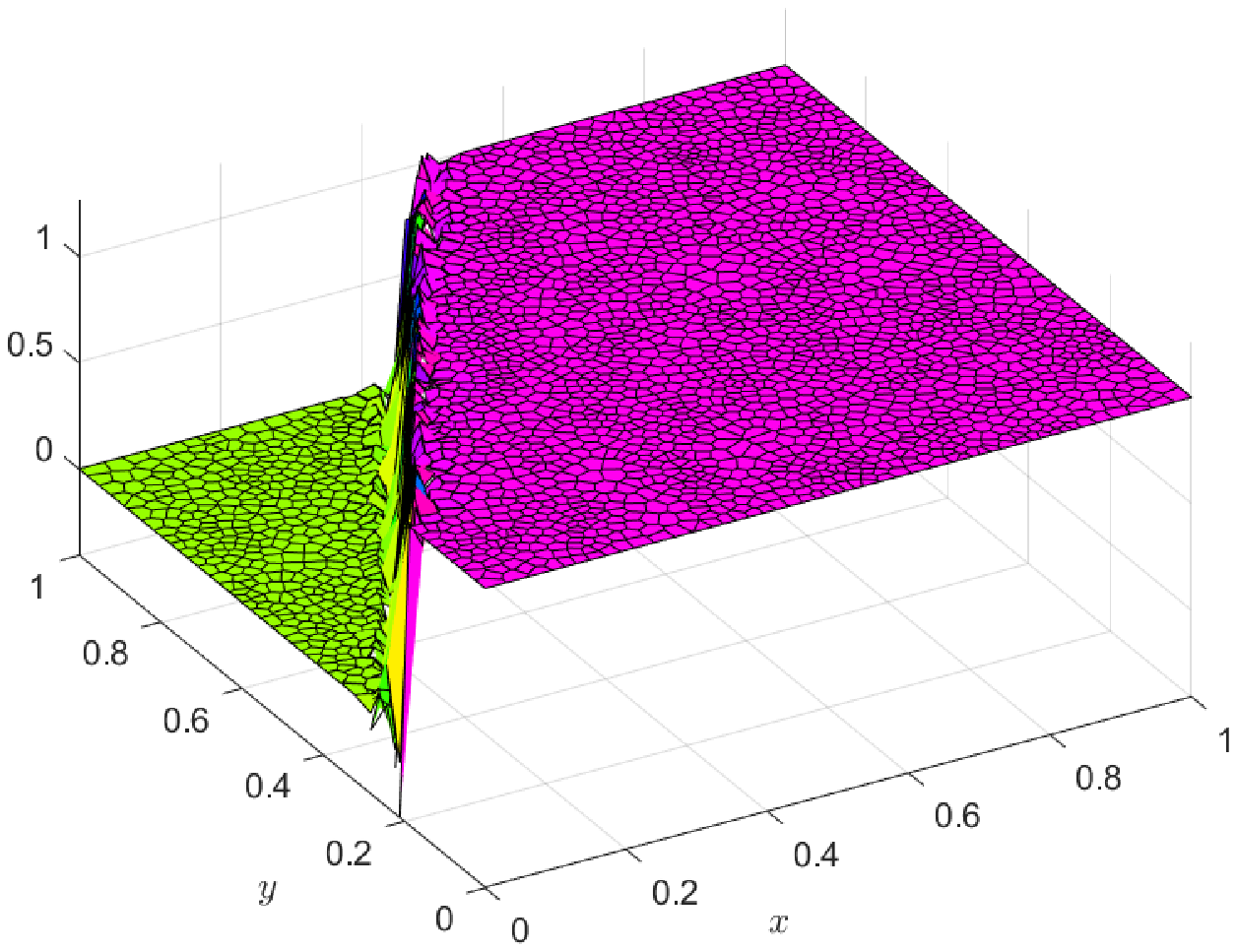}
}
\centering
\caption{\emph{Example 4.} Numerical solutions  with $\nu=10^{-9}$ and different reconstruction order $k$ over the general Voronoi meshes: (a) 2560 elements, $k=1$; (b) 2560 elements, $k=2$; (c) 640 elements, $k=3$; (d) 2560 elements, $k=3$;}\label{fig:example4:1e-9:ns}
\end{figure}

Figure \ref{fig:example4:1e-9:ns} shows the numerical solutions with different reconstruction order $k$ over the general Voronoi meshes for $\nu=10^{-9}$. Unfortunately, the boundary layer is not captured for any $k$ ($k=1,2,3$) although all the approximate solutions have no spurious oscillations along both lines $x=1$ and $y=1$. However, it is also observed that we have succeeded capturing the internal layer over such complicated meshes with small overshooting/undershooting near it, and thickness of the approximated internal layer is thinner as the value of $k$ increases, which reaches our expectations.

\begin{figure}[htbp]
\centering
\subfigure[]{
\includegraphics[width=7.5cm]{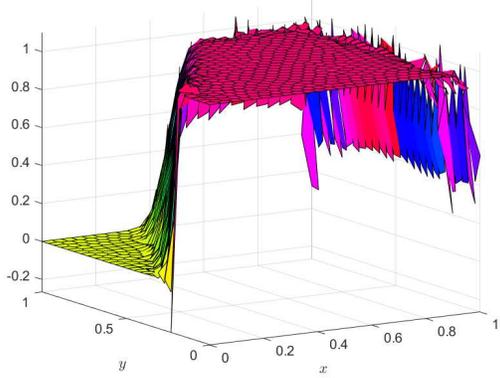}
}
\quad
\subfigure[]{
\includegraphics[width=7.5cm]{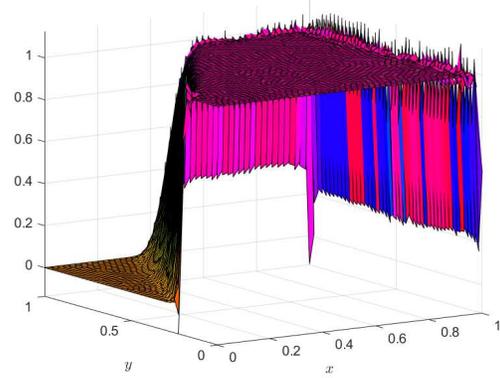}
}

\subfigure[]{
\includegraphics[width=7.5cm]{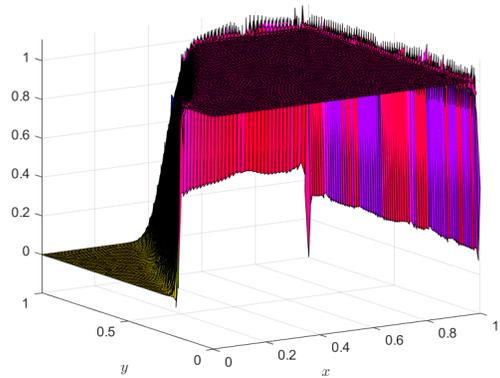}
}
\quad
\subfigure[]{
\includegraphics[width=7.5cm]{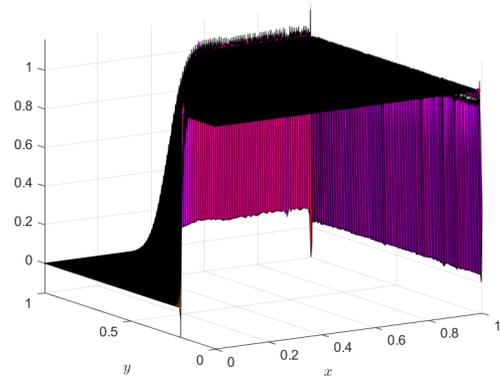}
}
\centering
\caption{\emph{Example 4.} Numerical solutions with $k=2$ and $\nu=10^{-3}$ over regular polygonal meshes: (a) 640 elements; (b) 2560 elements; (c) 6520 elements; (d) 25800 elements.}\label{fig:example4:1e-3:2ns}
\end{figure}

In the intermediate regime, $\nu=10^{-3}$, Figure \ref{fig:example4:1e-3:2ns} shows the performance of reconstruction order $k=2$ over regular polygonal meshes. It captures the internal layers with very small overshooting/undershooting. Further as degree of freedom becomes larger, the boundary layer can be well captured and overshooting/undershooting gradually disappears. Some similar results are also shown in \cite{Kim2019207}.

\section{Conclusions}
\label{conclusion}

In this paper, we provide the comprehensive process of constructing a patch reconstruction finite element space, and solve a weighted discrete least-square problem to derive a novel discontinuous finite element method over polytopic meshes. Together with the standard
symmetric interior penalty Galerkin method introduced in Section 2.6 of \cite{Riviere2008}, we propose a new DG approximation for the steady-state convection-diffusion-reaction problems. The most prominent characteristic is that arbitrary order accuracy can be achieved with only one degree of freedom per element. Not only reduce the scale of finial linear systems, it provides enough flexibility of designing polytopic meshes.

A variety of numerical experiments are introduced to demonstrate the theoretical results and to show performance of the reconstruction operator. In particular, it approximates the singularities well over refined polygonal meshes for the problem in a convection-dominated regime. It would be interesting to consider an improved operator to approximate boundary and internal layers well over polytopic meshes without showing spurious oscillations. We shall leave this issue to further exploration.

\section*{Acknowledgements}
This research was supported in part by National Science Foundation Grant DMS-11771348, Natural Science Foundation of Shaan Xi Province in 2019 (No.2019JQ--755) and Natural Science Foundation of Shaanxi Provincial Department of Education in 2019 (19JK0462).

\appendix{\section{Interpretations of (A4)}\label{app:proof:A4}}

In fact, the assumption (A4) is proposed literally from the problem \eqref{discrete:ls}. More specifically, (A4) requires that the number $\#\mathcal I_K$ (note that $\#\mathcal I_K=\# S(K)=M$) cannot be small, and it is up to at least $n_k=\dim(P_k)$ to ensure the unisolvence of \eqref{discrete:ls}. However, compared with the above qualitative statement, we prefer employing matrix analysis to discuss the unisolvence of the equivalent problem \eqref{discrete:wls:2}.

The existence of $(\underline{\bm{X}}_k^\top\underline{\bm{X}}_k)^{-1}$ is equivalent to that all singular values of $\underline{\bm{X}}_k$ are positive. We gather all singular values of $\underline{\bm{X}}_k$ into a sequence $\{\sigma_i\}_{i=1}^l$, $l=\min\{M, n_k\}$, and the sequence satisfies that
\[\sigma_1\geqslant\cdots\geqslant
\sigma_r>0=\sigma_{r+1}=\cdots=\sigma_l\]
for some $r\leqslant l$. So roughly speaking,
the demand $M\geqslant n_k$ ensures the possibility that each singular value of $\underline{\bm{X}}_k$ is positive. Then for the reason that bigger value of $M$ is better, we propose the following lemma to be an appropriate interpretation:
\begin{lemma}\label{svd}
Suppose that $m\geqslant n+1$, $n\geqslant 2$, and $\sigma_1,\cdots,\sigma_n$ are all singular values of a matrix $\bm{A}\in\mathbb{R}^{m\times n}$ with
$\sigma_1\geqslant \cdots\geqslant\sigma_n\geqslant 0$. Then we randomly delete a row of $\bm{A}$ to obtain another matrix $\widetilde{\bm{A}}\in\mathbb{R}^{(m-1)\times n}$. In addition, we denote all singular values of $\widetilde{\bm{A}}$ by $\widetilde{\sigma}_1,\cdots,\widetilde{\sigma}_n$ with $\widetilde{\sigma}_1\geqslant \cdots\geqslant\widetilde{\sigma}_n\geqslant 0$. Then the following results hold:
\begin{equation}\label{sv}
\sigma_i\geqslant\widetilde{\sigma}_i
\geqslant\sigma_{i+1},~ i=1,\cdots,n-1,\quad\sigma_n\geqslant
\widetilde{\sigma}_n.
\end{equation}
\end{lemma}
\begin{pf}
We denote by $\mathcal U_r$ the set of all $r\times r$ unitary matrices.
Let
\begin{equation*}
\bm{H}=
\left[
  \begin{array}{cc}
    \bm{0} & \bm{A} \\
    \bm{A}^\top & \bm{0} \\
  \end{array}
\right],
\end{equation*}
then $\bm{H}\in\mathbb{R}^{(m+n)\times(m+n)}$ is a Hermite matrix. We let $\bm{A}=\bm{U}_1\bm{\Sigma}\bm{V}^\top$ be the singular value decomposition (SVD) of $\bm{A}$, where $\bm{U}=[\bm{U}_1,\bm{U}_2]\in
\mathcal U_m$, $\bm{V}\in\mathcal U_n$, $\bm{U}_1$ is the first $n$ columns of $\bm{U}$, and $\bm{\Sigma}=\text{diag}\{\sigma_1,\cdots,\sigma_n\}$.
Let
\begin{equation*}
\bm{Q}=\frac{1}{\sqrt{2}}
\left[
  \begin{array}{ccc}
    \bm{U}_1 & \sqrt{2}\bm{U}_2 & \bm{U}_1 \\
    \bm{V} & \bm{0} & -\bm{V} \\
  \end{array}
\right],
\end{equation*}
and we note that $\bm{Q}\in\mathcal U_{m+n}$. Then by some derivations, we obtain that
\begin{equation}\label{pf:svd}
\bm{Q}^\top\bm{H}\bm{Q}
=\left[
  \begin{array}{ccc}
    \bm{\Sigma} & \bm{0} & \bm{0} \\
    \bm{0} & \bm{0} & \bm{0} \\
    \bm{0} & \bm{0} & -\bm{\Sigma} \\
  \end{array}
\right].
\end{equation}
It follows from \eqref{pf:svd} that the eigenvalue set of $\bm{H}$ is
\begin{equation}\label{eigen:H}
\lambda(\bm{H})=\big\{
\sigma_1,\cdots,\sigma_n,0,\cdots,0,-\sigma_n,
\cdots,-\sigma_1\big\}.
\end{equation}
In \eqref{eigen:H} it is easy to observe that the first $n$ numbers of $\lambda(\bm{H})$ are all singular values of $\bm{A}$.

On the other hand, we let
\begin{equation*}
\widetilde{\bm{H}}=
\left[
  \begin{array}{cc}
    \bm{0} & \widetilde{\bm{A}} \\
    \widetilde{\bm{A}}^\top & \bm{0} \\
  \end{array}
\right],
\end{equation*}
then $\widetilde{\bm{H}}\in\mathbb{R}^{(m+n-1)\times(m+n-1)}$ is a Hermite matrix. We note that $\widetilde{\bm{H}}$ can be regarded as a matrix given by deleting the $j$th row and the $j$th column from $\bm{H}$ for a random number $j$. Following the same techniques of deriving \eqref{eigen:H}, we obtain the eigenvalue set of $\widetilde{\bm{H}}$ is
\begin{equation}\label{eigen:tH}
\lambda(\widetilde{\bm{H}})=\big\{
\widetilde{\sigma}_1,\cdots,
\widetilde{\sigma}_n,0,\cdots,0,
-\widetilde{\sigma}_n,
\cdots,-\widetilde{\sigma}_1\big\}.
\end{equation}
Similarly, in \eqref{eigen:tH} the first $n$ numbers of $\lambda(\widetilde{\bm{H}})$ are all singular values of $\widetilde{\bm{A}}$.

Finally, it follows from \eqref{eigen:H}--\eqref{eigen:tH}, and the separation theorem originating from the famous ``min-max" theorem for eigenvalues of Hermite matrices that the result \eqref{sv} holds.

\end{pf}

The above lemma implies that the smallest singular value (always nonnegative) of $\underline{\bm{X}}_k$ becomes bigger as the value of $M$ increasing. So the demand is reasonable and necessary that the number $M$ should be greater than $n_k$ to approach the assumption (A4).

\section{Proofs of all lemmas in Section \ref{useful:pk}}\label{append:pf:24}
\subsection{Proof of Lemma \ref{estimates:op}}
\begin{pf}
The proof of Lemma \ref{estimates:op} can be found in Lemma 4 of \cite{Li2019268}. Compared with the results in Lemma 4 of \cite{Li2019268}, our results \eqref{estimate:L2op} and \eqref{estimate:H1op} have more concise expressions. More specifically, in \cite{Li2019268}, the primal result of \eqref{estimate:L2op} is
\begin{equation}\label{primal:estimate:L2op}
\|v-\mathcal P^kv\|_{0,K}
\leqslant Ch_Kd_K^k|v|_{k+1,S(K)},
\end{equation}
where $d_K$ is the diameter of $S(K)$, and $C$ is a constant depends only on $N$, $\gamma$ and $k$. Based on (A2), (A3), and the process of forming the element patch $S(K)$ for any $K\in\mathcal T_h$ in Section \ref{subsec:ep}, we can find that the geometric relationship $d_K\leqslant Ch$ holds with $C$ a positive constant depending only on $k$. By substituting the case $d_K\leqslant Ch$ into \eqref{primal:estimate:L2op}, we obtain the result \eqref{estimate:L2op}. Further it follows the same techniques that \eqref{estimate:H1op} also holds.
\end{pf}

\subsection{Proof of Lemma \ref{estimates:op2}}
\begin{pf}
The result \eqref{estimate:trace} follows directly from \eqref{re:polygontr}, Lemma \ref{estimates:op} and \eqref{primal:estimate:L2op}, so we omit the details.
\end{pf}

\section{Proofs of all lemmas and theorems in Section \ref{sec:DG:appro}}\label{append:pf:stable}
\subsection{Proof of Lemma \ref{useful:inequality}}
\begin{pf}
For any $e\in\mathcal E_h^o$, we denote by $K_e^1$, $K_e^2$ the two adjacent elements of $\mathcal T_h$ such that $K_e^1\cap K_e^2=e$. Using the triangle inequality and \eqref{inverse:trace2}, for any $v_h\in V_h^k$ we obtain that
\begin{equation}\label{average:trace}
\begin{split}
\|\{\!\{\nabla_h v_h\}\!\}\|_{0,e}
&\leqslant\frac{1}{2}\|(\nabla v_h)|_{K_e^1}\|_{0,e}
+\frac{1}{2}\|(\nabla v_h)|_{K_e^2}\|_{0,e}\\
&\leqslant \frac{\kappa_5}{2}\Big(
h_{K_e^1}^{-\frac{1}{2}}\|\nabla v_h\|_{0,K_e^1}+h_{K_e^2}^{-\frac{1}{2}}
\|\nabla v_h\|_{0,K_e^2}
\Big).
\end{split}
\end{equation}
Then for any $e\in\mathcal E_h^o$, it follows from the Cauchy-Schwarz inequality, \eqref{average:trace}, and the fact $|e|\leqslant h_{K_e^i}^{d-1}$, $i=1,2$, that
\begin{align*}
&\int_e\Big|
\{\!\{\nu\nabla_h v_h\}\!\}
\cdot\llbracket w \rrbracket\Big|\\
&\qquad\leqslant\nu|e|^{\frac{1}{2(d-1)}}\|
\{\!\{\nabla_h v_h\}\!\}\|_{0,e}\cdot
\bigg(\frac{1}{|e|^{\frac{1}{2(d-1)}}}
\|\llbracket w \rrbracket\|_{0,e}\bigg)\\
&\qquad\leqslant\nu\frac{\kappa_5}{2}
|e|^{\frac{1}{2(d-1)}}\Big(
h_{K_e^1}^{-\frac{1}{2}}\|\nabla v_h\|_{0,K_e^1}+h_{K_e^2}^{-\frac{1}{2}}
\|\nabla v_h\|_{0,K_e^2}
\Big)\cdot\bigg(\frac{1}{|e|^{\frac{1}{2(d-1)}}}
\|\llbracket w \rrbracket\|_{0,e}\bigg)\\
&\qquad\leqslant\nu\kappa_5
\Big(
\|\nabla v_h\|_{0,K_e^1}^2+
\|\nabla v_h\|_{0,K_e^2}^2
\Big)^{\frac{1}{2}}
\cdot
\bigg(\frac{1}{|e|^{\frac{1}{2(d-1)}}}
\|\llbracket w \rrbracket\|_{0,e}\bigg).
\end{align*}
Hence, using the Cauchy-Schwarz inequality again and the assumption (A1), we obtain
\begin{align*}
&\sum_{e\in\mathcal E_h}\int_e\Big|
\{\!\{\nu\nabla_h v_h\}\!\}
\cdot\llbracket w \rrbracket\Big|\\
&\qquad\leqslant\nu\kappa_5
\bigg(\sum_{e\in\mathcal E_h}\frac{1}{|e|^{1/(d-1)}}
\|\llbracket w \rrbracket\|_{0,e}^2\bigg)^{\frac{1}{2}}\cdot\bigg[
\sum_{e\in\mathcal E_h^o}
\Big(
\|\nabla v_h\|_{0,K_e^1}^2+\|\nabla v_h\|_{0,K_e^2}^2
\Big)+\sum_{e\in\Gamma}\|\nabla v\|_{0,K_e}^2
\bigg]^{\frac{1}{2}}\\
&\qquad\leqslant\nu\kappa_5\sqrt{\mathcal N_0}~|v_h|_{1,h}\bigg(\sum_{e\in\mathcal E_h}\frac{1}{|e|^{1/(d-1)}}
\|\llbracket w \rrbracket\|_{0,e}^2\bigg)^{\frac{1}{2}}.
\end{align*}
We choose $\kappa_9=\kappa_5\sqrt{\mathcal N_0}$, and the result \eqref{dg:inequality} then follows.
\end{pf}

\subsection{Proof of Lemma \ref{lemma:coercivity}}
\begin{pf}
We note that the reaction-convection part \eqref{rc:part} is the same as one shown in \cite{Ayuso20091391}. So we just prove \eqref{coer:diff} and \eqref{bound:norm}. For convenience, we set
\[
\|v\|_{p}=\bigg(\sum_{e\in\mathcal E_h}\frac{1}{|e|^{1/(d-1)}}
\|\llbracket v \rrbracket\|_{0,e}^2\bigg)^{\frac{1}{2}}.
\]
Following the same technique of deriving \eqref{dg:inequality}, for any continuous and bounded vector-valued function $\bm{\varphi}$ defined on $\Omega$, we have
\begin{equation}\label{dg:inequality2}
\sum_{e\in\mathcal E_h}\int_e\big|\{\!\{\bm{\varphi}v_h\}\!\}
\cdot\llbracket w \rrbracket\big|\leqslant \varphi_0\kappa_9\|v_h\|_{0,\Omega}\|w\|_{p}
\quad\forall v_h\in V_h^k,~\forall w\in V(\mathcal T_h),
\end{equation}
where $\varphi_0$ is a constant such that $|\bm{\varphi}(\bm{x})|\leqslant\varphi_0$, $\forall \bm{x}\in\Omega$.
Then using the Cauchy-Schwarz inequality, \eqref{wc:property}, Lemma \ref{useful:inequality} and \eqref{dg:inequality2}, we obtain
\begin{equation}\label{coercivity:diffusion}
\begin{split}
a_h^{D}(v_h,\mu v_h)\geqslant&\nu
\Big[
(\mu_1+\delta)\big(|v_h|_{1,h}^2+\eta_0
\|v_h\|_{p}^2\big)-2(\mu_2+\delta)\kappa_9
|v_h|_{1,h}\|v_h\|_{p}\\
&\quad-\mu_3\kappa_{10}\|v_h\|_{0,\Omega}\big(
|v_h|_{1,h}+\|v_h\|_{p}\big)
\Big]\\
=&\nu(I_1-I_2-I_3).
\end{split}
\end{equation}
It follows \eqref{wc:property2}, \eqref{assum:eta0} and the Young's inequality that
\begin{equation}\label{I2:estimate}
I_2\leqslant
\frac{2}{\eta_2\theta_0}(\mu_1+\delta)\big(|v_h|_{1,h}^2+
\eta_0
\|v_h\|_{p}^2\big).
\end{equation}
Then using the Poincar\'{e}--Friedrichs inequality \eqref{Poincare:dg}, the Young's inequality again and finally \eqref{wc:property2}, we have
\begin{equation}\label{I3:estimate}
\begin{split}
I_3
&\leqslant\mu_3\kappa_{10}LC_P\big(|v_h|_{1,h}^2
+\|v_h\|_{p}^2\big)^{\frac{1}{2}}\big(|v_h|_{1,h}
+\|v_h\|_{p}\big)\\
&\leqslant\frac{\sqrt{2}}{\eta_1}(\mu_1+\delta)\big(|v_h|_{1,h}^2+
\|v_h\|_{p}^2\big).
\end{split}
\end{equation}
Collecting \eqref{coercivity:diffusion}--\eqref{I3:estimate} with $\eta_0\geqslant 1$ yields \eqref{coer:diff}.

Now we turn our attention to \eqref{bound:norm}. According to the Young's inequality, \eqref{Poincare:dg} and \eqref{wc:property2}, it is a simple matter to check that
\begin{align*}
\normdg{\mu v_h}^2&=\nu|\mu v_h|_{1,h}^2+\nu\|\mu v_h\|_{p}^2+\normdg{\mu v_h}_{RC}^2\\
&\leqslant
2\mu_3^2\nu\|v_h\|_{0,\Omega}^2+
(\mu_2+\delta)^2\nu\big(2|v_h|_{1,h}^2+
\|\mu v_h\|_{p}^2\big)+(\mu_2+\delta)^2\normdg{v_h}_{RC}^2\\
&\leqslant
\bigg(\frac{2}{\eta_1^2\kappa_{10}^2}+
\frac{2}{\eta_2^2}\bigg)\big(
\normdg{v_h}_{D}^2+\normdg{v_h}_{RC}^2
\big),
\end{align*}
which leads to the result \eqref{bound:norm}.
\end{pf}

\subsection{Proof of Lemma \ref{super:appro}}
\begin{pf}
Because the proof of \eqref{project:H1} exactly follows the process of proving \eqref{project:L2}, we shall just prove \eqref{project:L2} and \eqref{project:trace}. Recalling $\mu=\overline{\xi}+\delta$ and $\delta$ is a constant, we have
\begin{equation}\label{project:constant}
\mu v_h-\widetilde{\mu v_h}\equiv \overline{\xi}v_h-\widetilde{\overline{\xi}v_h}.
\end{equation}
We note that $v_h|_K\in P_k(K)$, so $\partial^{\bm{\alpha}} v_h=0$ for $|\bm{\alpha}|=k+1$. In addition, $C^\ast$ is denoted as the maximum combinatorial number of $k+1$. Then for any $K\in\mathcal T_h$, it follows from \eqref{polygon:appro}, \eqref{project:constant}, \eqref{Linfty:norm}, \eqref{polygon:inverse}, and the fact $h_K<L$ that
\begin{equation}\label{temp:element}
\begin{split}
&\|\mu v_h-\widetilde{\mu v_h}\|_{0,K}\leqslant
\kappa_2 h_K^{k+1}\big|\overline{\xi} v_h\big|_{k+1,K}\\
&\quad\leqslant \kappa_2 C^\ast h_K^{k+1}
\bigg(\sum_{j=0}^kL^{k+1-j}\big|\overline{\xi}\big|_{k+1-j,\infty,K}\bigg)
\cdot\bigg(\sum_{j=0}^k\frac{|v_h|_{j,K}}{L^{k+1-j}}\bigg)\\
&\quad\leqslant \kappa_2\kappa_3C^\ast \frac{\big\|\overline{\xi}\big\|_{k+1,\infty,\Omega}}{L}\cdot
\sum_{j=0}^k\frac{h_K^{k+1-j}\|v_h\|_{0,K}}{L^{k-j}}\\
&\quad\leqslant \kappa_2\kappa_3C^\ast (k+1)
\frac{\big\|\overline{\xi}\big\|_{k+1,\infty,\Omega}}{L}
h_K\|v_h\|_{0,K}.
\end{split}
\end{equation}
After choosing $\kappa_{14}=\kappa_2\kappa_3C^\ast (k+1)$ and summing over all elements $K\in\mathcal T_h$, we obtain the inequality \eqref{project:L2}.

To prove \eqref{project:trace}, for any $K\in\mathcal T_h$ and any $e\subset \partial K$, using \eqref{re:polygontr} and \eqref{temp:element}, we have
\begin{equation}\label{temp:trace}
\begin{split}
\|\mu v_h-\widetilde{\mu v_h}\|_{0,e}&\leqslant \kappa_1^\ast
\frac{\big\|\overline{\xi}\big\|_{k+1,\infty,\Omega}}{L}
|e|^{\frac{1}{2(d-1)}}\big(\kappa_{14} \|v_h\|_{0,K}+\kappa_{15} \|v_h\|_{0,K}\big)\\
&\leqslant
\kappa_1^\ast(\kappa_{14}+\kappa_{15})
\frac{\big\|\overline{\xi}\big\|_{k+1,\infty,\Omega}}{L}
|e|^{\frac{1}{2(d-1)}}\|v_h\|_{0,K}.
\end{split}
\end{equation}
Hence, if we use the fact $|e|\leqslant h_K^{d-1}$, sum over all edges $e\in\mathcal E_h$ for \eqref{temp:trace}, and choose $\kappa_{16}=\kappa_1^\ast(\kappa_{14}+\kappa_{15})
\sqrt{\mathcal N_0}$, the inequality \eqref{project:trace} then follows.
\end{pf}

\subsection{Proof of Lemma \ref{lemma:contin}}
\begin{pf}
First, using the \eqref{temp:trace} in the process of proving \eqref{project:trace}, we have
\begin{equation}\label{temp:tracenorm}
\|\mu v_h-\widetilde{\mu v_h}\|_{p}= \bigg(\sum_{e\in\mathcal E_h}\frac{1}{|e|^{1/(d-1)}}
\big\|\llbracket \mu v_h-\widetilde{\mu v_h} \rrbracket\big\|_{0,e}^2\bigg)^{\frac{1}{2}}\leqslant
\kappa_{16}
\frac{\big\|\overline{\xi}\big\|_{k+1,\infty,\Omega}}{L}
\|v_h\|_{0,\Omega}.
\end{equation}
Then, it follows from Lemma \ref{useful:inequality}, \eqref{project:constant}, \eqref{polygon:appro}, \eqref{re:polygontr} and \eqref{temp:element} that for any $e\in \mathcal E_h^o$,
\begin{equation}\label{no:inverse:inequa}
\begin{split}
\big\|\{\!\{\nabla_h(\mu v_h-\widetilde{\mu v_h})\}\!\}\big\|_{0,e}&\leqslant
\frac{1}{2}\widetilde{\kappa}_1\kappa_2
\Big(h_{K_e^1}^{k-\frac{1}{2}}
\big|\overline{\xi} v_h\big|_{k+1,K_e^1}
+h_{K_e^2}^{k-\frac{1}{2}}
\big|\overline{\xi} v_h\big|_{k+1,K_e^2}\Big)\\
&\leqslant
\frac{1}{2}\widetilde{\kappa}_1\kappa_{14}
\frac{\big\|\overline{\xi}\big\|_{k+1,\infty,\Omega}}{L}
\Big(h_{K_e^1}^{-\frac{1}{2}}\|v_h\|_{0,K_e^1}
+h_{K_e^2}^{-\frac{1}{2}}\|v_h\|_{0,K_e^2}\Big).
\end{split}
\end{equation}
Following the technique to derive \eqref{dg:inequality} and replacing \eqref{average:trace} by \eqref{no:inverse:inequa}, we have
\begin{equation}\label{part:diffusion:contiu}
\sum_{e\in\mathcal E_h}\int_e\Big|
\{\!\{\nu\nabla_h(\mu v_h-\widetilde{\mu v_h})\}\!\}
\cdot\llbracket v_h \rrbracket\Big|
\leqslant
\nu\widetilde{\kappa}_1\kappa_{14}
\sqrt{\mathcal N_0}
\frac{\big\|\overline{\xi}\big\|_{k+1,\infty,\Omega}}{L}
\|v_h\|_{0,\Omega}\|v_h\|_{p}.
\end{equation}
Collecting \eqref{dg:inequality}, \eqref{project:H1}, \eqref{temp:tracenorm}, \eqref{part:diffusion:contiu}, the Cauchy-Schwarz inequality and \eqref{Poincare:dg}, we obtain
\begin{align*}
&\big|a_h^{D}(v_h,\mu v_h-\widetilde{\mu v_h})\big|\\
&\quad\leqslant
\nu|v_h|_{1,h}|\mu v_h-\widetilde{\mu v_h}|_{1,h}+
c_0\nu\|\mu v_h-\widetilde{\mu v_h}\|_{p}\|v_h\|_{p}\\
&\qquad+\nu\kappa_9
|v_h|_{1,h}\|\mu v_h-\widetilde{\mu v_h}\|_{p}+\nu\widetilde{\kappa}_1\kappa_{14}
\sqrt{\mathcal N_0}
\frac{\big\|\overline{\xi}\big\|_{k+1,\infty,\Omega}}{L}
\|v_h\|_{0,\Omega}\|v_h\|_{p}\\
&\quad\leqslant
\nu\frac{\big\|\overline{\xi}\big\|_{k+1,\infty,\Omega}}{L}
\|v_h\|_{0,\Omega}
\big(\kappa_{15}|v_h|_{1,h}
+c_0\kappa_{16}\|v_h\|_{p}+\kappa_9\kappa_{16}|v_h|_{1,h}
+\widetilde{\kappa}_1\kappa_{14}
\sqrt{\mathcal N_0}\|v_h\|_{p}\big)\\
&\quad\leqslant
2\sqrt{2}\max\big\{\kappa_{15}+\kappa_9\kappa_{16},
c_0\kappa_{16}+\widetilde{\kappa}_1\kappa_{14}\sqrt{\mathcal N_0}\big\}C_P\big\|\overline{\xi}\big\|_{k+1,\infty,\Omega}
\normdg{v_h}_{D}^2.
\end{align*}
Hence, if we choose
\[
\kappa_{17}=2\sqrt{2}\max\big\{\kappa_{15}+\kappa_9\kappa_{16},
c_0\kappa_{16}+\widetilde{\kappa}_1\kappa_{14}\sqrt{\mathcal N_0}\big\}C_P\big\|\overline{\xi}\big\|_{k+1,\infty,\Omega},
\]
the result \eqref{quasi:contin:df} then follows.

The proof of \eqref{quasi:contin:rc} has been shown in Lemma 4.3 of \cite{Ayuso20091391}, but here the right-hand constant is
\begin{equation}\label{kappa18}
\kappa_{18}=
\kappa_{14} C_r+\frac{\kappa_{14}}{2C_b}+\frac{\kappa_2\kappa_3\kappa_{14}}{C_b}+2\kappa_{16}.
\end{equation}
\end{pf}

\subsection{Proof of Theorem \ref{first:stability}}
\begin{pf}
For any $v_h\in V_h^k$, we set $w_h=\widetilde{\mu v_h}\in V_h^k$. It is so difficult to prove \eqref{first:infsup} directly that we shall prove the following equivalent results:
\begin{align}
\label{first:infsup1}
\mathcal A_h(v_h,w_h)&\geqslant \alpha_{11}\normdg{v_h}^2,\\
\label{first:infsup2}
\normdg{w_h}&\leqslant \alpha_{12}\normdg{v_h},
\end{align}
where $\alpha_{11}$, $\alpha_{12}$ are two positive constants.

We first prove \eqref{first:infsup1}. It follows from \eqref{coer:diff}, \eqref{quasi:contin:df} and the case $\kappa_{11}(\mu_1+\delta)\geqslant 2\kappa_{17}$ that
\begin{equation}\label{first:infsup1:temp1}
\begin{split}
a_h^{D}(v_h,\widetilde{\mu v_h})&=a_h^{D}(v_h,\widetilde{\mu v_h}-\mu v_h)+a_h^{D}(v_h,\mu v_h)\\
&\geqslant \kappa_{11}(\mu_1+\delta)\normdg{v_h}_{D}^2-\kappa_{17}\normdg{v_h}_{D}^2\\
&\geqslant \kappa_{17}\normdg{v_h}_{D}^2.
\end{split}
\end{equation}
Let $h_0$ be a positive constant satisfying $h_0<\big(\frac{\mu_1}{2\kappa_{18}}\big)^2 L$.
In the same way, using \eqref{coer:rc} and \eqref{quasi:contin:rc}, we obtain that for $h<h_0$,
\begin{equation}\label{first:infsup1:temp2}
\begin{split}
a_h^{RC}(v_h,\widetilde{\mu v_h})&=a_h^{RC}(v_h,\widetilde{\mu v_h}-\mu v_h)+a_h^{RC}(v_h,\mu v_h)\\
&\geqslant \frac{\mu_1}{2}\normdg{v_h}_{RC}^2-\kappa_{18}\bigg(\frac{h}{L}
\bigg)^{\frac{1}{2}}\normdg{v_h}_{RC}^2\\
&\geqslant \frac{\mu_1-2\kappa_{18}\sqrt{h_0}}{2\sqrt{L}}\normdg{v_h}_{RC}^2.
\end{split}
\end{equation}
Hence, if choosing $\alpha_{11}=\max\big\{\kappa_{17},\frac{\mu_1-2\kappa_{18}\sqrt{h_0}}{2\sqrt{L}}\big\}$, and collecting \eqref{first:infsup1:temp1}--\eqref{first:infsup1:temp2}, we have the first result \eqref{first:infsup1}.

Next, based on \eqref{bound:norm}, Lemma \ref{super:appro} and $h_K\leqslant h<L$, we have
\[
\normdg{\widetilde{\mu v_h}}
\leqslant\normdg{\widetilde{\mu v_h}-\mu v_h}
+\normdg{\mu v_h},
\]
and
\begin{align*}
\normdg{\widetilde{\mu v_h}-\mu v_h}^2&=\normdg{\widetilde{\mu v_h}-\mu v_h}_{D}^2
+\normdg{\widetilde{\mu v_h}-\mu v_h}_{RC}^2\\
&\leqslant\big\|\overline{\xi}\big\|_{k+1,\infty,\Omega}^2
\Big[\big(\kappa_{15}^2+\kappa_{16}^2\big)C_P^2
\normdg{v_h}_{D}^2+\big(\kappa_{14}^2+2\kappa_{16}^2\big)
\|(\overline{r}+b_0)^{\frac{1}{2}}v_h\|_{0,\Omega}^2\Big]\\
&\leqslant\big\|\overline{\xi}\big\|_{k+1,\infty,\Omega}^2
\max\big\{\big(\kappa_{15}^2+\kappa_{16}^2\big)C_P^2,\kappa_{14}^2+2\kappa_{16}^2\big\}
\normdg{v_h}^2.
\end{align*}
So, choosing
\[
\alpha_{12}=\big\|\overline{\xi}\big\|_{k+1,\infty,\Omega}
\sqrt{\max\big\{\big(\kappa_{15}^2+\kappa_{16}^2\big)C_P^2,\kappa_{14}^2+2\kappa_{16}^2\big\}}
+\kappa_{13}(\mu_1+\delta),
\]
we then obtain the result \eqref{first:infsup2}.

Further from the definitions of $\alpha_{11}$ and $\alpha_{12}$, it implies that $\alpha_1$ depends only on the velocity field $\bm{b}$. In addition, we shall provide enough reasons to claim that the condition $h<h_0$ is necessary. In \eqref{first:infsup1:temp2}, if using $h<L$ instead of $h<h_0$, we then obtain
\[
a_h^{rc}(v_h,\widetilde{\mu v_h})\geqslant \frac{\mu_1-2\kappa_{18}}{2}\normdg{v_h}_{rc}^2.
\]
So clearly, it requires $\mu_1>2\kappa_{18}$. However, from (H1) we note that the value of $\mu_1$ that depends only on the velocity field $\bm{b}$. Once the model problem \eqref{p1}--\eqref{p2} is given, the upper bound of $\mu_1$ is determinate, but from \eqref{kappa18} the value of $\kappa_{18}$ can not change either. It implies that the condition $\mu_1>2\kappa_{18}$ will be not always satisfied. Fortunately, \eqref{quasi:contin:rc} and the condition $h<h_0$ ensure that the result \eqref{first:infsup1:temp2} is always correct. Hence, the condition $h<h_0$ is necessary.
\end{pf}

\subsection{Proof of Lemma \ref{lemma:supg:preinfsup}}
\begin{pf}
For any $v_h\in V_h^k$, applying the reconstruction operator $\mathcal P^k$ defined in \eqref{global:op} to the piecewise term $\bm{b}\cdot\nabla_h v_h$, we know that $\mathcal P^k(\bm{b}\cdot\nabla_h v_h)\in V_h^k$. Then we set
\[
w_h=\sum_{K\in \mathcal T_h}c_K\big(\mathcal P^k(\bm{b}\cdot\nabla_h v_h)\big)|_K,
\]
where
\begin{equation*}
c_K=
\begin{cases}
\frac{h_K}{\|\bm{b}\|_{0,\infty,K}} &\text{if~convection~dominates~in}~K,\\
0 & \text{otherwise}.
\end{cases}
\end{equation*}
Following the technique of proving \eqref{first:infsup}, if the convection dominates, we shall prove that
\begin{align}
\label{supg:preinfsup1}
\normdg{w_h}&\leqslant \alpha_{21}\|v_h\|_b,\\
\label{supg:preinfsup2}
\mathcal A_h(v_h,w_h)&\geqslant \|v_h\|_b^2-\alpha_{22}\normdg{v_h}\|v_h\|_b,
\end{align}
where $\alpha_{21}$, $\alpha_{22}$ are two positive constants.

For \eqref{supg:preinfsup1}, the local inverse inequality \eqref{polygon:inverse} and two versions of the local inverse trace inequality \eqref{inverse:trace}, \eqref{inverse:trace2} will be extensively employed. It follows from \eqref{polygon:inverse} and \eqref{convection:dominate:assum} that
\begin{equation}\label{supg:preinfsup1:temp1}
\begin{split}
&\nu |w_h|_{1,h}^2=\sum_{K\in\mathcal T_h}\nu\bigg(
\frac{h_K}{\|\bm{b}\|_{0,\infty,K}}\bigg)^2|
\mathcal P^k(\bm{b}\cdot\nabla_h v_h)|_{1,K}^2\\
&\leqslant\sum_{K\in\mathcal T_h}
\frac{h_K\|\bm{b}\|_{0,\infty,K}}{2}\cdot
\frac{h_K^2}{\|\bm{b}\|_{0,\infty,K}^2}
\cdot \frac{\kappa_3^2}{h_K^2}\|\mathcal P^k(\bm{b}\cdot\nabla_h v_h)\|_{0,K}^2\leqslant \frac{\kappa_3^2}{2}\|v_h\|_b^2.
\end{split}
\end{equation}
Then, using \eqref{inverse:trace}, \eqref{polygon:area} and \eqref{convection:dominate:assum}, we have
\begin{align}
&\nu \|w_h\|_{p}^2=\sum_{e\in\mathcal E_h}\frac{\nu}{|e|^{1/(d-1)}}\big\|\llbracket c_K\mathcal P^k(\bm{b}\cdot\nabla_h v_h)\rrbracket\big\|_{0,e}^2\notag\\
\leqslant&\kappa_4^2C_g\bigg[\sum_{e\in\mathcal E_h^o}
\Big(c_{K_e^1}\|\mathcal P^k(\bm{b}\cdot\nabla_h v_h)\|_{0,K_e^1}^2+
c_{K_e^2}\|\mathcal P^k(\bm{b}\cdot\nabla_h v_h)\|_{0,K_e^2}^2\Big)+\sum_{e\in\Gamma}c_{K_e}
\|\mathcal P^k(\bm{b}\cdot\nabla_h v_h)\|_{0,K_e}^2
\bigg]\notag\\
\label{supg:preinfsup1:temp2}
\leqslant&2\kappa_5^2\mathcal N_0\|v_h\|_b^2.
\end{align}
Similarly, from \eqref{inverse:trace2} and \eqref{convection:dominate:assum}, we deduce that
\begin{equation}\label{supg:preinfsup1:temp3}
\sum_{e\in\mathcal E_h}\big\||\bm{b}\cdot\bm{n}|^{\frac{1}{2}}\llbracket w_h
\rrbracket\big\|_{0,e}^2
=\sum_{e\in\mathcal E_h}
\big\||\bm{b}\cdot\bm{n}|^{\frac{1}{2}}\llbracket c_K\mathcal P^k(\bm{b}\cdot\nabla_h v_h)\rrbracket\big\|_{0,e}^2
\leqslant
2\kappa_5^2\mathcal N_0\|v_h\|_b^2.
\end{equation}
Following the technique to estimate the term $\|(\overline{r}+b_0)^{\frac{1}{2}}w_h\|_{0,\Omega}$ in Lemma 4.5 of \cite{Ayuso20091391}, we get
\begin{equation}\label{supg:preinfsup1:temp4}
\|(\overline{r}+b_0)^{\frac{1}{2}}w_h\|_{0,\Omega}^2
=\|(\overline{r}+b_0)^{\frac{1}{2}}
\mathcal P^k(\bm{\omega}\cdot\nabla_h v_h)\|_{0,\Omega}^2
\leqslant \bigg(1+\frac{3}{2C_b}\bigg)
\|v_h\|_b^2.
\end{equation}
Finally, collecting \eqref{supg:preinfsup1:temp1}--\eqref{supg:preinfsup1:temp4} and let
\begin{equation}\label{alpha:21}
\alpha_{21}=\sqrt{\frac{\kappa_3^2}{2}
+\frac{3}{2C_b}+4\kappa_5^2\mathcal N_0+1},
\end{equation}
we obtain the first result \eqref{supg:preinfsup1}.

Next, we shall prove \eqref{supg:preinfsup2}. The rest proof is almost similar to the corresponding part in Lemma 4.5 of \cite{Ayuso20091391}, but instead of always using the simple general constant $C$ in \cite{Ayuso20091391}, we shall distinguish each constant in all inequalities because of our general polytopic meshes. For the diffusion term $a_h^{D}(v_h,w_h)$, it follows the Cauchy-Schwarz inequality,  \eqref{supg:preinfsup1:temp1} and \eqref{supg:preinfsup1:temp2} that
\[
\int_\Omega\big|\nu
\nabla_h v_h\cdot\nabla_h w_h\big|
\leqslant\nu^{\frac{1}{2}}|v_h|_{1,h}
\cdot\nu^{\frac{1}{2}}|w_h|_{1,h}
\leqslant\frac{\kappa_3}{\sqrt{2}}
\nu^{\frac{1}{2}}|v_h|_{1,h}\|v_h\|_b,
\]
and
\begin{equation*}
\sum_{e\in\mathcal E_h}\frac{\sigma_e\nu}{|e|^{1/(d-1)}}
\int_e\big|\llbracket v_h\rrbracket
\cdot\llbracket w_h\rrbracket\big|
\leqslant c_0\nu^{\frac{1}{2}}\|v_h\|_{p}
\cdot\nu^{\frac{1}{2}}\|w_h\|_{p}
\leqslant\kappa_5c_0\sqrt{\mathcal N_0}
\nu^{\frac{1}{2}}\|v_h\|_{p}\|v_h\|_b.
\end{equation*}
In addition, from Lemma \ref{useful:inequality}, \eqref{supg:preinfsup1:temp1} and \eqref{supg:preinfsup1:temp2}, we have
\begin{align*}
&\sum_{e\in\mathcal E_h}\int_e
\big|\{\!\{\nu\nabla_h v_h\}\!\}
\cdot\llbracket w_h \rrbracket\big|+\int_e\big|\{\!\{\nu\nabla_h w_h\}\!\}
\cdot\llbracket v_h \rrbracket\big|\\
&\leqslant\kappa_9\nu^{\frac{1}{2}}
|v_h|_{1,h}\cdot\nu^{\frac{1}{2}}\|w_h\|_{p}
+\kappa_9\nu^{\frac{1}{2}}
|w_h|_{1,h}\cdot\nu^{\frac{1}{2}}\|v_h\|_{p}\\
&\leqslant
\bigg(
\kappa_5\kappa_9\sqrt{\mathcal N_0}\nu^{\frac{1}{2}}
|v_h|_{1,h}
+\frac{\kappa_3\kappa_9}{\sqrt{2}}
\nu^{\frac{1}{2}}\|v_h\|_{p}
\bigg)\|v_h\|_b.
\end{align*}
Hence, it follows from the Cauchy-Schwarz inequality that
\begin{equation}\label{supg:preinfsup2:temp1}
a_h^{D}(v_h,w_h)\geqslant
-\widetilde{\alpha}_{22}\normdg{v_h}\|v_h\|_b,
\end{equation}
where
\begin{equation}\label{alpha:tilde:22}
\widetilde{\alpha}_{22}=
\bigg[
\Big(\frac{\kappa_9}{\sqrt{2}}+\kappa_5\kappa_9
\sqrt{\mathcal N_0}\Big)^2
+\Big(\frac{\kappa_3\kappa_9}{\sqrt{2}}
+\kappa_5c_0
\sqrt{\mathcal N_0}\Big)^2
\bigg]^{\frac{1}{2}}.
\end{equation}
On the other hand, for the reaction-convection term $a_h^{RC}(v_h,w_h)$, we first note that the $L^2$-projection property of $\mathcal P^k$ mentioned in \eqref{poly:projection} implies
\begin{equation}\label{supg:preinfsup2:temp2}
\int_\Omega(\bm{b}\cdot\nabla_h v_h)w_h
=\sum_{K\in\mathcal T_h}
\int_K \mathcal P^k(\bm{b}\cdot\nabla v_h)w_h
=\|v_h\|_b^2.
\end{equation}
Then using integration by parts, \eqref{magic:formula} and \eqref{upwind:value}, we shall rewrite $a_h^{RC}(v_h,w_h)$. To fix ideas, collecting
\begin{align*}
&-\int_\Omega v_h\bm{b}\cdot\nabla_h w_h
=\int_\Omega (\bm{b}\cdot\nabla_h v_h)w_h+
\int_\Omega \nabla\cdot\bm{b}v_hw_h
-\sum_{K\in\mathcal T_h}\int_{\partial K}
(\bm{b}\cdot\bm{n})v_hw_h,\\
&
-\sum_{K\in\mathcal T_h}\int_{\partial K}
(\bm{b}\cdot\bm{n})v_hw_h
=-\sum_{e\in\mathcal E_h}\int_e
\{\!\{\bm{b}v_h\}\!\}\cdot\llbracket w_h \rrbracket-\sum_{e\in\mathcal E_h^o}\int_e\llbracket \bm{b}v_h \rrbracket
\{\!\{w_h\}\!\},\\
&\sum_{e\in\mathcal E_h\backslash\Gamma^-}
\int_e\{\!\{\bm{b}v_h\}\!\}_{up}\cdot
\llbracket w_h \rrbracket
=\sum_{e\in\mathcal E_h^o}
\int_e\{\!\{\bm{b}v_h\}\!\}_{up}\cdot
\llbracket w_h \rrbracket+
\sum_{e\in\Gamma^+}
\int_e(\bm{b}\cdot\bm{n})v_hw_h,
\end{align*}
and
\[
\sum_{e\in\mathcal E_h^o}
\int_e\{\!\{\bm{b}v_h\}\!\}_{up}\cdot
\llbracket w_h \rrbracket
=\sum_{e\in\mathcal E_h^o}
\int_e\{\!\{\bm{b}v_h\}\!\}\cdot
\llbracket w_h \rrbracket+
\sum_{e\in\mathcal E_h^o}
\int\frac{\bm{b}\cdot\bm{n}^+}{2}
\llbracket v_h \rrbracket\cdot
\llbracket w_h \rrbracket,
\]
we rewrite $a_h^{RC}(v_h,w_h)$ as:
\begin{equation*}
\begin{split}
a_h^{RC}(v_h,w_h)&=\int_\Omega r v_h w_h+
\int_\Omega (\bm{b}\cdot\nabla_h v_h)w_h+\frac{1}{2}
\int_\Omega \nabla\cdot\bm{b}v_hw_h\\
&\quad+\sum_{e\in\mathcal E_h^o}
\int\frac{\bm{b}\cdot\bm{n}^+}{2}
\llbracket v_h \rrbracket\cdot
\llbracket w_h \rrbracket
-\sum_{e\in\mathcal E_h^o\cup\Gamma^{-}}\int_e\llbracket \bm{b}v_h \rrbracket
\{\!\{w_h\}\!\}.
\end{split}
\end{equation*}
It follows from $h<L$, (H2), (H3), \eqref{velocity:assumption}, \eqref{supg:preinfsup1:temp3} and \eqref{supg:preinfsup1:temp4} that
\[
\int_\Omega\big|r v_h w_h\big|
\leqslant C_r\sqrt{1+\frac{3}{2C_b}}
~\big\|(\overline{r}+b_0)^{\frac{1}{2}}
v_h\big\|_{0,\Omega}\big\|v_h\big\|_b,
\]
\[
\int_\Omega\big|\nabla\cdot\bm{b}v_hw_h
\big|\leqslant\frac{1}{C_b}
\big\|(\overline{r}+b_0)^{\frac{1}{2}}
v_h\big\|_{0,\Omega}\big\|v_h\big\|_b,
\]
and
\[
\bigg|\sum_{e\in\mathcal E_h^o}
\int\frac{\bm{b}\cdot\bm{n}^+}{2}
\llbracket v_h \rrbracket\cdot
\llbracket w_h \rrbracket
-\sum_{e\in\mathcal E_h^o\cup\Gamma^{-}}\int_e\llbracket \bm{b}v_h \rrbracket
\{\!\{w_h\}\!\}\bigg|\leqslant 2\kappa_5\sqrt{\mathcal N_0}\sum_{e\in\mathcal E_h}\Big(\big\||\bm{b}\cdot\bm{n}|^{\frac{1}{2}}
\llbracket v_h
\rrbracket\big\|_{0,e}^2\Big)^{\frac{1}{2}}
\big\|v_h\big\|_b.
\]
Hence, using the Cauchy-Schwarz inequality, \eqref{supg:preinfsup2:temp2} and the above inequalities, we obtain
\begin{equation}\label{supg:preinfsup2:temp3}
a_h^{RC}(v_h,w_h)\geqslant \|v_h\|_b^2
-\overline{\alpha}_{22}\normdg{v_h}\|v_h\|_b,
\end{equation}
where
\begin{equation}\label{alpha:22}
\overline{\alpha}_{22}
=\bigg[\bigg(
C_r\sqrt{1+\frac{3}{2C_b}}
+\frac{1}{2C_b}
\bigg)^2+4\kappa_5^2\mathcal N_0
\bigg]^{\frac{1}{2}}.
\end{equation}
Let $\alpha_{22}=\widetilde{\alpha}_{22}
+\overline{\alpha}_{22}$, and then we have \eqref{supg:preinfsup2} from \eqref{supg:preinfsup2:temp1} and \eqref{supg:preinfsup2:temp3}. Finally, \eqref{supg:preinfsup1} and \eqref{supg:preinfsup2} give \eqref{supg:preinfsup} with $\alpha_2=1/\alpha_{21}$ and
$\alpha_3=\alpha_{22}/\alpha_{21}$.

Further \eqref{alpha:21}, \eqref{alpha:tilde:22} and \eqref{alpha:22} ensure that both $\alpha_2$ and $\alpha_3$ are independent of $h$, $\nu$, $\bm{b}$, and $c$.
\end{pf}

\subsection{Proof of Theorem \ref{second:stability}}
\begin{pf}
According to Theorem \ref{first:stability} and Lemma \ref{lemma:supg:preinfsup}, for any $v_h\in V_h^k$, we have
\begin{equation}\label{first:infsup:var}
\frac{1}{\alpha_1}\cdot\sup_{w_h\in V_h^k}\frac{\mathcal A_h(v_h, w_h)}{\normdg{w_h}}\geqslant
\normdg{v_h},
\end{equation}
and
\begin{equation}\label{sup:preinfsup:var}
\frac{1}{\alpha_3}\cdot\sup_{w_h\in V_h^k}\frac{\mathcal A_h(v_h, w_h)}{\normdg{w_h}}\geqslant
\frac{\alpha_2}{\alpha_3}\|v_h\|_b
-\normdg{v_h}.
\end{equation}
Then, by adding \eqref{first:infsup:var} and \eqref{sup:preinfsup:var}, we obtain
\[
\bigg(\frac{\alpha_3}{\alpha_1\alpha_2}
+\frac{1}{\alpha_2}\bigg)\cdot\sup_{w_h\in V_h^k}\frac{\mathcal A_h(v_h, w_h)}{\normdg{w_h}}\geqslant
\|v_h\|_b,
\]
which, together with \eqref{first:infsup:var}, gives \eqref{second:infsup} with
\[
\alpha_4=\frac{\sqrt{(\alpha_1+\alpha_3)^2+
\alpha_2^2}}{\alpha_1\alpha_2}.
\]
\end{pf}

\section{Proofs of all lemmas and theorems in Section \ref{prior:error}}\label{append:pf:4}

\subsection{Proof of Theorem \ref{priori:dg}}
\label{append:pf:th41}
\begin{pf}
Define
\[
\zeta_1=u-\mathcal P^k u,\quad \zeta_2=\mathcal P^k u_h-\mathcal P^k u,
\]
and we then have $u-\mathcal P^k u_h=\zeta_1-\zeta_2$. Noting that $\zeta_2\in V_h^k$, and then from Theorem \ref{first:stability} and the consistency \eqref{con:dg}, we have
\begin{equation}\label{stability:app}
\alpha_1\normdg{\zeta_2}
\leqslant\frac{a_h(\zeta_2, v_h)}{\normdg{v_h}}
=\frac{a_h(\zeta_1,v_h)}{\normdg{v_h}},\quad\forall
v_h\in V_h^k.
\end{equation}
In addition, for future purposes, it follows from our discussions of Section \ref{subsec:ep}, Lemma 6 of \cite{Li2019268} and Lemma 3.4 of \cite{Li2012259} that there exists a positive constant $\kappa_{19}$ depending only on $N$, $\gamma$, $k$ and $d$ such that for each element $K\in\mathcal T_h$, $\#S(K)$ has a uniform upper bound, i.e., $M\leqslant\kappa_{19}$. Next, we shall estimate $a_h^{D}(u-\mathcal P^k u,v_h)$ and $a_h^{RC}(u-\mathcal P^k u,v_h)$, respectively.

For the diffusion part $a_h^{D}(\zeta_1,v_h)$, using the Cauchy-Schwarz inequality, \eqref{dg:inequality}, \eqref{estimate:H1op} and \eqref{estimate:trace}, we deduce
\begin{align*}
a_h^{D}(\zeta_1,v_h)&\leqslant
(2c_0+1)\kappa_8\sqrt{\mathcal N_0\kappa_{19}}\nu^{\frac{1}{2}}h^k
|u|_{k+1,\Omega}\cdot\nu^{\frac{1}{2}}\|v_h\|_{p}\\
&\quad+\big(2c_0\kappa_9\sqrt{\mathcal N_0\kappa_{19}}+\kappa_7\sqrt{\kappa_{19}}\big)
\nu^{\frac{1}{2}}h^k
|u|_{k+1,\Omega}\cdot\nu^{\frac{1}{2}}|v_h|_{1,h}.
\end{align*}
Hence, $a_h^{D}(\zeta_1,v_h)$ can be estimated by
\begin{equation}\label{diffusion:error:estimate}
a_h^{D}(\zeta_1,v_h)\leqslant\alpha_5\nu^{\frac{1}{2}}h^k
|u|_{k+1,\Omega}\normdg{v_h}_{D},
\end{equation}
where
\[
\alpha_5=
\Big[(2c_0+1)^2\kappa_8^2\mathcal N_0\kappa_{19}+\big(2c_0\kappa_9\sqrt{\mathcal N_0\kappa_{19}}+\kappa_7\sqrt{\kappa_{19}}\big)^2
\Big]^{\frac{1}{2}}.
\]

For the reaction--convection part $a_h^{RC}(\zeta_1,v_h)$, let us first recall
\[
a_h^{RC}(\zeta_1,v_h)=\int_\Omega c\zeta_1v_h
-\int_\Omega(\bm{b}\cdot\nabla_h v_h)\zeta_1
+\sum_{e\in\mathcal E_h\cup \Gamma^+}\int_e\{\!\{\bm{b}\zeta_1\}\!\}
\cdot\llbracket v_h\rrbracket
+\sum_{e\in\mathcal E_h^o}\int_e\frac{\bm{b}\cdot\bm{n}^+}{2}\llbracket \zeta_1\rrbracket
\cdot\llbracket v_h\rrbracket,
\]
and the standard $L^2$-projection operator $\textsf{P}_h^k$ on $V_h^k$. Because $v_h=\mathcal P^k v\in V_h^k$, we then have
$\textsf{P}_h^0\bm{b}\cdot\nabla_hv_h\in V_h^k$. It follows from our projection property \eqref{poly:projection} that
\begin{equation}\label{rc:error:projection}
\int_\Omega
\big(\textsf{P}_h^0\bm{b}
\cdot\nabla_hv_h\big)\zeta_1=0.
\end{equation}
Hence, using \eqref{rc:error:projection},
the Cauchy-Schwarz inequality,
the classical interpolation results, \eqref{velocity:assumption}, \eqref{polygon:inverse} and \eqref{estimate:L2op}, we have
\begin{equation}\label{rc:error:estimate:temp1}
\begin{split}
&\int_\Omega-\big(\bm{b}\cdot\nabla_h v_h\big)\zeta_1
=\int_\Omega\big(\textsf{P}_h^0\bm{b}-\bm{b}\big)
\cdot\nabla_hv_h\zeta_1\\
&\leqslant C
\kappa_3\kappa_6\sqrt{\frac{\kappa_{19}}
{C_b}}\bigg(
\frac{\|\bm{b}\|_{0,\infty,\Omega}}
{L}\bigg)^{\frac{1}{2}}h^{k+1}|u|_{k+1,\Omega}
\normdg{v_h}_{RC},
\end{split}
\end{equation}
where $C$ is a standard interpolation constant. Then, using the Cauchy-Schwarz inequality and \eqref{estimate:trace}, we obtain
\begin{equation}\label{rc:error:estimate:temp2}
\sum_{e\in\mathcal E_h\cup \Gamma^+}\int_e\{\!\{\bm{b}\zeta_1\}\!\}
\cdot\llbracket v_h\rrbracket
\leqslant\kappa_8\sqrt{\mathcal N_0\kappa_{19}}
\|\bm{b}\|_{0,\infty,\Omega}^{\frac{1}{2}}h^{k+\frac{1}{2}}|u|_{k+1,\Omega}
\normdg{v_h}_{RC}.
\end{equation}
Similarly, we have
\begin{equation}\label{rc:error:estimate:temp3}
\sum_{e\in\mathcal E_h^o}\int_e\frac{\bm{b}\cdot\bm{n}^+}{2}\llbracket \zeta_1\rrbracket
\cdot\llbracket v_h\rrbracket
\leqslant\kappa_8\sqrt{\mathcal N_0\kappa_{19}}
\|\bm{b}\|_{0,\infty,\Omega}^{\frac{1}{2}}h^{k+\frac{1}{2}}|u|_{k+1,\Omega}
\normdg{v_h}_{RC}.
\end{equation}
It follows from \eqref{asm:r}, (H3), \eqref{velocity:assumption} and \eqref{estimate:L2op} that
\begin{equation}\label{rc:error:estimate:temp4}
\int_\Omega c\zeta_1v_h
\leqslant\big(C_r^{\frac{1}{2}}+C_b^{-1}\big)
\kappa_6\sqrt{\kappa_{19}}h^{k+1}\bigg[\|r\|_{0,\infty,\Omega}^{\frac{1}{2}}
+\bigg(\frac{\|\bm{\omega}\|_{0,\infty,\Omega}}
{L}\bigg)^{\frac{1}{2}}\bigg]|u|_{k+1,\Omega}\normdg{v_h}_{RC}.
\end{equation}
Finally, if we collect \eqref{rc:error:estimate:temp1}--\eqref{rc:error:estimate:temp4} and note that $h<L$, $a_h^{RC}(\zeta_1,v_h)$ can be estimated by
\begin{equation}\label{rc:error:estimate}
a_h^{RC}(\zeta_1,v_h)\leqslant\alpha_6\big(\|r\|_{0,\infty,\Omega}^{\frac{1}{2}}h
+\|\bm{b}\|_{0,\infty,\Omega}^{\frac{1}{2}}h^{\frac{1}{2}}\big)h^k
|u|_{k+1,\Omega}\normdg{v_h}_{RC},
\end{equation}
where
\[
\alpha_6=C
\kappa_3\kappa_6\sqrt{\frac{\kappa_{19}}
{C_b}}+2\kappa_8\sqrt{\mathcal N_0\kappa_{19}}+\big(C_r^{\frac{1}{2}}+C_b^{-1}\big)
\kappa_6\sqrt{\kappa_{19}}.
\]

Now from the estimates \eqref{diffusion:error:estimate} and \eqref{rc:error:estimate}, we obtain
\begin{equation}\label{total:error:estimate}
\mathcal A_h(\zeta_1,v_h)\leqslant(\alpha_5+\alpha_6)\big(\nu^{\frac{1}{2}}+
\|\bm{b}\|_{0,\infty,\Omega}^{\frac{1}{2}}h^{\frac{1}{2}}+
\|r\|_{0,\infty,\Omega}^{\frac{1}{2}}h
\big)h^k
|u|_{k+1,\Omega}\normdg{v_h}.
\end{equation}
Hence, substituting \eqref{total:error:estimate} into \eqref{stability:app} yields
\[
\normdg{\zeta_2}\leqslant
\frac{\alpha_5+\alpha_6}{\alpha_1}
\big(\nu^{\frac{1}{2}}+
\|\bm{b}\|_{0,\infty,\Omega}^{\frac{1}{2}}h^{\frac{1}{2}}+
\|r\|_{0,\infty,\Omega}^{\frac{1}{2}}h
\big)h^k
|u|_{k+1,\Omega},
\]
and the result \eqref{priori:dg:ineq} then follows by the triangle inequality.
\end{pf}

\subsection{Proof of Lemma \ref{auxiliary:regularity}}
\begin{pf}
Using the standard Poincar\'{e} inequality, we obtain that
\begin{equation}\label{con:poin}
\|\psi\|_{0,\Omega}\leqslant C(\Omega,d)\|\nabla\psi\|_{0,\Omega}.
\end{equation}
Then, it follows from \eqref{var:aux}, integration by parts, \eqref{asm:r} and \eqref{con:poin} that
\begin{equation}\label{temp:estimate:H1}
\nu\|\nabla\psi\|_{0,\Omega}
\leqslant C(\Omega,d)\|u-\mathcal P^k u_h\|_{0,\Omega}.
\end{equation}
Using \eqref{var:aux}, \eqref{con:poin} and \eqref{temp:estimate:H1}, we derive that
\begin{align*}
\nu\|\psi\|_{2,\Omega}
&\leqslant C(\Omega,d)(\|\bm{b}\|_{0,\infty,\Omega}
\|\nabla\psi\|_{0,\Omega}+\|r\|_{0,\infty,\Omega}
\|\nabla\psi\|_{0,\Omega}+\|u-\mathcal P^k u_h\|_{0,\Omega})\\
&\leqslant C(\Omega,d)(\nu^{-1}\|\bm{b}\|_{0,\infty,\Omega}
+\nu^{-1}\|r\|_{0,\infty,\Omega}+1)\|u-\mathcal P^k u_h\|_{0,\Omega},
\end{align*}
and \eqref{var:stable} then follows.
\end{pf}

\subsection{Proof of Theorem \ref{priori:l2}}
\label{append:pf:th44}
\begin{pf}
For any $v_h\in V_h^k$, it follows from \eqref{var:aux}, integration by parts, the fact that $\psi\in H^2(\Omega)$, and \eqref{con:dg} that
\begin{equation}\label{L2}
\begin{split}
\|u-\mathcal P^ku_h\|_{0,\Omega}^2&=(-\nu\Delta \psi+\bm{b}\cdot\nabla\psi+c\psi,u-\mathcal P^ku_h)_\Omega=\mathcal A_h(u-\mathcal P^ku_h,\psi)\\
&=\mathcal A_h(u-\mathcal P^ku_h,\psi-v_h)\\
&=a_h^{D}(u-\mathcal P^ku_h,\psi-v_h)+a_h^{RC}(u-\mathcal P^ku_h,\psi-v_h),\quad \forall v_h\in V_h^k.
\end{split}
\end{equation}
Let $\zeta=u-\mathcal P^ku_h$ and $v_h=\widetilde{\psi}$, where with respect to $\psi$, $\widetilde{\psi}|_K\in P_k(K)$ is an approximation polynomial satisfying \eqref{polygon:appro}.

For the diffusion part, using the Cauchy-Schwarz inequality, \eqref{polygon:appro} and \eqref{polygon:tr}, we obtain
\begin{equation}\label{ahd}
\begin{split}
&a_h^{D}(u-\mathcal P^ku_h,\psi-v_h)
=a_h^{D}(\zeta,\psi-\widetilde{\psi})\\
=&\nu\sum_{K\in\mathcal T_h}
\int_K\nabla\zeta\cdot\nabla(\psi-\widetilde{\psi})
-\sum_{e\in\mathcal E_h}\int_e\{\!\{\nu\nabla(\psi-\widetilde{\psi})\}\!\}\cdot
\llbracket \zeta\rrbracket\\
\leqslant&\big(\kappa_2+\kappa_1\kappa_2\sqrt{C_g\mathcal N_0}\big)\nu^{\frac{1}{2}}h\|\psi\|_{2,\Omega}\normdg{\zeta}.
\end{split}
\end{equation}

For the reaction-convection part, following the technique to rewrite $a_h^{RC}(v_h,w_h)$ in Lemma \ref{lemma:supg:preinfsup}, and using the fact that $\bm{b}$ is divergence free, the Cauchy-Schwarz inequality, \eqref{polygon:appro} and \eqref{polygon:tr}, we obtain
\begin{equation}\label{ahrc}
\begin{split}
&a_h^{RC}(u-\mathcal P^ku_h,\psi-v_h)
=a_h^{RC}(\zeta,\psi-\widetilde{\psi})\\
=&\int_\Omega r\zeta(\psi-\widetilde{\psi})
+\int_\Omega(\bm{b}\cdot\nabla_h \zeta)(\psi-\widetilde{\psi})
-\sum_{e\in\mathcal E_h^o\cup\Gamma^-}\int_e\llbracket \bm{b}\zeta\rrbracket\{\psi-\widetilde{\psi}\}\\
\leqslant& \kappa_2\|r\|_{0,\infty,\Omega}^{\frac{1}{2}}h^2\|\psi\|_{2,\Omega}\normdg{\zeta}
+\kappa_2\|\bm{b}\|_{0,\infty,\Omega}
\nu^{-\frac{1}{2}}h^2\|\psi\|_{2,\Omega}\normdg{\zeta}\\
&+\kappa_1\kappa_2\sqrt{C_g\mathcal N_0}A(h)\|\psi\|_{2,\Omega}\normdg{\zeta},
\end{split}
\end{equation}
where the term $A(h)$ depends on $h$, and it is written as:
\begin{equation*}
A(h)=
\begin{cases}
\|\bm{b}\|_{0,\infty,\Omega}\nu^{-\frac{1}{2}}
h&\text{if diffusion dominates},\\
\|\bm{b}\|_{0,\infty,\Omega}^{\frac{1}{2}}
h^{\frac{1}{2}}&\text{if convection dominates}.
\end{cases}
\end{equation*}

If the diffusion dominates, i.e., the value of $\nu$ is not small, we can achieve the optimal error order $O(h^{k+1})$. However when the convection dominates, following \eqref{convection:dominate:assum} about the relationship between $h$ and $\nu$, and Lemma \ref{auxiliary:regularity}, we can just obtain a sub-optimal error order $O(h^{k+\frac{1}{2}})$.

Hence, combining \eqref{var:stable}, \eqref{L2}--\eqref{ahrc}, Theorem \ref{priori:dg} and Lemma \ref{auxiliary:regularity}, we then obtain the error estimates \eqref{L2result} in the $L^2$-norm.
\end{pf}
\section*{References}
\bibliography{mybibfile}

\end{document}